\documentclass[12pt,onecolumn,draftclsnofoot]{IEEEtran}

\usepackage[T1]{fontenc}
\usepackage[utf8]{inputenc}
\usepackage{lmodern}
\usepackage{amsmath,amssymb,amsthm,bm,mathtools}
\usepackage{microtype}
\usepackage{enumitem}
\usepackage{xcolor}
\usepackage{hyperref}
\hypersetup{hidelinks}
\usepackage{bbm}
\usepackage{float}
\usepackage[section]{placeins}
\usepackage{perpage}
\MakePerPage{footnote}
\usepackage{cite}

\usepackage{tikz}
\usetikzlibrary{arrows.meta,positioning,calc,fit,backgrounds,shapes.geometric}

\newcommand{\E}{\mathbb{E}}

\newcommand{\Ld}{\Lambda_{\mathrm{death}}}
\newcommand{\1}{\mathbbm{1}}

\newcommand{\Prob}{\mathbb{P}}

\definecolor{bleudefrance}{rgb}{0.19, 0.55, 0.91}

\theoremstyle{definition}
\newtheorem{theorem}{Theorem}
\newtheorem{lemma}{Lemma}

\newtheorem{corollary}{Corollary}
\newtheorem{remark}{Remark}

\title{Fundamental Limits of \\ Self-Regulating Random Walks}

\author{Ali~Khalesi and Rawad~Bitar%
\thanks{A.~Khalesi is an Assistant Professor at Institut Polytechnique des Sciences Avanc\'ees (IPSA) and a member of LINCS Lab, Paris, France (e-mail: ali.khalesi@ipsa.fr). R.~Bitar is a Research Group Leader in the Chair of Communications Engineering, Technical University of Munich, Germany (e-mail: rawad.bitar@tum.de).}%
\thanks{Parts of the results were presented at the 2026 IEEE International Symposium on Information Theory~\cite{khalesi2026fundamental}.}
}

\begin{document}
\maketitle
\vspace{-50pt}
\begin{abstract}
We study self-regulating random walks (SRRWs), a decentralized mechanism for
maintaining a stable population of mobile tokens on a graph. Tokens move by
random walks, may be lost at faulty or malicious trap nodes, and may be locally
forked or terminated according to node visit ages. The goal is to prevent both
extinction and runaway growth without a central coordinator. We develop a
graph-aware, controller-agnostic theory for SRRWs on finite connected graphs
with lazy reversible random-walk dynamics. The main technical tool is a set of
node-dependent return-time tail bounds, which yield population-dependent
envelope functions describing how often local ages are large enough to trigger
control actions. These envelopes translate graph structure into drift
certificates: when the population is small, the achieved fork rate must
compensate trap-induced absorption; when the population is large, trap losses
and deliberate terminations must dominate forks. Using these certificates, we
characterize population corridors for which the full system-state process is
positive recurrent. We also identify two fundamental limits: stable
finite-communication operation requires sufficient graph-limited forking
capacity to offset absorption, and recovery from burst deletions or insertions
is limited by the graph-dependent frequency of eligible local visits. Finally,
we describe a minimal hysteresis controller with two age thresholds, one for
replenishment and one for suppression, and show that its recovery behavior
follows from the same certified margin conditions. Worked examples and
numerical experiments on heterogeneous graphs illustrate how the stationary
distribution, return-time envelopes, drift margins, finite-cost balance,
reaction-time limits, and operating corridors can be computed and validated in
concrete SRRW systems.
\end{abstract}

\section{Introduction}

Modern decentralized systems increasingly rely on lightweight \emph{tokens} that
move across a network and trigger local actions as they move. A token may
represent a model update in decentralized learning, a gradient carrier, a task
thread in a distributed computation, a relay message in a communication
network, or a data/task item carried by a decentralized drone swarm. In all
these settings, the network is represented by a graph: nodes may be processors,
servers, agents, waypoints, or spatial regions, and edges represent feasible
communication, computation, or mobility links. The key feature is that decisions
are local. A node does not know the full global state of the system; it only
reacts to information available when a token visits it.

A recent line of work has proposed \emph{self-regulating random walks} (SRRWs)
as a simple abstraction for such token-based systems~\cite{egger2024self}.
This abstraction is particularly relevant in the presence of trap nodes, captured
by the ``Pac-Man'' attack model introduced in~\cite{chen2025random}. In an SRRW,
tokens move as random walks on a graph. Tokens may be deleted when they visit
unreliable or malicious trap nodes. To compensate for losses and to prevent
excessive growth, ordinary nodes may also locally \emph{fork} a token, creating
one additional copy, or deliberately \emph{terminate} a token. The goal is to
maintain a stable population of tokens without a central coordinator
\cite{egger2024srrw,ayache2019rwg,ayache2021pwrwsgd,ayache2022walk,liu2024entrapment,khalesi2025}.
This type of mechanism is relevant, for example, in decentralized learning and
optimization over graphs, where mobile tokens repeatedly sample data or
gradients at visited nodes and collectively shape the learning dynamics
\cite{mao2020walkman,dimakis2010gossip,boyd2006randomized}. It is also natural
for decentralized drone swarms, where tokens can represent sensing tasks, map
tiles, detections, relay messages, or local model updates carried across a
mobility graph whose nodes are waypoints or regions of the operating area.

From a system-design perspective, SRRWs are attractive because they require very
little state at each node. Each node maintains a local \emph{age}, namely the
time elapsed since that node was last visited by any token. When a token arrives,
the node reads this age and chooses one of three actions:
\(\textsf{pass}\), \(\textsf{fork}\), or \(\textsf{terminate}\).
The action \(\textsf{pass}\) simply lets the token continue its random walk. The
action \(\textsf{fork}\) creates exactly one additional copy of the token. The
action \(\textsf{terminate}\) deliberately removes the token in order to reduce
the population. Trap-induced deletions are different: they are exogenous losses
caused by designated trap nodes. Thus the model separates two effects: losses
created by the environment, and corrective actions chosen by the local
controller.

The basic control problem is easy to state but difficult to analyze. If too few
tokens remain, the system may go extinct. If too many tokens are created, the
population may grow without control and overload the network. A good local
policy must therefore create enough tokens when the population is small and
suppress growth when the population is large. The difficulty is that each node
only observes local information, mainly its own visit age, while the global
population evolves through many interacting random walks, forks, trap deletions,
and deliberate terminations.

This motivates the main question of the paper: \emph{What graph-dependent limits
govern decentralized self-regulation?} We address this question from a
\emph{controller-agnostic} perspective. That is, we do not begin by fixing one
particular heuristic controller. Instead, we ask what any decentralized
age-based policy can achieve on a given graph, under a given trap configuration
and a given cap on the probability of forking at a visit.

The central quantity in our analysis is the \emph{population drift}. Informally,
the drift is the expected direction in which the population moves. At a
population level \(z\), the relevant per-visit drift is
\(p_{\mathrm{fork}}(z)-\Ld-K_{\mathrm{terminate}}(z)\). Here
\(p_{\mathrm{fork}}(z)\) is the actual mean probability that a token visit
creates a new token, \(K_{\mathrm{terminate}}(z)\) is the actual mean probability
that a non-trap visit deliberately terminates a token, and \(\Ld\) is the
per-token absorption pressure caused by traps. We call these quantities
\emph{achieved rates} because they are the rates actually produced by the local
policy at population level \(z\), after accounting for the graph, the random
walk, the age process, and the policy itself. Thus, an achieved-rate drift
condition simply means a condition on the actual net creation rate of tokens:
positive drift means that creation dominates deletion, while negative drift
means that deletion and suppression dominate creation.

Our analysis is built on the following intuition. A node can act only when it is
visited by a token. Moreover, an age-based controller can act only when the
node's local age is large enough. Therefore, the graph controls how often local
ages become large. If the random walk returns quickly to a node, that node's
age is frequently reset and large-age events are rare. If the random walk
returns slowly, large-age events are more common. This return-time behavior is
the bridge between the graph topology and the population-level control limits.

We formalize this bridge through graph-dependent return-time envelopes. For
each node \(u\), let \(T_u^+\) be the return time of the lazy random walk to
\(u\). We show that the tail probability \(\Pr\{T_u^+\ge A\}\) admits
node-dependent exponential bounds controlled by the stationary mass
\(\boldsymbol{\pi}(u)\). Since a node's age is reset by visits of any live token,
these return-time bounds lead to population-dependent non-trap envelope
functions \(\mathcal L_{\pi,z,\mathrm{nt}}^{+}(A)\) and
\(\mathcal L_{\pi,z,\mathrm{nt}}^{-}(A)\). These functions bound how often
non-trap nodes are old enough to trigger local actions at population level
\(z\). In simple terms, they quantify the amount of age-based control
opportunity available on the graph.

The envelopes then give graph-aware certificates for viability and safety.
Viability means that, when the population is too small, the achieved fork rate
is large enough to compensate trap-induced losses. Safety means that, when the
population is too large, trap losses and deliberate terminations dominate forks.
The same envelope functions also define margin functions that identify an
operating corridor \([Z_{\mathrm{low}},Z_{\mathrm{high}}]\), where the
population is intended to live. Below the corridor, the certified drift points
upward; above the corridor, the certified drift points downward.

In this paper, we make these ideas precise for finite connected graphs under
lazy reversible random-walk dynamics. The system model includes a trap set
\(\mathcal P_{\mathrm{trap}}\), trap-kill probabilities \(\zeta(u)\), a global
per-visit fork cap \(q\), local age-based control, and deliberate terminations.
The absorption pressure is
\(\Ld:=\sum_{u\in\mathcal P_{\mathrm{trap}}}\zeta(u)\boldsymbol{\pi}(u)\),
which is the stationary per-token rate at which traps delete tokens. This
quantity captures how strongly the trap set pulls the population toward
extinction.

\subsection*{Contributions}

We consider decentralized local age-based policies. When a token visits a
non-trap node, the node observes its local age, and possibly a coarse indication
of the population regime, and then chooses exactly one action: pass, fork once,
or terminate. Trap nodes delete visiting tokens according to their prescribed
trap probabilities. Forking is uniformly capped: at any eligible non-trap visit,
the fork probability is at most \(q\). Our main contributions are as follows.

\begin{itemize}
  \item \textbf{Graph-aware viability and safety envelopes:}
  We derive population-dependent envelope functions that bound how often
  non-trap nodes are eligible for age-based actions. These envelopes translate
  return-time behavior of the graph into population-level control limits. Using
  them, we give achieved-rate drift conditions for viability and safety. In
  particular, the actual fork rate must compensate the absorption pressure when
  the population is small, and the combined deletion and suppression rate must
  dominate forks when the population is large. When these inequalities hold
  with positive margins below and above a corridor, the full system-state
  process is positive recurrent to that corridor
  (Theorem~\ref{Theorem~1}).

  \item \textbf{Margin-based characterization of the operating corridor:}
  The envelopes define two transparent margin functions. The viability margin
  measures how much certified replenishment exceeds trap pressure below the
  corridor. The safety margin measures how much certified suppression exceeds
  forking above the corridor. These margins give a graph-aware way to choose
  \(Z_{\mathrm{low}}\) and \(Z_{\mathrm{high}}\): below
  \(Z_{\mathrm{low}}\) the drift is certified to be upward, while above
  \(Z_{\mathrm{high}}\) it is certified to be downward
  (Remark~\ref{rem:corridor-characterization}).

  \item \textbf{Finite-communication feasibility:}
  Each live token makes one random-walk transition per time step. Therefore the
  communication cost is the population size. We show that any viable and safe
  finite-communication policy must have enough graph-limited forking capability
  to balance the absorption pressure generated by traps. This gives a necessary
  feasibility condition in terms of the population-weighted upper fork envelope
  (Theorem~\ref{thm:comm-cost}).

  \item \textbf{Recovery limits after burst disturbances:}
  We study external deletion and insertion shocks. A deletion shock suddenly
  pushes the population below the corridor, while an insertion shock pushes it
  above the corridor. Because corrective actions can occur only through
  age-eligible token visits, the graph limits how quickly the population can
  return to the corridor. We show that the recovery time scales at least logarithmically in the
  shock ratio, with constants controlled by the graph-dependent return-time
  behavior (Theorem~\ref{thm:overshoot}).

  \item \textbf{Minimal hysteresis control:}
  We introduce a simple two-threshold policy, called minimal hysteresis control
  (MHC). Below the corridor, MHC uses an aggressive age threshold
  \(A_{\downarrow}\) to replenish the population. Above the corridor, it uses a
  conservative age threshold \(A_{\uparrow}\) together with suppression or
  deliberate termination. Inside the corridor, it remains idle. Under the
  certified margin conditions, MHC is positive recurrent to the corridor and
  its recovery behavior matches the graph-induced reaction limitation up to
  constant factors
  (Corollary~\ref{cor:overshoot-strategy}).

  \item \textbf{Worked examples:}
  We include explicit calculations on two small graphs: a three-node path and
  an irregular five-node graph. These examples show how to compute the
  stationary distribution, trap pressure, return-time tails, age-eligibility
  quantities, drift margins, and the resulting operating corridor. They make the
  proof mechanism concrete and show how the abstract envelope conditions can be
  checked directly on finite graphs (Section~\ref{sec:worked-example}).

  \item \textbf{Numerical validation on heterogeneous graphs:}
  We complement the worked examples with numerical experiments on a
  heterogeneous modular graph with trap nodes. These experiments validate the
  drift-corridor mechanism of Theorem~\ref{Theorem~1}, the finite-cost balance
  condition of Theorem~\ref{thm:comm-cost}, the reaction-time limitation of
  Theorem~\ref{thm:overshoot}, and the hysteresis recovery behavior of
  Corollary~\ref{cor:overshoot-strategy}. The simulations show how the
  theoretical quantities translate into operating corridors, exposure balances,
  recovery lower bounds, and logarithmic hysteresis recovery scaling in a
  concrete non-regular network (Section~\ref{sec:numerical-validation}).
\end{itemize}

Overall, the paper gives a graph-dependent design principle for decentralized
population regulation. Once the random-walk graph, trap profile, and fork cap
are fixed, the envelope functions determine where replenishment is possible,
where suppression is possible, which corridors can be certified, and how fast
the system can recover after a shock. This provides a rigorous foundation for
using simple local age-based policies in failure-prone distributed systems,
including decentralized learning networks, distributed computation systems,
communication networks, and decentralized drone swarms.

Section~\ref{sec:model} presents the system model and problem formulation.
Section~\ref{sec:main-results} states the main results: viability and safety
envelopes, the corridor characterization, the finite-communication feasibility
condition, the recovery limitation, and the hysteresis recovery guarantee.
Section~\ref{sec:worked-example} gives worked examples on small graphs.
Section~\ref{sec:numerical-validation} provides numerical validation on a
heterogeneous modular graph. Section~\ref{sec:conclusion} concludes the paper
and discusses future directions. The appendices contain the proofs of the main
results and the supporting lemmas.

\paragraph*{Notation}
Scalars are plain (e.g., \(z,q,A\)); vectors are bold lowercase
(e.g., \(\boldsymbol{\pi}\)); matrices are bold uppercase
(e.g., \(\mathbf P,\mathbf P'\)); and sets are calligraphic
(e.g., \(\mathcal V,\mathcal E,\mathcal P_{\mathrm{trap}}\)). The graph is
\(\mathcal G=(\mathcal V,\mathcal E)\), and tokens move according to the lazy
random-walk kernel \(\mathbf P'=\epsilon\mathbf I+(1-\epsilon)\mathbf P\), with
stationary distribution \(\boldsymbol{\pi}\), where $\mathbf P$ is a transition matrix and $0<\epsilon<1$. Time is discrete, \(t\in\mathbb N\),
and \(Z_t\) denotes the number of live tokens. We use \(\E[\cdot]\) and
\(\Prob(\cdot)\) for expectation and probability, \(\1\{\cdot\}\) for
indicators, \(x_+=\max\{x,0\}\), and all logarithms are natural.

\section{System Model and Problem Formulation}
\label{sec:model}

We consider a finite connected graph on which tokens move as lazy random walks.
Time is discrete, with \(\mathbb N=\{0,1,2,\ldots\}\). We denote by
\((Z_t)_{t\in\mathbb N}\) the population process, where \(Z_t\) is the number of
live tokens in the network at time \(t\). The system is initialized with a
single live token, \(Z_0=1\), whose initial location is either fixed
arbitrarily or drawn from an arbitrary distribution on \(\mathcal V\). The
choice of this initial location affects only the transient phase: all
stationary rates and envelope quantities below are evaluated after a burn-in
period exceeding the mixing time of the underlying lazy walk. Before analyzing
traps and control, we first characterize this underlying random walk. Its
stationary distribution determines the long-run visit frequencies of nodes,
which in turn govern how often tokens encounter traps and how often local ages
are refreshed.
\subsection{Graph and Random Walks}

Let \(\mathcal{G}=(\mathcal{V},\mathcal{E})\) be a finite, connected,
undirected graph with \(|\mathcal{V}|=n\). Each token moves on
\(\mathcal G\) according to a lazy reversible random walk with transition matrix
\begin{equation}
  \mathbf{P}'
  =
  \epsilon \mathbf{I}+(1-\epsilon)\mathbf{P},
  \qquad 0<\epsilon<1,
  \label{eq:model-lazy-kernel}
\end{equation}
where \(\mathbf{P}\) is the transition matrix of an underlying reversible random
walk on \(\mathcal{G}\).\footnote{Throughout the paper, all random walks are
lazy unless stated otherwise.} Since \(\mathcal{G}\) is connected and the walk
is lazy, the Markov chain with transition matrix \(\mathbf P'\) is irreducible
and aperiodic. It therefore admits a unique stationary distribution
\(\boldsymbol{\pi}\in\Delta^{n-1}\) satisfying
\(\boldsymbol{\pi}^{\top}\mathbf{P}'=\boldsymbol{\pi}^{\top}\).
Here
\[
  \Delta^{n-1}
  :=
  \Big\{
  \mathbf{x}\in\mathbb{R}^{n}: x_u\ge0,\ 
  \sum_{u\in\mathcal V}x_u=1
  \Big\}
\]
is the probability simplex over \(\mathcal V\).

At time \(t\), the network contains \(Z_t\) live tokens. Between control
actions, each live token evolves according to \(\mathbf P'\). The population may
increase through forks and may decrease through trap-induced deletions or
deliberate terminations. The following lemma recalls the stationary distribution
and convergence properties of the random walk.

\begin{lemma}[Stationary distribution on an undirected graph]
\label{lem:stationary}
Let \(\mathcal{G}=(\mathcal{V},\mathcal{E})\) be a finite, connected,
undirected graph. Consider either:
\begin{enumerate}[label=(\roman*),leftmargin=*]
  \item \textbf{Simple random walk (SRW):}
  \(\mathbf{P}_{uv}=1/\deg(u)\) if \((u,v)\in\mathcal E\), and \(0\) otherwise;

  \item \textbf{Weighted reversible random walk:} edges have symmetric weights
  \(w_{uv}=w_{vu}\ge0\), and
  \[
    \mathbf{P}_{uv}=\frac{w_{uv}}{w_u},
    \qquad
    w_u:=\sum_{v:(u,v)\in\mathcal E}w_{uv}.
  \]
\end{enumerate}
Define
\[
  \boldsymbol{\pi}(u)=
  \begin{cases}
    \dfrac{\deg(u)}{2|\mathcal E|}, & \text{for (i)},\\[6pt]
    \dfrac{w_u}{\sum_{x\in\mathcal V}w_x}, & \text{for (ii)}.
  \end{cases}
\]
Then \(\boldsymbol{\pi}\) is a stationary distribution for \(\mathbf P\):
\[
  \boldsymbol{\pi}^{\top}\mathbf{P}
  =
  \boldsymbol{\pi}^{\top},
  \qquad
  \sum_{u\in\mathcal V}\boldsymbol{\pi}(u)=1.
\]
Moreover, since \(\mathcal G\) is connected, the chain is irreducible. If it is
made lazy through \eqref{eq:model-lazy-kernel}, then it becomes aperiodic and
preserves the same stationary distribution. Consequently, for every initial
distribution \(\alpha\),
\[
  \lim_{t\to\infty}\alpha^{\top}(\mathbf{P}')^t
  =
  \boldsymbol{\pi}^{\top}.
\]
Finally, if \(Z_t\) independent tokens have locations distributed according to
\(\boldsymbol{\pi}\), then, conditioned on \(Z_t\), the node-occupancy vector
follows a multinomial distribution,
\(\mathrm{Multinomial}(Z_t,\boldsymbol{\pi})\).\footnote{The multinomial
distribution is the natural generalization of the binomial distribution to more
than two possible outcomes. Here, conditioned on \(Z_t\), each live token is
independently located at one node of \(\mathcal V\), and the probability that a
given token is at node \(u\) is \(\boldsymbol{\pi}(u)\). The node-occupancy
vector counts how many tokens are located at each node. Thus, if
\(\mathcal V=\{1,\ldots,n\}\), the vector
\((N_1(t),\ldots,N_n(t))\) satisfies
\(\sum_{u=1}^n N_u(t)=Z_t\), and for any nonnegative integers
\(m_1,\ldots,m_n\) with \(\sum_{u=1}^n m_u=Z_t\),
\[
\Pr\{N_1(t)=m_1,\ldots,N_n(t)=m_n\mid Z_t\}
=
\frac{Z_t!}{m_1!\cdots m_n!}
\prod_{u=1}^n \boldsymbol{\pi}(u)^{m_u}.
\]
In particular, \(\mathbb E[N_u(t)\mid Z_t]=Z_t\boldsymbol{\pi}(u)\).}

\end{lemma}

The proof of Lemma~\ref{lem:stationary} is given in
Appendix~\ref{Proof~of~Lemma~1}. The stationary distribution
\(\boldsymbol{\pi}\) links graph topology to system-level rates: it describes
the long-run frequency with which tokens visit each node. These visit
frequencies determine how often tokens encounter trap nodes and how often local
ages are refreshed by the token population.

The laziness in \eqref{eq:model-lazy-kernel} is used to avoid periodicity. On
time scales larger than the mixing time, token locations are approximately
\(\boldsymbol{\pi}\)-distributed. More precisely, if \(\gamma\) denotes the
spectral gap of \(\mathbf P'\), then the mixing time satisfies the standard
bound
\[
T_{\mathrm{mix}}(\varepsilon)
\le
\frac{1}{\gamma}
\log\!\Big(\frac{1}{\varepsilon\,\boldsymbol{\pi}_{\min}}\Big),
\qquad
\boldsymbol{\pi}_{\min}:=\min_{u\in\mathcal V}\boldsymbol{\pi}(u).
\]
In particular, when \(\boldsymbol{\pi}_{\min}\asymp 1/n\), this gives
\(T_{\mathrm{mix}}(\varepsilon)=O(\gamma^{-1}\log n)\) for fixed
\(\varepsilon\).\footnote{For a reversible Markov chain with transition matrix
\(\mathbf P'\), let \(\rho_2\) denote the second-largest eigenvalue of
\(\mathbf P'\) in magnitude, and define the spectral gap
\(\gamma:=1-\rho_2\); see, e.g., \cite[Ch.~12]{LevinPeres2017}.}
In the proofs, this stationary approximation is made precise through a
block-stationary comparison process after a burn-in period of length at least
the mixing time.

\subsection{Traps and Absorption Pressure}
A subset \(\mathcal P_{\mathrm{trap}}\subseteq\mathcal V\) contains unreliable
nodes. Whenever a token visits a trap node
\(u\in\mathcal P_{\mathrm{trap}}\), it is terminated with probability
\(\zeta(u)\in[0,1]\), independently across visits. Once the underlying random
walk has mixed, a typical token location is approximately distributed according
to \(\boldsymbol{\pi}\). In this regime, the expected per-token, per-step
trap-deletion probability is
\begin{equation}
  \Ld
  :=
  \sum_{u\in\mathcal P_{\mathrm{trap}}}
  \zeta(u)\boldsymbol{\pi}(u).
  \label{eq:model-Ldeath}
\end{equation}
We call \(\Ld\) the \emph{absorption pressure}. Equivalently, under stationarity,
a token experiences trap deletion at rate \(\Ld\) per step.

If the only live token is deleted before any new token is created, then the
population becomes zero. In that case, \(Z_t=0\) means that no live tokens remain
in the network. Without an additional mechanism, this zero-population state is
absorbing: once the system reaches it, no token remains to move, fork, or restart
the dynamics. To avoid making the analysis depend on this trivial collapse
event, all recurrence statements below are interpreted in one of two standard
ways. Either we restrict attention to the non-extinct communicating class, which
excludes the absorbing zero-population state, or we use a regeneration
convention, meaning that whenever the population reaches zero, the system is
restarted by injecting a new token. The viability conditions below are designed
to prevent collapse in the low-population regime by ensuring that, when the
population is small, the expected creation rate of new tokens dominates the
expected deletion rate.

\subsection{Local Observations and Control}

Each node \(u\) maintains a \emph{local visit age}
\begin{equation}
  A_u(t)=t-L_u(t),
  \label{eq:model-local-age}
\end{equation}
where \(L_u(t)\) is the most recent time strictly before \(t\) at which node
\(u\) was visited by any token, with a suitable initialization at \(t=0\). Thus
\(A_u(t)\) is a local variable at node \(u\), but its distribution depends on
the current population size. When there are many live tokens, nodes are visited
more frequently and local ages tend to be smaller. When there are few live
tokens, nodes are refreshed less frequently and local ages tend to be larger.

For technical recurrence statements, we may use the saturated age
\begin{equation}
  \bar A_u(t):=\min\{A_u(t),A_{\max}\},
  \label{eq:model-saturated-age}
\end{equation}
where \(A_{\max}\) is chosen larger than all triggering ages used by the policy.

Upon a token arrival to node \(u\) at time \(t\), the node observes its local age
and selects a single action
\[
  a_u(t)\in
  \{\textsf{fork},\,\textsf{terminate},\,\textsf{pass}\}.
\]
In the local-control class considered here, the decision is made at the node
where the token arrives. The node uses \(A_u(t)\) to decide whether to fork,
terminate, or pass the token. Some controllers may also use a coarse estimate of
the current population regime. For example, the node may know, or estimate,
whether the population is below the target corridor, inside the corridor, or
above the corridor. This is an idealized but useful abstraction: it does not
require nodes to know the exact global state, only a coarse population-regime
signal.

Thus, the most general local rule considered here can be written as
\begin{equation}
  a_u(t)
  \sim
  \mu_u\bigl(\cdot\,\big|\,A_u(t),\mathsf{Reg}(Z_t)\bigr),
  \label{eq:model-local-policy}
\end{equation}
where \(\mathsf{Reg}(Z_t)\) is the population-regime indicator, for example
\begin{equation}
  \mathsf{Reg}(Z_t)
  \in
  \{\textsf{low},\textsf{in},\textsf{high}\}.
  \label{eq:model-regime-indicator}
\end{equation}
Here \(\textsf{low}\) means \(Z_t<Z_{\mathrm{low}}\),
\(\textsf{in}\) means \(Z_{\mathrm{low}}\le Z_t\le Z_{\mathrm{high}}\), and
\(\textsf{high}\) means \(Z_t>Z_{\mathrm{high}}\). When no population-regime
information is used, the rule reduces to
\[
  a_u(t)\sim\mu_u(\cdot\,|\,A_u(t)).
\]
For notational simplicity, we often suppress the optional population-regime
argument and write only the dependence on the local age.

\paragraph*{Weak Average--Crossing (W--AC) rule}
For each non-trap node
\(u\in\mathcal V\setminus\mathcal P_{\mathrm{trap}}\), fix an age trigger
\(\theta_u\ge0\) and a per-visit fork-cap function
\(q_u:\mathbb R_+\to[0,1]\). The threshold \(\theta_u\) specifies when forking
is permitted, while \(q_u(A_u(t))\) specifies the probability that an eligible
visit to node \(u\) at time \(t\) creates one additional token.

\medskip
\noindent\textbf{Fork semantics:}
On a visit to a non-trap node \(u\) at time \(t\), the token either does not fork
or forks once, creating exactly one additional copy. Thus the population
increases by \(+1\) when a fork occurs. The quantity \(q_u(\cdot)\) is therefore
a probability cap, not a multiplicity factor. We define the global per-visit
fork cap
\begin{equation}
  q
  :=
  \sup_{u\in\mathcal V\setminus\mathcal P_{\mathrm{trap}},\,A\ge0}
  q_u(A)
  \in(0,1].
  \label{eq:model-global-fork-cap}
\end{equation}

Under W--AC, a visit to a non-trap node \(u\) at time \(t\) triggers a fork with
probability
\begin{equation}
  p_{\mathrm{fork}}(u\,|\,A_u(t))
  =
  q_u(A_u(t))\mathbf 1\{A_u(t)\ge\theta_u\}.
  \label{eq:model-WAC-fork-prob}
\end{equation}
The rule is fully local, since its decision depends only on the observed age at the visited node and on local randomization.

The graph dependence enters when we evaluate the induced fork intensity after
mixing. Because local age is reset by visits of any token, the probability that
a node remains old depends on the current population level \(z\). We define the achieved mean fork probability per token-visit at population level
\(z\) as
\begin{equation}
  p_{\mathrm{fork}}(z)
  :=
  \E[F_t^{\mathrm{fork}}\mid Z_t = \lceil z \rfloor],
  \label{eq:model-pfork}
\end{equation}
where \(F_t^{\mathrm{fork}}=1\) if a typical token-step triggers a fork and
\(F_t^{\mathrm{fork}}=0\) otherwise.\footnote{The notation
\(Z_t = \lceil z \rfloor\) indicates evaluation in a population regime concentrated near
the nominal level \(z\). Since \(Z_t\) is integer-valued, this may be interpreted
as conditioning on \(Z_t=z\) when \(z\in\mathbb N\), or on a small integer
neighborhood of \(z\) in the local block-stationary approximation.} For a common threshold
\(\theta_u=\theta\) and \(q_u(\cdot)\equiv q\), this becomes
\[
  p_{\mathrm{fork}}(z,\theta)
  :=
  q
  \sum_{u\in\mathcal V\setminus\mathcal P_{\mathrm{trap}}}
  \boldsymbol{\pi}(u)
  \Pr\{A_u(t)\ge\theta\mid U_t=u,\ Z_t= \lceil z \rfloor\},
\]
where \(U_t\) denotes the node visited by a typical token-step.

\paragraph*{Effective triggering ages}
For a local W--AC policy, the relevant age scale is summarized by effective
triggering ages. In the general case, these ages must be distinguished because
a policy may use node-dependent thresholds and may fork eligible visits with
probability strictly below the global cap \(q\).

For a population regime under consideration, define the true non-trap
age-eligibility probability at population level \(z\) and threshold \(A\) by
\begin{equation}
  R_{\pi,z,\mathrm{nt}}(A)
  :=
  \sum_{u\in\mathcal V\setminus\mathcal P_{\mathrm{trap}}}
  \boldsymbol{\pi}(u)
  \Pr\{A_u(t)\ge A\mid U_t=u,\ Z_t= \lceil z \rfloor\}.
  \label{eq:model-true-eligibility}
\end{equation}
The population-dependent non-trap envelopes are
\begin{equation}
  \mathcal L_{\pi,z,\mathrm{nt}}^{\pm}(A)
  :=
  \sum_{u\in\mathcal V\setminus\mathcal P_{\mathrm{trap}}}
  \boldsymbol{\pi}(u)
  e^{-c_{\pm}(u)zA\boldsymbol{\pi}(u)},
  \label{eq:model-L-env}
\end{equation}
where \(0<c_-(u)\le c_+(u)\) are node-dependent constants. The return-time
analysis in Appendix~\ref{Proof-of-Theorem-1} shows that
\begin{equation}
  \mathcal L_{\pi,z,\mathrm{nt}}^{+}(A)
  \le
  R_{\pi,z,\mathrm{nt}}(A)
  \le
  \mathcal L_{\pi,z,\mathrm{nt}}^{-}(A).
  \label{eq:model-eligibility-envelope}
\end{equation}

The \emph{upper-bound effective triggering age}, denoted by
\(A_{\mathrm{eff}}^{\mathrm{ub}}\), is chosen so that every fork-eligible visit
under the policy is also counted by the event
\[
  \{A_u(t)\ge A_{\mathrm{eff}}^{\mathrm{ub}}\}.
\]
For threshold policies with node-dependent thresholds \(\{\theta_u\}\), a valid
choice is
\begin{equation}
  A_{\mathrm{eff}}^{\mathrm{ub}}
  :=
  \min_{u\in\mathcal V\setminus\mathcal P_{\mathrm{trap}}}
  \theta_u.
  \label{eq:model-Aeff-up}
\end{equation}
Indeed,
\[
  \mathbf 1\{A_u(t)\ge\theta_u\}
  \le
  \mathbf 1\{A_u(t)\ge A_{\mathrm{eff}}^{\mathrm{ub}}\}.
\]
Since every per-visit fork probability is at most \(q\), the achieved fork
intensity satisfies the universal upper-envelope bound
\begin{equation}
  p_{\mathrm{fork}}(z)
  \le
  q\,R_{\pi,z,\mathrm{nt}}(A_{\mathrm{eff}}^{\mathrm{ub}})
  \le
  q\,\mathcal L_{\pi,z,\mathrm{nt}}^{-}
  (A_{\mathrm{eff}}^{\mathrm{ub}}).
  \label{eq:model-fork-upper-envelope}
\end{equation}

A corresponding lower-envelope bound is not automatic for an arbitrary W--AC
policy, because a policy may fork eligible visits with probability below the
cap \(q\). It is used only for policies that actually realize the
lower-envelope fork intensity. For example, if a threshold policy forks each
eligible non-trap visit with probability \(q\), then a valid lower-bound effective
triggering age is
\begin{equation}
  A_{\mathrm{eff}}^{\mathrm{lb}}
  :=
  \max_{u\in\mathcal V\setminus\mathcal P_{\mathrm{trap}}}
  \theta_u.
  \label{eq:model-Aeff-low}
\end{equation}
Indeed,
\[
  \mathbf 1\{A_u(t)\ge A_{\mathrm{eff}}^{\mathrm{lb}}\}
  \le
  \mathbf 1\{A_u(t)\ge\theta_u\}.
\]
Therefore, for such lower-envelope-realizing policies,
\begin{equation}
  p_{\mathrm{fork}}(z)
  \ge
  q\,R_{\pi,z,\mathrm{nt}}(A_{\mathrm{eff}}^{\mathrm{lb}})
  \ge
  q\,\mathcal L_{\pi,z,\mathrm{nt}}^{+}
  (A_{\mathrm{eff}}^{\mathrm{lb}}).
  \label{eq:model-fork-lower-envelope}
\end{equation}

When the policy uses a common threshold \(\theta_u=\theta\) over all non-trap
nodes and forks eligible visits at the cap \(q\), the two effective ages
coincide:
\[
  A_{\mathrm{eff}}^{\mathrm{ub}}
  =
  A_{\mathrm{eff}}^{\mathrm{lb}}
  =
  \theta.
\]
Only in this common-threshold, cap-realizing case do we write simply
\(A_{\mathrm{eff}}\). Then
\begin{equation}
  q\,\mathcal L_{\pi,z,\mathrm{nt}}^{+}(A_{\mathrm{eff}})
  \le
  p_{\mathrm{fork}}(z)
  \le
  q\,\mathcal L_{\pi,z,\mathrm{nt}}^{-}(A_{\mathrm{eff}}).
  \label{eq:model-fork-envelope}
\end{equation}
Thus, throughout the results and proofs, upper-envelope statements use
\(A_{\mathrm{eff}}^{\mathrm{ub}}\), lower-envelope viability certificates use
\(A_{\mathrm{eff}}^{\mathrm{lb}}\), and the shorter notation
\(A_{\mathrm{eff}}\) is used only when these two ages coincide. For hysteresis
policies, the same convention applies separately in the low- and
high-population regimes.

\paragraph*{Intentional terminations and the global termination rate}
A policy may also permit nodes to deliberately terminate a visiting token in
order to curb excessive population growth. Let
\(q_u^{\mathrm{terminate}}(a)\in[0,1]\) denote the probability that non-trap node
\(u\) intentionally terminates a token upon a visit when the observed age equals
\(a\). After mixing, the achieved mean intentional termination probability per
token-visit at population level \(z\) is
\begin{equation}
   K_{\mathrm{terminate}}(z)
   :=
   \mathbb{E}\!\left[
   q_{U_t}^{\mathrm{terminate}}\!\big(A_{U_t}(t)\big)
   \,\middle|\, Z_t= \lceil z \rfloor
   \right].
   \label{eq:model-Kterminate}
\end{equation}
Equivalently,
\[
   K_{\mathrm{terminate}}(z)
   \approx
   \sum_{u\in\mathcal V\setminus\mathcal P_{\mathrm{trap}}}
   \boldsymbol{\pi}(u)
   \mathbb{E}\!\left[
   q_u^{\mathrm{terminate}}\!\big(A_u(t)\big)
   \,\middle|\, U_t=u,\ Z_t= \lceil z \rfloor
   \right].
\]
Thus, \(K_{\mathrm{terminate}}(z)\) is distinct from trap-induced deletion: it
is the average probability that a non-trap visit deliberately removes a token.

\subsection{Population Dynamics}

Let \(\Phi_{\mathrm{fork}}(t)\), \(\Phi_{\mathrm{terminate}}(t)\), and
\(\Phi_{\mathrm{death}}(t)\) denote respectively the total numbers of intentional
forks, deliberate terminations, and trap-induced deletions at time \(t\). Then
\begin{equation}
  Z_{t+1}
  =
  Z_t
+
  \Phi_{\mathrm{fork}}(t)
  -
  \Phi_{\mathrm{death}}(t)
  -
  \Phi_{\mathrm{terminate}}(t).
  \label{eq:model-popdyn}
\end{equation}
Consequently,
\begin{equation}
  \E[Z_{t+1}-Z_t\mid Z_t]
  =
   \E[\Phi_{\mathrm{fork}}(t)\mid Z_t]
  -\E[\Phi_{\mathrm{death}}(t)\mid Z_t]
  -\E[\Phi_{\mathrm{terminate}}(t)\mid Z_t].
  \label{eq:model-popdyn-expectation}
\end{equation}

Under the one-action-per-visit rule, each non-trap visit can trigger at most one
of the three actions: fork, terminate, or pass. Different tokens may trigger
different actions at different nodes within the same time step. Trap-induced
deletions are accounted for separately through \(\Phi_{\mathrm{death}}(t)\).
If token locations at time \(t\) have empirical occupancy distribution
\(\alpha_t\) over \(\mathcal V\), then
\[
  \E[\Phi_{\mathrm{death}}(t)\mid Z_t]
  =
  Z_t
  \sum_{u\in\mathcal P_{\mathrm{trap}}}
  \zeta(u)\alpha_t(u).
\]
 After a transient period exceeding \(T_{\mathrm{mix}}\), we may approximate
\(\alpha_t\approx\boldsymbol{\pi}\). Hence, when \(Z_t= \lceil z \rfloor\),
\(\E[\Phi_{\mathrm{death}}(t)\mid Z_t= \lceil z \rfloor]= \lceil z \rfloor\Ld\). Likewise,
\(\E[\Phi_{\mathrm{fork}}(t)\mid Z_t= \lceil z \rfloor]
= \lceil z \rfloor\,p_{\mathrm{fork}}(z)\) and
\(\E[\Phi_{\mathrm{terminate}}(t)\mid Z_t= \lceil z \rfloor]
= \lceil z \rfloor\,K_{\mathrm{terminate}}(z)\). Thus, the stationary drift
approximation is
\begin{equation}
  \E[Z_{t+1}-Z_t\mid Z_t= \lceil z \rfloor]
  \approx
  z\big(
  p_{\mathrm{fork}}(z)
  -
  \Ld
  -
  K_{\mathrm{terminate}}(z)
  \big).
  \label{eq:model-stationary-drift}
\end{equation}

\subsection{Control Objective}
Given \((\mathcal G,\mathbf P',\boldsymbol{\pi})\) and the trap profile
\(\{\zeta(u)\}_{u\in\mathcal P_{\mathrm{trap}}}\), we study decentralized local
policies \(\{\mu_u\}_{u\in\mathcal V}\) that start from one token and regulate
the population process \((Z_t)_{t\in\mathbb N}\). The objective is to avoid
extinction, avoid runaway growth, and repeatedly return to a prescribed corridor
\([Z_{\mathrm{low}},Z_{\mathrm{high}}]\). The endpoints of this corridor may be
chosen directly as design parameters, or characterized through the
population-dependent viability and safety margins introduced in
Theorem~\ref{Theorem~1} and Remark~\ref{rem:corridor-characterization}. We seek
conditions under which the resulting dynamics are viable and safe, and we
characterize the fundamental limits on steady-state communication cost and
shock-recovery performance within this policy class.

\section{Main Results}
\label{sec:main-results}

We now state the main results of the paper. The first result gives the basic
viability and safety conditions for decentralized age-based policies. It shows
when the population has positive drift below the desired corridor and negative
drift above it, and it expresses these conditions through graph-dependent
age-eligibility envelopes.
The second result gives a finite-cost feasibility converse. Since trap losses
and fork opportunities both scale with the number of live tokens in the
per-token trap model, the converse takes the form of a necessary balance
condition on the available fork intensity, rather than a direct lower bound on
the mean population size.
The third result studies burst disturbances. It shows that after a sudden
deletion or insertion of tokens, the mean recovery time is limited by how often
the graph allows eligible local corrective visits.
Finally, we give a recovery guarantee for a two-threshold hysteresis policy.
This policy applies replenishment below the corridor, suppression above the
corridor, and remains idle inside the corridor.

\subsection{Viability and Safety Envelopes}

We now state the basic drift criterion. All quantities below are defined in
Section~\ref{sec:model}: the absorption pressure \(\Ld\), the global fork cap
\(q\), the achieved fork and termination rates \(p_{\mathrm{fork}}(z)\) and
\(K_{\mathrm{terminate}}(z)\), the upper-bound and lower-bound effective triggering ages
\(A_{\mathrm{eff}}^{\mathrm{ub}}\) and \(A_{\mathrm{eff}}^{\mathrm{lb}}\), and
the non-trap envelope functions
\(\mathcal L_{\pi,z,\mathrm{nt}}^{\pm}(\cdot)\). For a given population level
\(z\), we write the achieved per-visit drift as
\[
  d(z)
  :=
  p_{\mathrm{fork}}(z)-\Ld-K_{\mathrm{terminate}}(z).
\]
The sign of \(d(z)\) determines the local mean direction of the population:
positive drift means that forks dominate losses, while negative drift means that
trap deaths and deliberate terminations dominate forks.

The theorem below separates three points. First, viability and safety are
fundamentally achieved-rate drift conditions. Second, the graph-dependent
envelopes give sufficient certificates for these achieved-rate conditions.
Third, if the achieved drift has enough margin to dominate the block-averaging
error, then the full system-state process is recurrent to the prescribed
population corridor.

\begin{theorem}[Viability and safety envelopes]
\label{Theorem~1}
Consider any decentralized W--AC policy from Section~\ref{sec:model}. Then the
following statements hold.

\begin{enumerate}[label=\textup{(\roman*)},leftmargin=*]

\item \textbf{Achieved-rate viability condition.}
In a low-population regime in which the policy disables deliberate
terminations, so that \(K_{\mathrm{terminate}}(z)=0\), recovery requires that
the achieved fork rate compensate for the trap-induced absorption pressure.
Thus, a necessary achieved-rate condition is
\begin{equation}
  p_{\mathrm{fork}}(z)\ge \Ld .
  \label{eq:thm1-V0}
\end{equation}
\item \textbf{Achieved-rate safety condition.}
In a high-population regime, safety requires the net per-visit drift to be
nonpositive. Equivalently,
\begin{equation}
  d(z)
  =
  p_{\mathrm{fork}}(z)-\Ld-K_{\mathrm{terminate}}(z)
  \le 0 .
  \label{eq:thm1-S0}
\end{equation}

\item \textbf{Universal upper-envelope bound.}
For any W--AC policy, the upper-bound effective age
\(A_{\mathrm{eff}}^{\mathrm{ub}}\) yields the universal fork-rate bound
\[
  p_{\mathrm{fork}}(z)
  \le
  q\,\mathcal L_{\pi,z,\mathrm{nt}}^{-}
  (A_{\mathrm{eff}}^{\mathrm{ub}}).
\]
Consequently, high-population safety is certified whenever
\begin{equation}
  q\,\mathcal L_{\pi,z,\mathrm{nt}}^{-}
  (A_{\mathrm{eff}}^{\mathrm{ub}})
  \le   \Ld
  +
  K_{\mathrm{terminate}}(z).
  \tag{S}
  \label{eq:thm1-S-envelope}
\end{equation}

\item \textbf{Lower-envelope viability certificate.}
If the policy realizes the lower-envelope fork intensity with lower-bound
effective age \(A_{\mathrm{eff}}^{\mathrm{lb}}\), then low-population viability
is certified whenever
\begin{equation}
  q\,\mathcal L_{\pi,z,\mathrm{nt}}^{+}
  (A_{\mathrm{eff}}^{\mathrm{lb}})
  \ge
  \Ld .
  \tag{V}
  \label{eq:thm1-V-envelope}
\end{equation}
For a general W--AC policy that does not necessarily fork at the cap on eligible
visits, the low-population condition must instead be checked directly at the
achieved-rate level in \eqref{eq:thm1-V0}.

\item \textbf{Infeasibility obstructions.}
The same envelopes also give impossibility certificates. If
\[
  q\,\mathcal L_{\pi,z,\mathrm{nt}}^{-}
  (A_{\mathrm{eff}}^{\mathrm{ub}})
  <
  \Ld
\]
in a low-population regime, then no W--AC policy with fork cap \(q\) and that
upper-bound effective age can satisfy \eqref{eq:thm1-V0} at that population
level. Similarly, if \(d(z)>0\) in a high-population regime, then the safety
condition \eqref{eq:thm1-S0} fails. If the policy is known to realize a
lower-envelope fork intensity, then the stronger condition
\[
  q\,\mathcal L_{\pi,z,\mathrm{nt}}^{+}
  (A_{\mathrm{eff}}^{\mathrm{lb}})
  >  
  \Ld
  +
  K_{\mathrm{terminate}}(z)
\]
also certifies high-population drift infeasibility.

\item 
\textbf{Positive recurrence to a population corridor.}
Fix a mixing accuracy \(\varepsilon\in(0,1/8]\), and let
\(T_{\mathrm{mix}}=T_{\mathrm{mix}}(\varepsilon)\) be the corresponding
total-variation mixing time of the lazy walk. Define the largest effective
triggering age used in the block construction as
\(A_{\mathrm H}:=\max\{1,A_{\mathrm{eff}}^{\mathrm{ub}},
A_{\mathrm{eff}}^{\mathrm{lb}}\}\), with the convention
\(A_{\mathrm{eff}}^{\mathrm{lb}}
:=A_{\mathrm{eff}}^{\mathrm{ub}}\) when no lower-envelope certificate is
invoked. For a fixed constant \(\kappa\ge4\), choose the block length
\(B:=T_{\mathrm{mix}}+\kappa A_{\mathrm H}\). Thus, each block contains an
initial mixing period and an interval sufficiently long relative to every
effective triggering age used in the analysis. By Lemma~\ref{lem:block-Oz}, when a block starts with population \(z\), the
actual expected population drift over that block differs from the drift of the
frozen-population stationary comparison process by at most \(c_1(B)z\), where
\(c_1(B)<\infty\) is independent of \(z\). Let
\(\varepsilon_B:=c_1(B)/B\). At a high level, \(\varepsilon_B\) is the
normalized block-averaging error: it measures, per token and per time step, the
discrepancy caused by incomplete mixing and by population changes occurring
within the block. The drift margins below are required to exceed this
approximation error.

Suppose that there exist
\(0<Z_{\mathrm{low}}<Z_{\mathrm{high}}<\infty\) and
\(\eta_-,\eta_+>0\) such that
\begin{equation}
  p_{\mathrm{fork}}(z)-\Ld
  \ge
  \varepsilon_B+\eta_-,
  \qquad 0<z<Z_{\mathrm{low}},
  \label{eq:thm1-low-margin}
\end{equation}
and
\begin{equation}
  p_{\mathrm{fork}}(z)-\Ld-K_{\mathrm{terminate}}(z)
  \le
  -\varepsilon_B-\eta_+,
  \qquad z>Z_{\mathrm{high}} .
  \label{eq:thm1-high-margin}
\end{equation}
The first condition is implied by the lower-envelope certificate when
\(q\,\mathcal L_{\pi,z,\mathrm{nt}}^{+}
(A_{\mathrm{eff}}^{\mathrm{lb}})-\Ld
\ge\varepsilon_B+\eta_-\) for every
\(0<z<Z_{\mathrm{low}}\). The second condition is implied by the
upper-envelope certificate when
\(q\,\mathcal L_{\pi,z,\mathrm{nt}}^{-}
(A_{\mathrm{eff}}^{\mathrm{ub}})-\Ld-K_{\mathrm{terminate}}(z)
\le-\varepsilon_B-\eta_+\) for every \(z>Z_{\mathrm{high}}\).

Under finite age saturation, with saturation level larger than all triggering
ages used by the policy, and under the non-extinction/regeneration convention
of Section~\ref{sec:model}, the full system-state process is positive recurrent
to the corridor \([Z_{\mathrm{low}},Z_{\mathrm{high}}]\). Equivalently, the set
of full states whose population component lies in this corridor is a positive
recurrent petite set, and the process returns to it infinitely often with finite
mean return time.
\end{enumerate}
\end{theorem}

The proof is provided in Appendix~\ref{Proof-of-Theorem-1}.

Theorem~\ref{Theorem~1} can be used in two ways. If a specific controller is
already fixed, one may work directly with the achieved rates
\(p_{\mathrm{fork}}(z)\) and \(K_{\mathrm{terminate}}(z)\) and check
\eqref{eq:thm1-low-margin}--\eqref{eq:thm1-high-margin}. If one wants a
graph-level design rule before fully specifying the controller, the envelope
conditions \textup{(V)} and \textup{(S)} provide conservative certificates based
only on the graph, the trap profile, the fork cap, and the chosen triggering
ages.

This also gives a natural way to choose the operating corridor. The lower edge
of the corridor should be placed just above the population levels where
certified replenishment is still positive. The upper edge should be placed just
below the population levels where certified suppression becomes positive. The
following remark formalizes this construction.

\begin{remark}[Envelope-induced characterization of the operating corridor]
\label{rem:corridor-characterization}
The corridor endpoints \(Z_{\mathrm{low}}\) and \(Z_{\mathrm{high}}\) in
Theorem~\ref{Theorem~1} need not be arbitrary design parameters. They can be
chosen from the population-dependent envelope margins.

For a low-population replenishment age \(A_{\downarrow}\), define
\[
  M_{\mathrm V}(z;A_{\downarrow})
  :=
  q\,\mathcal L_{\pi,z,\mathrm{nt}}^{+}(A_{\downarrow})-\Ld .
\]
Thus \(M_{\mathrm V}(z;A_{\downarrow})>0\) means that the lower certified fork
intensity is larger than the trap-induced absorption pressure at population
level \(z\). This margin is relevant for policies that realize the lower
envelope, such as threshold policies that fork each eligible non-trap visit
with probability \(q\).

For a high-population suppression age \(A_{\uparrow}\), define
\[
  M_{\mathrm S}(z;A_{\uparrow})
  :=
  \Ld+K_{\mathrm{terminate}}(z)
  -
  q\,\mathcal L_{\pi,z,\mathrm{nt}}^{-}(A_{\uparrow}) .
\]
Thus \(M_{\mathrm S}(z;A_{\uparrow})>0\) means that trap-induced losses together
with deliberate terminations dominate the upper certified fork intensity.

Because local ages are reset by visits of any live token, the envelope functions
\(\mathcal L_{\pi,z,\mathrm{nt}}^{\pm}(A)\) decrease as the population \(z\)
increases. Hence \(M_{\mathrm V}(z;A_{\downarrow})\) is nonincreasing in \(z\).
If \(K_{\mathrm{terminate}}(z)\) is nondecreasing in \(z\), then
\(M_{\mathrm S}(z;A_{\uparrow})\) is nondecreasing in \(z\).

With the block correction \(\varepsilon_B\) from Theorem~\ref{Theorem~1}, and
for prescribed positive margins \(\eta_-\) and \(\eta_+\), define
\[
  Z_{\mathrm{low}}(\eta_-)
  :=
  \sup\big\{
  z\ge1:
  M_{\mathrm V}(z;A_{\downarrow})\ge \varepsilon_B+\eta_-
  \big\},
\]
and
\[
  Z_{\mathrm{high}}(\eta_+)
  :=
  \inf\big\{
  z\ge1:
  M_{\mathrm S}(z;A_{\uparrow})\ge \varepsilon_B+\eta_+
  \big\}.
\]
Whenever \(Z_{\mathrm{low}}(\eta_-)<Z_{\mathrm{high}}(\eta_+)\), the interval
\[
  \big[
  Z_{\mathrm{low}}(\eta_-),
  Z_{\mathrm{high}}(\eta_+)
  \big]
\]
is an envelope-certified operating corridor. Below this interval, the
block-corrected viability margin is positive; above it, the block-corrected
safety margin is positive.

For common-threshold policies, the low- and high-population ages may coincide,
and one may write \(A_{\downarrow}=A_{\uparrow}=A_{\mathrm{eff}}\). For
hysteresis policies such as MHC, the two ages are chosen separately:
\(A_{\downarrow}\) controls low-population replenishment, while
\(A_{\uparrow}\) controls high-population suppression. Since \(Z_t\) is
integer-valued, the practical corridor may be taken as the nearest integer
interval for which the two block-corrected margin inequalities remain valid.
\end{remark}

\subsection{Communication--Cost Feasibility Converse}
\label{subsec:comm-cost}

We now connect population regulation to communication cost. In our model, each
live token makes one random-walk transition per time step. Hence the
instantaneous communication cost is simply the number of live tokens $C_{\mathrm{comm}}(t)=Z_t$.
Therefore, the long-run average communication cost is the same as the long-run
mean population:
\begin{equation}
\label{Comm}
  \overline C_{\mathrm{comm}}
  :=
  \limsup_{T\to\infty}
  \frac1T\sum_{t=0}^{T-1}\E[C_{\mathrm{comm}}(t)]
  =
  \limsup_{T\to\infty}
  \frac1T\sum_{t=0}^{T-1}\E[Z_t]
  =:
  \overline Z .
\end{equation}

A finite-cost policy can be viable only if it creates enough new tokens to
replace, on average, the tokens lost at traps. The loss side is measured by the
long-run trap exposure per live token. In the stationary-envelope formulation,
this exposure is represented by the absorption pressure \(\Ld\). The creation
side is limited by the graph-dependent fork opportunities. In particular, the
upper fork envelope in \eqref{eq:model-fork-upper-envelope} bounds how often
the policy can create new tokens.

Both losses and fork opportunities scale with the number of live tokens.
Therefore the relevant balance is obtained after normalizing by the cumulative
population exposure \(\sum_{t=0}^{T-1}\E[Z_t]\). This gives the following
necessary feasibility condition.

\begin{theorem}[Finite-cost feasibility converse]
\label{thm:comm-cost}
Consider any decentralized W--AC policy from Section~\ref{sec:model} satisfying
the achieved-rate viability and safety conditions of
Theorem~\ref{Theorem~1}. Suppose that the policy has finite positive long-run
mean population, \(0<\overline Z<\infty\). Let
\(A_{\mathrm{eff}}^{\mathrm{ub}}\) be the upper-bound effective triggering age in
the operating regime. Then the following holds.

\begin{enumerate}[label=\textup{(\roman*)},leftmargin=*]

\item \textbf{Exposure-weighted upper fork envelope.}
The long-run fork opportunities available to the policy are bounded {from above}
by the exposure-weighted upper envelope
\begin{equation}
  \overline{\mathcal L}_{\mathrm{fork}}^{-}
  :=
  \limsup_{T\to\infty}
  \frac{
  \sum_{t=0}^{T-1}
  \E\!\left[
  Z_t\,
  \mathcal L_{\pi,Z_t,\mathrm{nt}}^{-}
  (A_{\mathrm{eff}}^{\mathrm{ub}})
  \right]
  }{
  \sum_{t=0}^{T-1}\E[Z_t]
  } .
  \label{eq:avg-upper-fork-envelope}
\end{equation}

\item \textbf{Actual trap exposure per live token.}
The corresponding long-run trap exposure per live token is
\[
  \widehat{\Lambda}_{\mathrm{death}}
  :=
  \liminf_{T\to\infty}
  \frac{
  \sum_{t=0}^{T-1}\E[\Phi_{\mathrm{death}}(t)]
  }{
  \sum_{t=0}^{T-1}\E[Z_t]
  } .
\]

\item \textbf{Necessary finite-cost balance.}
Every viable and safe finite-cost policy must satisfy
\[
  q\,\overline{\mathcal L}_{\mathrm{fork}}^{-}
  \ge
  \widehat{\Lambda}_{\mathrm{death}} .
\]
Equivalently, the long-run graph-limited fork intensity must be large enough to
cover the actual long-run trap exposure.

\item \textbf{Stationary-envelope form.}
In the stationary-envelope formulation used in Theorem~\ref{Theorem~1}, the
actual trap exposure is identified with the stationary absorption pressure
\(\Ld\), up to the block-stationary averaging error already included in the
proof. Hence, in this formulation,
\begin{equation}
  q\,\overline{\mathcal L}_{\mathrm{fork}}^{-}
  \ge
  \Ld .
  \label{eq:commcost-lower}
\end{equation}

\item \textbf{Concentrated operating regime.}
In particular, suppose that the steady-state population is concentrated around
a nominal operating level \(z_\star\), so that
\[
  \overline{\mathcal L}_{\mathrm{fork}}^{-}
  =
  \mathcal L_{\pi,z_\star,\mathrm{nt}}^{-}
  (A_{\mathrm{eff}}^{\mathrm{ub}})
  +o(1).
\]
Then, the feasibility condition reduces to
\[
  q\,\mathcal L_{\pi,z_\star,\mathrm{nt}}^{-}
  (A_{\mathrm{eff}}^{\mathrm{ub}})
  \ge
  \Ld
  \quad
  \text{up to the corresponding }o(1)\text{ error}.
\]
Thus, at the operating population level, the graph must provide enough
fork-eligible visits to balance the absorption pressure generated by the trap
set.

\end{enumerate}
\end{theorem}

The proof is provided in Appendix~\ref{proof-theorem2}.

\subsection{Reaction Under Burst Disturbances}
We next study how fast the population can recover after an external shock. The
setting and the relevant quantities are as follows.

\begin{enumerate}[label=\textup{(\roman*)},leftmargin=*]

\item \textbf{Shock model.}
A deletion shock at time \(t_0\) leaves the system with
\(Z_{t_0}=z_-<Z_{\mathrm{low}}\) tokens and an insertion shock leaves it with
\(Z_{t_0}=z_+>Z_{\mathrm{high}}\) tokens. The corridor
\([Z_{\mathrm{low}},Z_{\mathrm{high}}]\) is assumed to satisfy the viability
and safety conditions of Theorem~\ref{Theorem~1}.

\item \textbf{Local reaction constraint.}
Corrective actions can occur only when tokens visit nodes. A node can fork or
deliberately terminate a token only when a live token visits that node and the
observed local age satisfies the corresponding triggering rule. Moreover, under
the one-action-per-visit rule, a single visit can create at most one new token
or remove at most one token. Therefore, after a burst disturbance, the fastest
possible mean response is controlled by the same population-dependent
age-eligibility envelopes used in Theorem~\ref{Theorem~1}.

\item \textbf{Mean-response recovery time.}
All recovery times in this subsection are defined in terms of the mean
population \(\E[Z_t]\). That is, we measure when the expected population crosses
a prescribed threshold. We do not claim that every individual realization of
the random process reaches the threshold by that time.

\item \textbf{Maximal low-population growth rate.}
Let \(A_{\downarrow}\) be the triggering age used for replenishment below the
corridor. Define
\begin{equation}
  \Gamma_{+}^{\max}
  :=
  \sup_{1\le z<Z_{\mathrm{low}}}
  q\,\mathcal{L}_{\pi,z,\mathrm{nt}}^{-}(A_{\downarrow}).
  \label{eq:Gamma-plus-max}
\end{equation}
This is an optimistic upper bound on the per-token growth rate available during
low-population recovery.

\item \textbf{Maximal high-population removal rate.}
Let \(A_{\uparrow}\) be the triggering age used for suppression above the
corridor. Let \(\overline K_{\mathrm{terminate}}(z)\) be an upper envelope for
the deliberate termination rate \(K_{\mathrm{terminate}}(z)\). For example, if
deliberate terminations are age-gated at \(A_{\uparrow}\) with per-visit cap
\(r\), one may take
\[
  \overline K_{\mathrm{terminate}}(z)
  =
  r\,\mathcal{L}_{\pi,z,\mathrm{nt}}^{-}(A_{\uparrow}).
\]
Define
\begin{equation}
  \Gamma_{-}^{\max}
  :=
  \sup_{z>Z_{\mathrm{high}}}
  \bigl(
  \Ld+\overline K_{\mathrm{terminate}}(z)
  \bigr),
  \qquad
  0\le \Gamma_{-}^{\max}<1 .
  \label{eq:Gamma-minus-max}
\end{equation}
This is an optimistic upper bound on the per-token removal rate available during
high-population recovery.

\item \textbf{Graph-dependent recovery scale.}
The graph dependence can be summarized through the return-time exponent from
Lemma~\ref{lem:ret-tail}. Define
\[
  \lambda_r:=\min_{u\in\mathcal V}c_-(u).
\]
Then, for every node \(u\) and every \(A\ge1\),
\[
  \Pr\{T_u^+\ge A\}
  \le
  \exp\{-\lambda_r A\boldsymbol{\pi}(u)\}.
\]
When the chosen triggering ages make \(\Gamma_{+}^{\max}\) and
\(\Gamma_{-}^{\max}\) proportional to \(\lambda_r\), the recovery time has the
natural graph scale \(\lambda_r^{-1}\).

\item \textbf{Inner recovery band.}
Fix \(0<\rho<1/2\) and define the inner band
\[
  \mathcal C_\rho
  :=
  [Z_\rho^-,Z_\rho^+],
  \qquad
  Z_\rho^-:=(1-\rho)Z_{\mathrm{low}}+\rho Z_{\mathrm{high}},
  \qquad
  Z_\rho^+:=\rho Z_{\mathrm{low}}+(1-\rho)Z_{\mathrm{high}} .
\]
Thus, \(\mathcal C_\rho\) lies strictly inside the operating corridor.

\end{enumerate}

The quantities introduced above have the following interpretation. The constant
\(\Gamma_{+}^{\max}\) is the largest per-token growth rate that any
age-triggered replenishment rule can use below the corridor. Therefore, after a
deletion shock, even the fastest possible mean recovery cannot grow faster than
a multiplicative factor \(1+\Gamma_{+}^{\max}\) per step. Similarly,
\(\Gamma_{-}^{\max}\) is the largest per-token removal rate available above the
corridor. Hence, after an insertion shock, even the fastest possible mean
recovery cannot contract faster than a multiplicative factor
\(1-\Gamma_{-}^{\max}\) per step.

The next theorem turns these two one-step limitations into lower bounds on the
time needed for the mean population to re-enter the inner band
\(\mathcal C_\rho=[Z_\rho^-,Z_\rho^+]\). The result is therefore a reaction-time
converse: it does not depend on the details of a particular controller, but only
on the graph-limited availability of eligible corrective visits.

\begin{theorem}[Reaction limitation under burst disturbances]
\label{thm:overshoot}
Consider the SRRW dynamics of Section~\ref{sec:model} under a controller
satisfying the corridor conditions of Theorem~\ref{Theorem~1}. Then the
following mean-response recovery lower bounds hold.

\begin{enumerate}[label=\textup{(\roman*)},leftmargin=*]

\item \textbf{Deletion-shock recovery.}
Suppose that a deletion shock at time \(t_0\) leaves the system with
\(Z_{t_0}=z_-<Z_{\mathrm{low}}\). Define the mean recovery time to the lower
edge of the inner band by
\[
  \tau_{\mathrm{rec}}^{-}
  :=
  \inf\{t\ge t_0:\E[Z_t]\ge Z_\rho^{-}\}.
\]
Then
\begin{equation}
  \tau_{\mathrm{rec}}^{-}-t_0
  \ge
  \frac{
  \log\!\left(Z_\rho^{-}/z_-\right)
  }{
  \log(1+\Gamma_{+}^{\max})
  } .
  \label{eq:deletion-recovery-lower}
\end{equation}
Thus, after a deletion shock, the mean population cannot return to the inner
band faster than the maximal envelope-certified growth rate allows.

\item \textbf{Insertion-shock recovery.}
Suppose that an insertion shock at time \(t_0\) leaves the system with
\(Z_{t_0}=z_+>Z_{\mathrm{high}}\). Define the mean recovery time to the upper
edge of the inner band by
\[
  \tau_{\mathrm{rec}}^{+}
  :=
  \inf\{t\ge t_0:\E[Z_t]\le Z_\rho^{+}\}.
\]
Then
\begin{equation}
  \tau_{\mathrm{rec}}^{+}-t_0
  \ge
  \frac{
  \log\!\left(z_+/Z_\rho^{+}\right)
  }{
  -\log(1-\Gamma_{-}^{\max})
  } .
  \label{eq:insertion-recovery-lower}
\end{equation}
Thus, after an insertion shock, the mean population cannot return to the inner
band faster than the maximal envelope-certified removal rate allows.

\item \textbf{Graph-scale consequence.}
In particular, suppose that there exist constants \(C_+,C_->0\) such that
\[
  \Gamma_{+}^{\max}\le C_+\lambda_r,
  \qquad
  \Gamma_{-}^{\max}\le C_-\lambda_r,
  \qquad
  C_-\lambda_r<1 .
\]
Then the deletion-shock recovery time satisfies
\[
  \tau_{\mathrm{rec}}^{-}-t_0
  =
  \Omega\!\left(
  \frac1{\lambda_r}
  \log\frac{Z_\rho^{-}}{z_-}
  \right),
\]
and the insertion-shock recovery time satisfies
\[
  \tau_{\mathrm{rec}}^{+}-t_0
  =
  \Omega\!\left(
  \frac1{\lambda_r}
  \log\frac{z_+}{Z_\rho^{+}}
  \right).
\]
Therefore, the return to the corridor is limited by the graph-dependent
availability of eligible corrective visits.

\end{enumerate}
\end{theorem}

The proof is provided in Appendix~\ref{proof-of-theorem-3}.

We now state the corresponding recovery guarantee for a hysteresis controller.
The controller uses two feasible ages: \(A_{\downarrow}\) for
low-population replenishment and \(A_{\uparrow}\) for high-population
suppression. When \(Z_t<Z_{\mathrm{low}}\), it permits forks only at non-trap
visits with \(A_u(t)\ge A_{\downarrow}\). When
\(Z_t>Z_{\mathrm{high}}\), it suppresses aggressive forking and permits
deliberate terminations according to the high-population rule with triggering
age \(A_{\uparrow}\). Inside the corridor, it remains idle.

\begin{corollary}[Hysteresis recovery]
\label{cor:overshoot-strategy}
Consider the hysteresis controller described above. Let \(\varepsilon_B\) be the
normalized block-averaging error for the corresponding block construction.
Assume that there exist margins \(r_{\downarrow},r_{\uparrow}>0\) such that
\[
  p_{\mathrm{fork}}(z)-\Ld
  \ge
  \varepsilon_B+r_{\downarrow},
  \qquad
  z<Z_{\mathrm{low}},
\]
and
\[
  p_{\mathrm{fork}}(z)-\Ld-K_{\mathrm{terminate}}(z)
  \le
  -\varepsilon_B-r_{\uparrow},
  \qquad
  z>Z_{\mathrm{high}}.
\]
Then the following statements hold.

\begin{enumerate}[label=\textup{(\roman*)},leftmargin=*]

\item \textbf{Positive recurrence.}
The hysteresis controller is positive recurrent to the corridor
\([Z_{\mathrm{low}},Z_{\mathrm{high}}]\).

\item \textbf{Recovery after deletion shocks.}
After a deletion shock \(Z_{t_0}=z_-<Z_{\mathrm{low}}\), the mean recovery time
to the lower edge of the inner band satisfies
\[
  \tau_{\mathrm{rec}}^{-}-t_0
  =
  O\!\left(
  \frac1{r_{\downarrow}}
  \log\frac{Z_\rho^{-}}{z_-}
  \right),
\]
where
\[
  \tau_{\mathrm{rec}}^{-}
  :=
  \inf\{t\ge t_0:\E[Z_t]\ge Z_\rho^{-}\}.
\]

\item \textbf{Recovery after insertion shocks.}
After an insertion shock \(Z_{t_0}=z_+>Z_{\mathrm{high}}\), the mean recovery
time to the upper edge of the inner band satisfies
\[
  \tau_{\mathrm{rec}}^{+}-t_0
  =
  O\!\left(
  \frac1{r_{\uparrow}}
  \log\frac{z_+}{Z_\rho^{+}}
  \right),
\]
where
\[
  \tau_{\mathrm{rec}}^{+}
  :=
  \inf\{t\ge t_0:\E[Z_t]\le Z_\rho^{+}\}.
\]

\item \textbf{Matching the reaction limitation.}
When the drift margins \(r_{\downarrow}\) and \(r_{\uparrow}\) are comparable
to the envelope-certified response rates in Theorem~\ref{thm:overshoot}, the
hysteresis recovery time matches the reaction limitation up to constant factors.
In particular, when these rates are of order \(\lambda_r\), the recovery scale
is \(\lambda_r^{-1}\log(\cdot)\).

\end{enumerate}
\end{corollary}

The proof is provided in Appendix~\ref{Proof-of-corollary-1}.
\section{Worked Example}
\label{sec:worked-example}

This section gives two finite, hand-checkable illustrations of the main
certificates, showing how the stationary distribution, trap pressure,
return-time tails, age eligibility, and drift margins are computed to determine
the operating corridor and to instantiate Theorem~\ref{Theorem~1}, the
finite-communication converse of Theorem~\ref{thm:comm-cost}, the reaction
limitation of Theorem~\ref{thm:overshoot}, and the hysteresis recovery guarantee
of Corollary~\ref{cor:overshoot-strategy}.

\begin{figure}[t]
\centering
\begin{tikzpicture}[scale=1.0, every node/.style={circle,draw,minimum size=18pt,inner sep=1pt,font=\small}]

\node (p1) at (0,0) {1};
\node (p2) at (1.5,0) {2};
\node[fill=red!15] (p3) at (3,0) {3};
\draw (p1)--(p2)--(p3);
\node[draw=none, rectangle] at (1.5,-0.75) {\small Three-node path};
\node[draw=none, rectangle] at (3,-0.45) {\small trap};

\node (g1) at (5.5,0.6) {1};
\node (g2) at (5.5,-0.6) {2};
\node (g3) at (7.0,0) {3};
\node (g4) at (8.5,0) {4};
\node[fill=red!15] (g5) at (10.0,0) {5};

\draw (g1)--(g2)--(g3)--(g1);
\draw (g3)--(g4)--(g5);

\node[draw=none, rectangle] at (7.75,-1.0) {\small Irregular five-node graph};
\node[draw=none, rectangle] at (10.0,-0.45) {\small trap};

\end{tikzpicture}
\caption{Two illustrative graphs. Left: a three-node path with trap node \(3\).
Right: an irregular five-node graph with trap node \(5\).}
\label{fig:worked-graphs}
\end{figure}

In the general theory, local age is reset by visits of any live token. Hence the
relevant probabilities are population dependent. To keep the arithmetic
transparent, the examples first compute exact single-token return tails. Then
we use a population-aware surrogate obtained by raising the single-token
non-return probability to the power \(z\). This surrogate has the same
monotonicity as the envelope functions in Theorem~\ref{Theorem~1}: increasing
the population makes old local ages less likely. The block correction
\(\varepsilon_B\) is omitted from the displayed arithmetic; in the theorem, the
same strict margins must be chosen large enough to dominate this correction.

\subsection{Three-node path}

Consider the three-node path in Fig.~\ref{fig:worked-graphs}, with node \(3\)
as the unique trap. The token follows the lazy simple random walk with laziness
\(\epsilon=1/2\). With nodes ordered as \(1,2,3\),
\[
\mathbf{P}'
=
\frac12\mathbf I
+
\frac12
\begin{pmatrix}
0 & 1 & 0\\[2pt]
1/2 & 0 & 1/2\\[2pt]
0 & 1 & 0
\end{pmatrix}
=
\begin{pmatrix}
1/2 & 1/2 & 0\\[2pt]
1/4 & 1/2 & 1/4\\[2pt]
0 & 1/2 & 1/2
\end{pmatrix}.
\]

\subsubsection*{Step 1: Stationarity and absorption pressure}

The ordinary path has degrees \(1,2,1\). Therefore, by
Lemma~\ref{lem:stationary},
\[
   \boldsymbol{\pi}(1)=\frac14,\qquad
   \boldsymbol{\pi}(2)=\frac12,\qquad
   \boldsymbol{\pi}(3)=\frac14.
\]
This validates the first step of the appendix proof: once the lazy walk has
mixed, node \(u\) is visited with long-run frequency \(\boldsymbol{\pi}(u)\).
If the trap-kill probability at node \(3\) is \(\zeta\), then the stationary
trap pressure is
\[
   \Ld
   =
   \zeta\boldsymbol{\pi}(3)
   =
   \frac{\zeta}{4}.
\]
For the numerical checks below we take
\[
   \zeta=\frac12,
   \qquad
   \Ld=\frac18.
\]
The eigenvalues of \(\mathbf P'\) are \(1,1/2,0\), so the spectral gap is
\(\gamma=1/2\). Thus the mixing time is constant on this graph, and the
block-stationary comparison used in Appendix~\ref{Proof-of-Theorem-1} has a
fixed finite burn-in scale.

\subsubsection*{Step 2: Return tails and age eligibility}

The non-trap nodes are \(1\) and \(2\). We compute the relevant return tails by
conditioning on the first few steps of the walk.

For \(u=1\), the event \(\{T_1^+\ge2\}\) means that the walk does not return to
node \(1\) at the next step. This has probability \(1/2\). Moreover, to have
\(T_1^+\ge3\), the walk must move from \(1\) to \(2\), and then avoid returning
from \(2\) to \(1\). Hence
\[
   \Pr_1\{T_1^+\ge2\}=\frac12,
   \qquad
   \Pr_1\{T_1^+\ge3\}
   =
   \frac12\cdot\frac34
   =
   \frac38 .
\]
Similarly, for \(u=2\),
\[
   \Pr_2\{T_2^+\ge2\}=\frac12,
   \qquad
   \Pr_2\{T_2^+\ge3\}
   =
   \frac14 .
\]
For the high-population check, we also use
\[
   \Pr_1\{T_1^+\ge4\}=\frac{5}{16},
   \qquad
   \Pr_2\{T_2^+\ge4\}=\frac18 .
\]

For a threshold \(A\), define the population-aware non-trap eligibility
functional
\[
   \mathcal L_{\mathrm{path,nt}}^{(z)}(A)
   :=
   \sum_{u\in\{1,2\}}
   \boldsymbol{\pi}(u)\Pr_u\{T_u^+\ge A\}^{z}.
\]
This is the finite-graph analogue of the population-dependent envelopes in
Lemma~\ref{lem:population-age-envelope}. It decreases in \(z\), exactly as the
general proof requires, because more live tokens refresh local ages more often.

For replenishment, choose
\[
   A_{\downarrow}=2.
\]
Then
\[
   \mathcal L_{\mathrm{path,nt}}^{(z)}(2)
   =
   \frac14\Big(\frac12\Big)^z
   +
   \frac12\Big(\frac12\Big)^z
   =
   \frac34\Big(\frac12\Big)^z .
\]
For high-population suppression, choose
\[
   A_{\uparrow}=4.
\]
Then
\[
   \mathcal L_{\mathrm{path,nt}}^{(z)}(4)
   =
   \frac14\Big(\frac{5}{16}\Big)^z
   +
   \frac12\Big(\frac18\Big)^z .
\]

\subsubsection*{Validation of Theorem~\ref{Theorem~1}}

Take
\[
   q=\frac12,
   \qquad
   r=\frac14,
\]
where \(q\) is the fork cap and \(r\) is the deliberate-termination cap above
the corridor.

Below the corridor, deliberate termination is inactive. At \(z=1\),
\[
   q\,\mathcal L_{\mathrm{path,nt}}^{(1)}(2)-\Ld
   =
   \frac12\cdot\frac38-\frac18
   =
   \frac1{16}>0.
\]
At \(z=2\),
\[
   q\,\mathcal L_{\mathrm{path,nt}}^{(2)}(2)-\Ld
   =
   \frac12\cdot\frac{3}{16}-\frac18
   =
   -\frac1{32}<0.
\]
Thus the positive low-population drift is certified only at the strictly low
level \(z=1\). We therefore take
\[
   Z_{\mathrm{low}}=2.
\]

Above the corridor, the safety margin associated with \(A_{\uparrow}=4\) is
\[
  M_{\mathrm S,path}(z;4)
  :=
  \Ld+r\,\mathcal L_{\mathrm{path,nt}}^{(z)}(4)
  -
  q\,\mathcal L_{\mathrm{path,nt}}^{(z)}(4)
  =
  \frac18-\frac14\,\mathcal L_{\mathrm{path,nt}}^{(z)}(4).
\]
Since \(\mathcal L_{\mathrm{path,nt}}^{(z)}(4)\) decreases in \(z\),
\(M_{\mathrm S,path}(z;4)\) increases in \(z\). In particular,
\[
   M_{\mathrm S,path}(3;4)
   =
   \frac{8051}{65536},
   \qquad
   M_{\mathrm S,path}(4;4)
   =
   \frac{130415}{1048576}.
\]
Choosing a safety margin \(\eta_+\) with
\[
   M_{\mathrm S,path}(3;4)<\eta_+\le M_{\mathrm S,path}(4;4)
\]
makes \(4\) the first population level at which the high-population safety
margin is certified. We include this first certified level as the upper edge of
the integer corridor and take
\[
   Z_{\mathrm{high}}=4.
\]
Therefore the certified operating corridor is
\[
   [Z_{\mathrm{low}},Z_{\mathrm{high}}]=[2,4].
\]

This is precisely the drift structure used in Appendix~\ref{Proof-of-Theorem-1}.
Below the corridor, the block drift has the form
\[
   \E[Z_{k+1}-Z_k\mid Z_k=z,\ z<2]
   =
   zB\cdot\frac1{16}\pm c_1z.
\]
For sufficiently small normalized block error \(\varepsilon_B=c_1/B\), this is
positive. Above the corridor, the margin is at least
\(M_{\mathrm S,path}(5;4)\) for integer \(z>4\), so
\[
   \E[Z_{k+1}-Z_k\mid Z_k=z,\ z>4]
   \le
   -zB\,M_{\mathrm S,path}(5;4)\pm c_1z,
\]
which is negative once the strict margin dominates the block error. Since the
full state space is finite on bounded population sets after age saturation, the
same Foster--Lyapunov argument as in Theorem~\ref{Theorem~1} gives positive
recurrence to \([2,4]\), under the non-extinction/regeneration convention.

\subsubsection*{Validation of Theorem~\ref{thm:comm-cost}}

Theorem~\ref{thm:comm-cost} is a necessary long-run balance condition. It says
that any viable and safe finite-cost policy must satisfy
\[
   q\,\overline{\mathcal L}_{\mathrm{fork}}^{-}
   \ge
   \Ld .
\]
On the path, the low-population fork envelope at \(z=1\) gives
\[
   q\,\mathcal L_{\mathrm{path,nt}}^{(1)}(2)
   =
   \frac12\cdot\frac38
   =
   \frac{3}{16}
   >
   \frac18
   =
   \Ld .
\]
Thus, whenever the process is in the low state \(z=1\), the graph provides
enough fork opportunity to compensate trap pressure. At \(z=2\), however,
\[
   q\,\mathcal L_{\mathrm{path,nt}}^{(2)}(2)
   =
   \frac12\cdot\frac{3}{16}
   =
   \frac{3}{32}
   <
   \frac18
   =
   \Ld .
\]
This does not contradict the converse. It illustrates why
Theorem~\ref{thm:comm-cost} is stated in terms of the exposure-weighted average
\[
  \overline{\mathcal L}_{\mathrm{fork}}^{-}
  =
  \limsup_{T\to\infty}
  \frac{
  \sum_{t=0}^{T-1}
  \E\!\left[
  Z_t\mathcal L_{\pi,Z_t,\mathrm{nt}}^{-}(A_{\mathrm{eff}}^{\mathrm{ub}})
  \right]
  }{
  \sum_{t=0}^{T-1}\E[Z_t]
  }.
\]
The path example therefore validates the proof mechanism of
Theorem~\ref{thm:comm-cost}: long-run viability requires enough cumulative
exposure to states where the available fork envelope compensates the accumulated
trap deaths.

\subsubsection*{Validation of Theorem~\ref{thm:overshoot}}

For a deletion shock, the proof of Theorem~\ref{thm:overshoot} upper-bounds the
fastest possible mean growth before recovery by
\[
   \E[Z_{t+1}]
   \le
   (1+\Gamma_+^{\max})\E[Z_t].
\]
On the path, since \(Z_{\mathrm{low}}=2\), the only integer population below the
corridor is \(z=1\). Hence
\[
   \Gamma_{+,\mathrm{path}}^{\max}
   =
   q\,\mathcal L_{\mathrm{path,nt}}^{(1)}(2)
   =
   \frac{3}{16}.
\]
The inner band endpoints are
\[
   Z_\rho^- = (1-\rho)2+\rho4 = 2+2\rho,
   \qquad
   Z_\rho^+ = \rho2+(1-\rho)4 = 4-2\rho.
\]
Thus after a deletion shock to \(z_-<2\),
\[
   \tau_{\mathrm{rec}}^{-}-t_0
   \ge
   \frac{
   \log\!\left((2+2\rho)/z_-\right)
   }{
   \log(1+3/16)
   }.
\]

For an insertion shock, the proof uses the fastest possible contraction rate
\[
   \Gamma_{-}^{\max}
   =
   \sup_{z>Z_{\mathrm{high}}}
   \bigl(\Ld+\overline K_{\mathrm{terminate}}(z)\bigr).
\]
Here \(Z_{\mathrm{high}}=4\), and the integer levels above the corridor satisfy
\(z\ge5\). Since the eligibility term decreases in \(z\), the supremum is
attained at \(z=5\). Therefore
\[
   \Gamma_{-,\mathrm{path}}^{\max}
   =
   \frac18
   +
   \frac14\,\mathcal L_{\mathrm{path,nt}}^{(5)}(4)
   <
   1.
\]
The reaction lower bound becomes
\[
   \tau_{\mathrm{rec}}^{+}-t_0
   \ge
   \frac{
   \log\!\left(z_+/(4-2\rho)\right)
   }{
   -\log\!\left(1-\Gamma_{-,\mathrm{path}}^{\max}\right)
   }.
\]
This is exactly the multiplicative comparison used in the proof of
Theorem~\ref{thm:overshoot}.

\subsubsection*{Validation of Corollary~\ref{cor:overshoot-strategy}}

The hysteresis controller for the path uses \(A_{\downarrow}=2\) below the
corridor and \(A_{\uparrow}=4\) above the corridor. The low-population
block-corrected margin can be chosen as any
\[
   0<r_{\downarrow}
   <
   \frac1{16}-\varepsilon_B,
\]
provided \(\varepsilon_B<1/16\). Above the corridor, since \(z>4\) implies
\(z\ge5\), the high-population block-corrected margin can be chosen as any
\[
   0<r_{\uparrow}
   <
   M_{\mathrm S,path}(5;4)-\varepsilon_B.
\]
With these choices, the hypotheses of Corollary~\ref{cor:overshoot-strategy}
are satisfied. Hence the hysteresis controller is positive recurrent to
\([2,4]\), and its mean-response recovery times satisfy
\[
   \tau_{\mathrm{rec}}^{-}-t_0
   =
   O\!\left(
   \frac1{r_{\downarrow}}
   \log\frac{2+2\rho}{z_-}
   \right),
\]
and
\[
   \tau_{\mathrm{rec}}^{+}-t_0
   =
   O\!\left(
   \frac1{r_{\uparrow}}
   \log\frac{z_+}{4-2\rho}
   \right).
\]
When these margins are comparable to the envelope-certified rates above, the
hysteresis recovery time matches the lower bounds of
Theorem~\ref{thm:overshoot} up to constant factors.

\subsection{Irregular five-node graph}

We now repeat the same verification on a graph without the symmetry of the
path. Let
\[
   \mathcal V=\{1,2,3,4,5\},
   \qquad
   \mathcal E=\{(1,2),(1,3),(2,3),(3,4),(4,5)\},
\]
and let node \(5\) be the unique trap. Again use the lazy simple random walk
with \(\epsilon=1/2\). The lazy transition matrix is
\[
\mathbf P'_{\mathrm{irr}}
=
\begin{pmatrix}
1/2 & 1/4 & 1/4 & 0   & 0\\[2pt]
1/4 & 1/2 & 1/4 & 0   & 0\\[2pt]
1/6 & 1/6 & 1/2 & 1/6 & 0\\[2pt]
0   & 0   & 1/4 & 1/2 & 1/4\\[2pt]
0   & 0   & 0   & 1/2 & 1/2
\end{pmatrix}.
\]

\subsubsection*{Step 1: Stationarity and absorption pressure}

The degrees are
\[
   d_1=2,\qquad d_2=2,\qquad d_3=3,\qquad d_4=2,\qquad d_5=1.
\]
Therefore
\[
   \boldsymbol{\pi}
   =
   \left(
   \frac15,\frac15,\frac3{10},\frac15,\frac1{10}
   \right).
\]
Since node \(5\) is the trap,
\[
   \Ld
   =
   \zeta\boldsymbol{\pi}(5)
   =
   \frac{\zeta}{10}.
\]
For the numerical checks below we again take
\[
   \zeta=\frac12,
   \qquad
   \Ld=\frac1{20}.
\]
This validates the same stationary trap-exposure step used in the proofs of
Theorems~\ref{Theorem~1} and~\ref{thm:comm-cost}.

\subsubsection*{Step 2: Return tails and age eligibility}

For each node \(u\), let \(Q^{(u)}\) be the matrix obtained from
\(\mathbf P'_{\mathrm{irr}}\) by deleting row \(u\) and column \(u\). Then
\[
   \Pr_u\{T_u^+\ge1\}=1,
\]
and, for every integer \(m\ge2\),
\[
   \Pr_u\{T_u^+\ge m\}
   =
   \sum_{v\neq u}
   \mathbf P'_{\mathrm{irr}}(u,v)
   \big[(Q^{(u)})^{m-2}\mathbf 1\big]_v.
\]
This identity is the finite-state version of the return-time argument in
Lemma~\ref{lem:ret-tail}. It follows by conditioning on the first step away
from \(u\) and then requiring the killed chain on
\(\mathcal V\setminus\{u\}\) to avoid \(u\) for the remaining steps.

Since node \(5\) is the trap, define the non-trap return-tail functional
\[
   \mathcal L_{\mathrm{irr,nt}}(m)
   :=
   \sum_{u\in\{1,2,3,4\}}
   \boldsymbol{\pi}(u)\Pr_u\{T_u^+\ge m\}.
\]
Direct evaluation gives
\[
   \mathcal L_{\mathrm{irr,nt}}(1)=\frac9{10},\qquad
   \mathcal L_{\mathrm{irr,nt}}(2)=\frac9{20},\qquad
   \mathcal L_{\mathrm{irr,nt}}(3)=\frac{27}{80},
\]
and
\[
   \mathcal L_{\mathrm{irr,nt}}(7)=\frac{1589}{10240}.
\]

As before, define the population-aware surrogate
\[
   \mathcal L_{\mathrm{irr,nt}}^{(z)}(m)
   :=
   \sum_{u\in\{1,2,3,4\}}
   \boldsymbol{\pi}(u)\Pr_u\{T_u^+\ge m\}^{z}.
\]
This quantity is nonincreasing in \(z\), just like the true population-dependent
envelopes in Lemma~\ref{lem:population-age-envelope}.

\subsubsection*{Validation of Theorem~\ref{Theorem~1}}

Take
\[
   q=\frac12,
   \qquad
   r=\frac14,
   \qquad
   A_{\downarrow}=3,
   \qquad
   A_{\uparrow}=7.
\]

For the low-population check, the relevant return tails at \(m=3\) are
\[
  \Pr_1\{T_1^+\ge3\}
  =
  \Pr_2\{T_2^+\ge3\}
  =
  \frac{19}{48},
\]
and
\[
  \Pr_3\{T_3^+\ge3\}=\frac38,
  \qquad
  \Pr_4\{T_4^+\ge3\}=\frac13.
\]
Therefore
\[
   \mathcal L_{\mathrm{irr,nt}}^{(1)}(3)=\frac{27}{80},
   \qquad
   \mathcal L_{\mathrm{irr,nt}}^{(2)}(3)=\frac{61}{480}.
\]
The low-population margins are
\[
  q\,\mathcal L_{\mathrm{irr,nt}}^{(1)}(3)-\Ld
  =
  \frac12\cdot\frac{27}{80}-\frac1{20}
  =
  \frac{19}{160}>0,
\]
and
\[
  q\,\mathcal L_{\mathrm{irr,nt}}^{(2)}(3)-\Ld
  =
  \frac12\cdot\frac{61}{480}-\frac1{20}
  =
  \frac{13}{960}>0.
\]
Choosing
\[
  \frac{13}{960}<\eta_-\le\frac{19}{160}
\]
makes only \(z=1\) strictly below the certified corridor after imposing the
chosen margin. Hence
\[
   Z_{\mathrm{low}}=2.
\]

For the high-population check, the non-trap return tails at \(m=7\) are
\[
  \Pr_1\{T_1^+\ge7\}
  =
  \Pr_2\{T_2^+\ge7\}
  =
  \frac{173}{864},
\]
and
\[
  \Pr_3\{T_3^+\ge7\}
  =
  \frac{147}{1024},
  \qquad
  \Pr_4\{T_4^+\ge7\}
  =
  \frac{2213}{13824}.
\]
The high-population safety margin is
\[
  M_{\mathrm{S,irr}}(z;7)
  :=
  \Ld+r\,\mathcal L_{\mathrm{irr,nt}}^{(z)}(7)
  -
  q\,\mathcal L_{\mathrm{irr,nt}}^{(z)}(7)
  =
  \frac1{20}
  -
  \frac14\,\mathcal L_{\mathrm{irr,nt}}^{(z)}(7).
\]
Since \(\mathcal L_{\mathrm{irr,nt}}^{(z)}(7)\) decreases in \(z\),
\(M_{\mathrm{S,irr}}(z;7)\) increases in \(z\). Choosing
\[
  M_{\mathrm{S,irr}}(4;7)<\eta_+\le M_{\mathrm{S,irr}}(5;7)
\]
makes \(5\) the first population level at which the high-population safety
margin is certified. We include this first certified level as the upper edge of
the integer corridor:
\[
   Z_{\mathrm{high}}=5.
\]
Therefore
\[
   [Z_{\mathrm{low}},Z_{\mathrm{high}}]=[2,5].
\]

The block drift signs are explicit. Below the corridor,
\[
   p_{\mathrm{fork}}^{\downarrow}-\Ld
   =
   q\,\mathcal L_{\mathrm{irr,nt}}(3)-\Ld
   =
   \frac12\cdot\frac{27}{80}-\frac1{20}
   =
   \frac{19}{160}>0.
\]
Above the corridor,
\[
   p_{\mathrm{fork}}^{\uparrow}
   =
   \frac12\cdot\frac{1589}{10240}
   =
   \frac{1589}{20480},
\]
and
\[
   K_{\mathrm{terminate}}^{\uparrow}
   =
   \frac14\cdot\frac{1589}{10240}
   =
   \frac{1589}{40960}.
\]
Thus
\[
   p_{\mathrm{fork}}^{\uparrow}
   -
   \Ld
   -
   K_{\mathrm{terminate}}^{\uparrow}
   =
   -\frac{459}{40960}<0.
\]
Consequently,
\[
   \E[Z_{k+1}-Z_k\mid Z_k=z,\ z<2]
   =
   zB\cdot\frac{19}{160}\pm c_1z,
\]
and
\[
   \E[Z_{k+1}-Z_k\mid Z_k=z,\ z>5]
   =
   -zB\cdot\frac{459}{40960}\pm c_1z.
\]
Once the strict margins dominate the block error, the Foster--Lyapunov proof of
Theorem~\ref{Theorem~1} applies exactly as in Appendix~\ref{Proof-of-Theorem-1}.
The bounded-population full-state sets are finite after age saturation, so the
process is positive recurrent to the corridor \([2,5]\), under the
non-extinction/regeneration convention.

\subsubsection*{Validation of Theorem~\ref{thm:comm-cost}}

For the finite-cost converse, the irregular graph gives a pointwise version of
the necessary fork-pressure balance near the lower corridor edge. At \(z=1\),
\[
  q\,\mathcal L_{\mathrm{irr,nt}}^{(1)}(3)
  =
  \frac12\cdot\frac{27}{80}
  =
  \frac{27}{160}
  >
  \frac1{20}
  =
  \Ld.
\]
At \(z=2\),
\[
  q\,\mathcal L_{\mathrm{irr,nt}}^{(2)}(3)
  =
  \frac12\cdot\frac{61}{480}
  =
  \frac{61}{960}
  >
  \frac1{20}
  =
  \Ld.
\]
Thus the irregular graph has enough graph-limited fork opportunity at the first
two population levels to balance trap pressure. Consequently, the
exposure-averaged condition
\[
   q\,\overline{\mathcal L}_{\mathrm{fork}}^{-}
   \ge
   \Ld
\]
is consistent with a viable finite-cost operation on this example. This is the
same balance obtained in Appendix~\ref{proof-theorem2} by summing the population
equation, dividing by cumulative population exposure, and comparing long-run
forks with long-run trap deaths.

\subsubsection*{Validation of Theorem~\ref{thm:overshoot}}

For deletion shocks, since \(Z_{\mathrm{low}}=2\), the only integer population
strictly below the corridor is \(z=1\). Hence
\[
   \Gamma_{+,\mathrm{irr}}^{\max}
   =
   q\,\mathcal L_{\mathrm{irr,nt}}^{(1)}(3)
   =
   \frac{27}{160}.
\]
The inner band endpoints are
\[
   Z_\rho^-=(1-\rho)2+\rho5=2+3\rho,
   \qquad
   Z_\rho^+=\rho2+(1-\rho)5=5-3\rho.
\]
Therefore, after a deletion shock to \(z_-<2\),
\[
   \tau_{\mathrm{rec}}^{-}-t_0
   \ge
   \frac{
   \log\!\left((2+3\rho)/z_-\right)
   }{
   \log(1+27/160)
   }.
\]

For insertion shocks, integer levels above the corridor satisfy \(z\ge6\). Since
\(\mathcal L_{\mathrm{irr,nt}}^{(z)}(7)\) decreases in \(z\), the supremum in
the contraction envelope is attained at \(z=6\). Thus
\[
   \Gamma_{-,\mathrm{irr}}^{\max}
   =
   \frac1{20}
   +
   \frac14\,\mathcal L_{\mathrm{irr,nt}}^{(6)}(7)
   <
   1.
\]
Hence, after an insertion shock to \(z_+>5\),
\[
   \tau_{\mathrm{rec}}^{+}-t_0
   \ge
   \frac{
   \log\!\left(z_+/(5-3\rho)\right)
   }{
   -\log\!\left(1-\Gamma_{-,\mathrm{irr}}^{\max}\right)
   }.
\]
These are exactly the two multiplicative comparisons used in
Appendix~\ref{proof-of-theorem-3}: even under the most optimistic local
response, the mean population cannot grow or contract faster than the
envelope-certified corrective-visit rates allow.

\subsubsection*{Illustration of Corollary~\ref{cor:overshoot-strategy}}
The hysteresis controller on the irregular graph uses
\(A_{\downarrow}=3\) for low-population replenishment and
\(A_{\uparrow}=7\) for high-population suppression. These triggering ages are
selected by evaluating the population-dependent eligibility and drift margins
over candidate integer ages: \(A_{\downarrow}=3\) gives the positive
low-population margins \(19/160\) at \(z=1\) and \(13/960\) at \(z=2\), while
\(A_{\uparrow}=7\) makes
\(\Ld+r\,\mathcal L_{\mathrm{irr,nt}}^{(z)}(7)
-q\,\mathcal L_{\mathrm{irr,nt}}^{(z)}(7)\) sufficiently positive above the
desired upper corridor edge, yielding the feasible corridor \([2,5]\). The pair
\((3,7)\) is a feasible graph-dependent choice and is not claimed to be unique
or globally optimal. The low-population block-corrected margin can be chosen as
any
\[
   0<r_{\downarrow}
   <
   \frac{19}{160}-\varepsilon_B,
\]
provided \(\varepsilon_B<19/160\). Above the corridor, since \(z>5\) implies
\(z\ge6\), the high-population margin can be chosen as any
\[
   0<r_{\uparrow}
   <
   M_{\mathrm{S,irr}}(6;7)-\varepsilon_B.
\]
With these choices, the assumptions of
Corollary~\ref{cor:overshoot-strategy} are satisfied. Therefore the hysteresis
controller is positive recurrent to \([2,5]\), and its mean-response recovery
times obey
\[
   \tau_{\mathrm{rec}}^{-}-t_0
   =
   O\!\left(
   \frac1{r_{\downarrow}}
   \log\frac{2+3\rho}{z_-}
   \right),
\]
and
\[
   \tau_{\mathrm{rec}}^{+}-t_0
   =
   O\!\left(
   \frac1{r_{\uparrow}}
   \log\frac{z_+}{5-3\rho}
   \right).
\]
When \(r_{\downarrow}\) and \(r_{\uparrow}\) are comparable to the corresponding
envelope-certified response rates, these upper bounds match the lower bounds of
Theorem~\ref{thm:overshoot} up to constant factors.

\subsection{Takeaway}

The two examples instantiate the full proof architecture of the paper. Lemma
\ref{lem:stationary} gives the stationary distribution, which determines the
trap pressure \(\Ld\). Lemmas~\ref{lem:ret-tail} and
\ref{lem:population-age-envelope} explain why return-time tails control
age-based eligibility. Theorem~\ref{Theorem~1} then converts the resulting
eligibility quantities into drift margins and a recurrent population corridor.
Theorem~\ref{thm:comm-cost} reads the same quantities as a necessary long-run
fork-pressure balance. Theorem~\ref{thm:overshoot} interprets them as speed
limits on recovery after burst disturbances. Finally,
Corollary~\ref{cor:overshoot-strategy} shows that a two-threshold hysteresis
policy using the same ages achieves recovery on the same logarithmic scale.

On the three-node path, the return-tail calculations are short enough to do by
hand, and the certified corridor is \([2,4]\). On the irregular five-node graph,
the same proof is carried out through the submatrix return-tail formula, and the
certified corridor is \([2,5]\). Thus the examples make concrete the central
message of the paper: graph geometry controls decentralized population
regulation through stationary visit frequencies, return-time tails, and the
resulting age-eligibility margins.

\section{Numerical Illustration}
\label{sec:numerical-validation}

This section illustrates the main theoretical mechanisms developed in
Section~\ref{sec:main-results}. The goal is to show how the graph-dependent
quantities appearing in the theorems predict the behavior of a concrete
self-regulating random-walk system.

All experiments use the same heterogeneous modular graph with two trap nodes.
The graph contains two dense communities connected through sparse bridges and
gateway nodes. It is therefore non-regular: the stationary distribution is
nonuniform, and the return-time behavior depends on the visited node. This makes
the example more informative than a regular graph, since the envelope quantities
must account for both topology and population-dependent age eligibility.

The numerical illustration proceeds in four steps. First, we compute the
stationary distribution, the trap pressure \(\Lambda_{\mathrm{death}}\), the
graph-dependent viability and safety margins, and the certified operating
corridor. This illustrates the drift-corridor mechanism of
Theorem~\ref{Theorem~1}. Second, we vary the fork cap \(q\) and check the
exposure-normalized balance condition of Theorem~\ref{thm:comm-cost}. Third, we
apply deletion and insertion shocks after stabilization and compare the observed
mean-response recovery times with the lower bounds in
Theorem~\ref{thm:overshoot}. Finally, we depict the constructive hysteresis
controller of Corollary~\ref{cor:overshoot-strategy} under repeated shocks and
illustrate the predicted logarithmic recovery scaling.

Throughout the section, all recovery times are mean-response recovery times:
they are defined through threshold crossings of \(\mathbb E[Z_t]\), not through
sample-path hitting times. This convention is the same as in
Theorem~\ref{thm:overshoot} and Corollary~\ref{cor:overshoot-strategy}.

\subsection{Illustration of Theorem~1: Drift-Corridor Mechanism}
\label{subsec:numerical-theorem1}

We first illustrate the drift-corridor mechanism of
Theorem~\ref{Theorem~1} and Remark~\ref{rem:corridor-characterization}. The
simulation follows the model of Section~\ref{sec:model}: tokens move according
to the lazy simple random walk, trap nodes delete visiting tokens with their
prescribed probabilities, and non-trap nodes apply local age-based fork or
termination rules. The purpose is to show how the envelope-induced margins determine the operating corridor and
the direction of the population drift on a non-regular graph.

\begin{figure}[t]
\centering
\includegraphics[width=\linewidth]{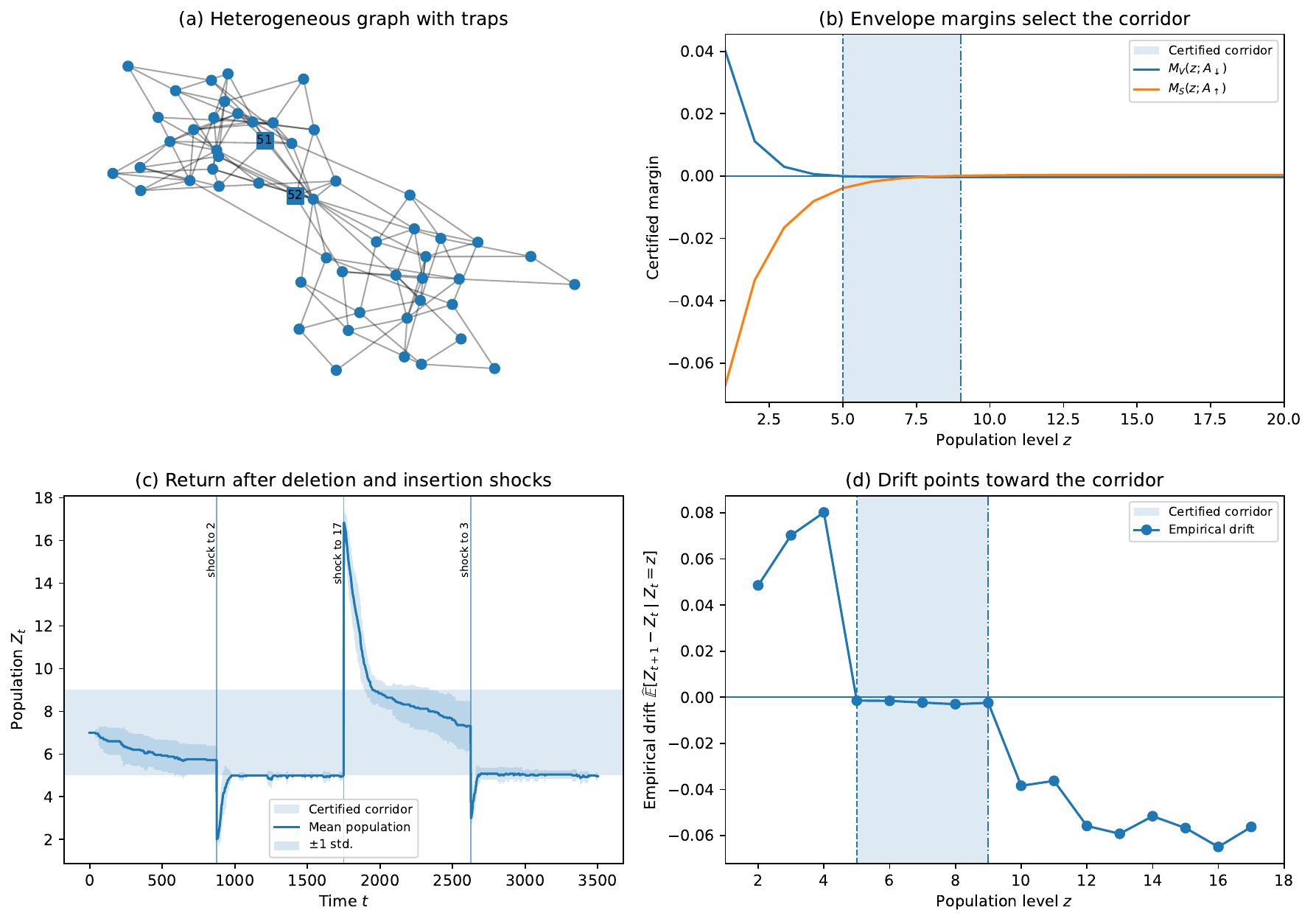}
\caption{Numerical illustration of the drift-corridor mechanism in
Theorem~\ref{Theorem~1}. Panel (a) shows the heterogeneous modular graph with
trap nodes \(51\) and \(52\). Panel (b) shows the computed viability and safety
margins, which select the certified corridor
\([Z_{\mathrm{low}},Z_{\mathrm{high}}]=[5,9]\). Panel (c) shows the mean
population trajectory under deletion and insertion shocks; after each shock,
the population is driven back toward the certified corridor. Panel (d) shows
the empirical conditional drift
\(\widehat{\mathbb E}[Z_{t+1}-Z_t\mid Z_t=z]\), which is positive below the
corridor and negative above it.}
\label{fig:numerical-theorem1-panel}
\end{figure}

The graph used in the experiment is shown in
Fig.~\ref{fig:numerical-theorem1-panel}(a). It is a heterogeneous modular graph
with two dense communities, sparse inter-community bridges, and two trap nodes,
labeled \(51\) and \(52\). Since the graph is non-regular, the stationary
distribution \(\boldsymbol{\pi}\) is nonuniform. The code computes
\(\boldsymbol{\pi}\), the absorption pressure
\(\Lambda_{\mathrm{death}}\)  from \eqref{eq:model-Ldeath}, together with finite-graph return-time tails
computed from the transition matrix restricted to
\(\mathcal V\setminus\{u\}\). More precisely, for each node \(u\), we form the
submatrix obtained by deleting the row and column of \(u\) from the lazy
transition matrix \(\mathbf P'\); powers of this submatrix give the probability
that the walk avoids returning to \(u\) for a prescribed number of steps, as
described in Section~\ref{sec:worked-example}. These return-time quantities are
then used to construct a population-dependent non-trap age-eligibility
surrogate consistent with the envelope definitions in
\eqref{eq:model-L-env} and \eqref{eq:model-eligibility-envelope}.

The controller is a two-threshold hysteresis policy. Below the corridor, a
non-trap node may fork a visiting token if its local age exceeds the
replenishment threshold \(A_{\downarrow}\). Above the corridor, a non-trap node
may deliberately terminate a visiting token if its local age exceeds the
suppression threshold \(A_{\uparrow}\). Inside the corridor, the controller
remains idle. Thus the simulation implements the same local mechanism captured
by the achieved fork and termination rates \(p_{\mathrm{fork}}(z)\) and
\(K_{\mathrm{terminate}}(z)\) in Theorem~\ref{Theorem~1}.

For each population level \(z\), the code evaluates the numerical viability
margin
\[
  M_V(z;A_{\downarrow})
  :=
  q\,\mathcal L_{\pi,z,\mathrm{nt}}(A_{\downarrow})
  -
  \Lambda_{\mathrm{death}},
\]
and the numerical safety margin
\[
  M_S(z;A_{\uparrow})
  :=
  \Lambda_{\mathrm{death}}
  +
  r\,\mathcal L_{\pi,z,\mathrm{nt}}(A_{\uparrow})
  -
  q\,\mathcal L_{\pi,z,\mathrm{nt}}(A_{\uparrow}).
\]
These are the finite-graph counterparts of the envelope margins in
Remark~\ref{rem:corridor-characterization}. The computed margins are shown in
Fig.~\ref{fig:numerical-theorem1-panel}(b).  They select the certified operating corridor
$[Z_{\mathrm{low}},Z_{\mathrm{high}}]=[5,9]. $\footnote{The endpoints are obtained by evaluating the block-corrected
viability and safety inequalities over the integer population levels:
\(Z_{\mathrm{low}}
=
1+\max\{z:
q\,\mathcal L_{\pi,z,\mathrm{nt}}^{+}(A_{\downarrow})
-\Lambda_{\mathrm{death}}
\ge
\varepsilon_B+\eta_-\}\)
and
\(Z_{\mathrm{high}}
=
\min\{z:
q\,\mathcal L_{\pi,y,\mathrm{nt}}^{-}(A_{\uparrow})
-\Lambda_{\mathrm{death}}
-K_{\mathrm{terminate}}(y)
\le
-\varepsilon_B-\eta_+,\;
\forall\,y>z\}\). For the present graph and controller parameters, these
evaluate to \(Z_{\mathrm{low}}=5\) and \(Z_{\mathrm{high}}=9\).}
Below this interval, the viability margin is positive, so the certified
replenishment pressure dominates trap absorption. Above this interval, the
safety margin is positive, so trap losses and deliberate terminations dominate
forking.

The resulting SRRW dynamics are shown in
Fig.~\ref{fig:numerical-theorem1-panel}(c). The system is subjected to three
external shocks: two deletion shocks, which push the population below the
corridor, and one insertion shock, which pushes it above the corridor. After
each shock, the mean population returns toward the certified interval
\([5,9]\). This agrees with the drift interpretation of
Theorem~\ref{Theorem~1}: the population is pushed upward below the corridor and
downward above it.

Finally, Fig.~\ref{fig:numerical-theorem1-panel}(d) reports the empirical
conditional drift
\[
  \widehat{\Delta}(z)
  :=
  \widehat{\mathbb E}
  \big[
  Z_{t+1}-Z_t
  \,\big|\,
  Z_t=z
  \big],
\]
computed from the simulated trajectories after excluding the external shock
times. The empirical drift is positive below \(Z_{\mathrm{low}}\) and negative
above \(Z_{\mathrm{high}}\). This displays the qualitative sign structure described by Theorem~\ref{Theorem~1}. In this example, graph-dependent age
eligibility determines viability and safety margins, these margins identify a
corridor, and the controlled SRRW dynamics are driven back toward this corridor
after disturbances.

\subsection{Illustration of Theorem~\ref{thm:comm-cost}: Finite-Cost Feasibility}
\label{subsec:numerical-theorem2}

We next illustrate the finite-cost feasibility mechanism of
Theorem~\ref{thm:comm-cost}. Theorem~\ref{Theorem~1} explains how a corridor can
be certified through viability and safety margins. Theorem~\ref{thm:comm-cost}
adds a complementary necessary condition: if the system is to remain viable,
safe, and finite-cost, then the long-run graph-limited fork opportunity must be
large enough to compensate the long-run exposure to traps.

The numerical experiment uses the same heterogeneous graph and trap profile as
in Fig.~\ref{fig:numerical-theorem1-panel}. Thus each live token performs one
random-walk transition per time step, and the instantaneous communication cost
is
\[
  C_{\mathrm{comm}}(t)=Z_t.
\]
Consequently, finite communication cost is equivalent to a finite long-run mean
population. The experiment varies the fork cap \(q\) and shows the resulting finite-population behavior.

For each value of \(q\), the simulation estimates the exposure-weighted upper
fork envelope
\[
  \overline{\mathcal L}_{\mathrm{fork}}^{-}
  =
  \frac{
  \sum_t
  \mathbb E\!\left[
  Z_t\,\mathcal L_{\pi,Z_t,\mathrm{nt}}^{-}(A_{\mathrm{eff}}^{\mathrm{ub}})
  \right]
  }{
  \sum_t \mathbb E[Z_t]
  },
\]
which is the finite-horizon analogue of the quantity appearing in
Theorem~\ref{thm:comm-cost}. It then compares the available graph-limited fork
capacity \(q\overline{\mathcal L}_{\mathrm{fork}}^{-}\) with two trap-exposure
quantities. The first is the actual empirical trap exposure per live token,
\[
  \widehat{\Lambda}_{\mathrm{death}}
  =
  \frac{\sum_t \mathbb E[\Phi_{\mathrm{death}}(t)]}
  {\sum_t \mathbb E[Z_t]},
\]
and the second is the stationary-envelope trap pressure
\(\Lambda_{\mathrm{death}}\) from \eqref{eq:model-Ldeath}. These correspond to
the two forms of the converse in Theorem~\ref{thm:comm-cost}: the
actual-exposure form and the stationary-envelope form.

\begin{figure}[t]
\centering
\includegraphics[width=\linewidth]{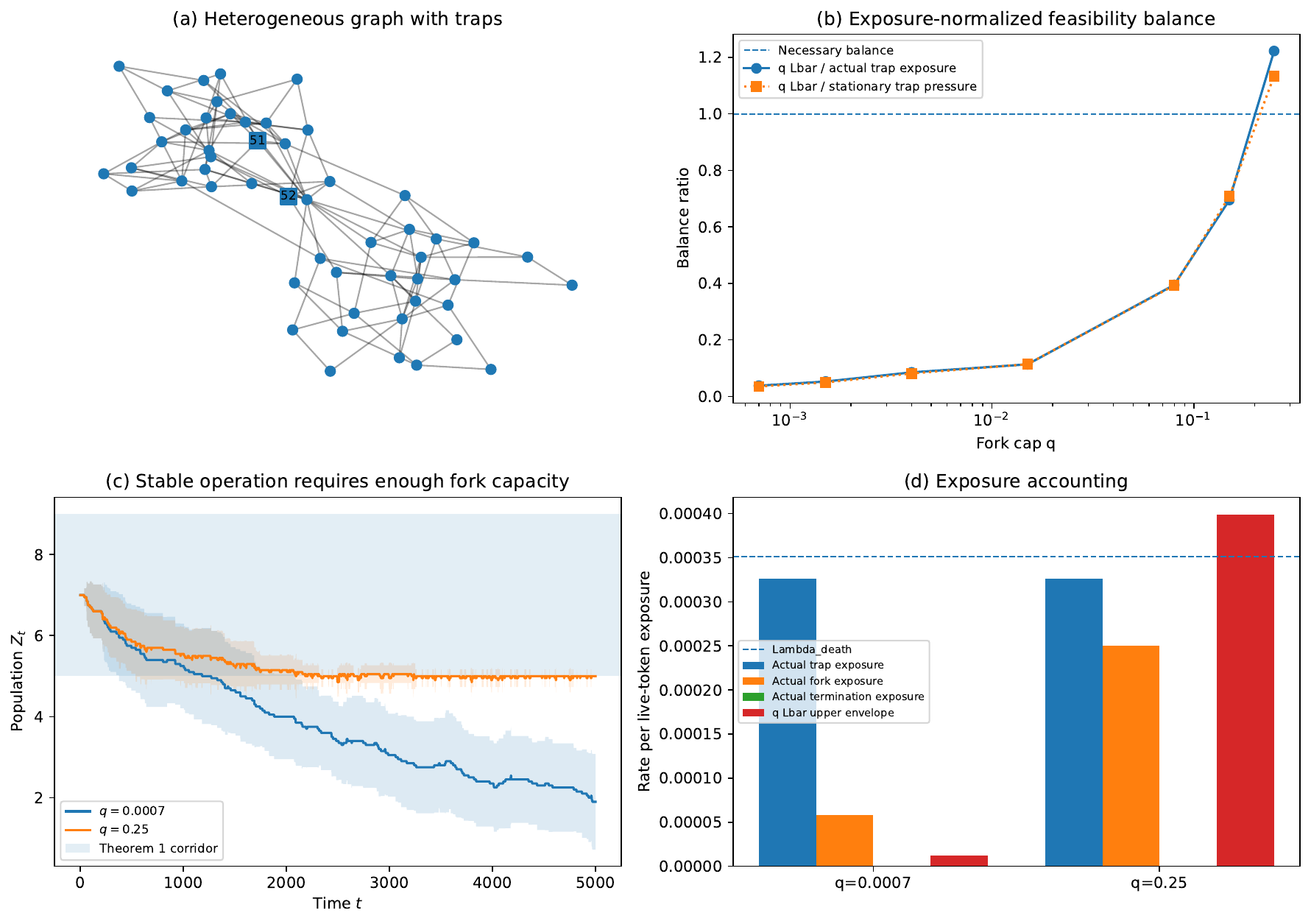}
\caption{Numerical illustration of the finite-cost feasibility converse in
Theorem~\ref{thm:comm-cost}. Panel (a) shows the heterogeneous graph with trap
nodes \(51\) and \(52\). Panel (b) plots the exposure-normalized balance ratios
\(q\overline{\mathcal L}_{\mathrm{fork}}^{-}/\widehat{\Lambda}_{\mathrm{death}}\)
and
\(q\overline{\mathcal L}_{\mathrm{fork}}^{-}/\Lambda_{\mathrm{death}}\)
as the fork cap \(q\) varies; the dashed horizontal line marks the necessary
threshold equal to one. Panel (c) compares representative population trajectories
for an insufficient fork cap \(q=0.0007\) and a feasible fork cap \(q=0.25\).
For the small fork cap, the mean population decays; for the larger fork cap, the
population remains in a finite operating range near the certified corridor.
Panel (d) reports the corresponding exposure accounting: actual trap exposure,
actual fork exposure, deliberate termination exposure, and the envelope upper
bound \(q\overline{\mathcal L}_{\mathrm{fork}}^{-}\), all normalized by cumulative
population exposure.}
\label{fig:numerical-theorem2-panel}
\end{figure}

Figure~\ref{fig:numerical-theorem2-panel}(b) shows the resulting balance ratios
as a function of the fork cap \(q\). The dashed horizontal line marks the
necessary threshold equal to one. For small fork caps, both ratios remain well
below one. In this regime, the graph-limited fork opportunity is too small to
balance trap exposure, and Theorem~\ref{thm:comm-cost} predicts that viable
finite-cost operation cannot be certified. As \(q\) increases, the ratios
increase, and for the largest fork cap shown in the experiment they cross the
necessary feasibility threshold.

The population trajectories in Fig.~\ref{fig:numerical-theorem2-panel}(c)
provide the corresponding dynamical interpretation. For the small fork cap
\(q=0.0007\), the mean population steadily decreases, showing that the system
does not maintain a positive operating level over the simulated horizon. By
contrast, for \(q=0.25\), the mean population remains in a finite operating
range near the certified corridor. This is consistent with
Theorem~\ref{thm:comm-cost}: when the exposure-normalized fork capacity is too
small, finite-cost viable operation fails; when the necessary balance is met,
stable finite-population behavior is no longer ruled out by the converse.

Figure~\ref{fig:numerical-theorem2-panel}(d) gives a direct exposure accounting
for the two representative fork caps. For \(q=0.0007\), the envelope upper bound
\(q\overline{\mathcal L}_{\mathrm{fork}}^{-}\) is far below the observed trap
exposure. This is the infeasible side of the converse. For \(q=0.25\), the
envelope upper bound exceeds the trap-pressure level, so the necessary
finite-cost balance is satisfied. The observed trajectory is then consistent
with finite-population operation.

This figure should be read together with the converse statement. Theorem~\ref{thm:comm-cost} does not say that satisfying
\(q\overline{\mathcal L}_{\mathrm{fork}}^{-}\ge\Lambda_{\mathrm{death}}\) alone
guarantees stability. Stability also requires the viability and safety drift
conditions of Theorem~\ref{Theorem~1}. Rather, the theorem says that if this
exposure balance is violated, then no viable and safe finite-cost policy can be
certified within the stated envelope formulation. The numerical results illustrate the operational message: insufficient graph-limited
fork capacity leads to population decay, while satisfying the exposure balance
is consistent with finite-cost stable operation.

\subsection{Illustration of Theorem~\ref{thm:overshoot}: Reaction Limitation}
\label{subsec:numerical-theorem3}

We now illustrate the reaction-time limitation of
Theorem~\ref{thm:overshoot}. Unlike the corridor illustration in
Subsection~\ref{subsec:numerical-theorem1}, where the emphasis was on the sign
of the drift, the goal here is to depict the speed-limit mechanism: after an
external burst disturbance, the population should not return to the corridor
faster than the graph-limited availability of eligible corrective visits allows.

The experiment uses the same heterogeneous graph and trap profile as in the
previous numerical illustrations. The process is initialized with one live token (\(Z_0=1\)) and evolves
under the non-extinction/regeneration convention of
Section~\ref{sec:model}. After
an initial transient phase, the mean population stabilizes inside the certified
operating corridor \([Z_{\mathrm{low}},Z_{\mathrm{high}}]=[5,9]\). Only after
this stabilization phase do we apply external shocks. This produces the full
sequence
\[
  \text{initialization}
  \;\longrightarrow\;
  \text{stabilization}
  \;\longrightarrow\;
  \text{shock}
  \;\longrightarrow\;
  \text{recovery},
\]
which is the operating regime addressed by Theorem~\ref{thm:overshoot}.

Recovery is measured in the mean-response sense, exactly as in
Theorem~\ref{thm:overshoot}. For a deletion shock that leaves
\(Z_{t_0}=z_-<Z_{\mathrm{low}}\), we define
\[
  \tau_{\mathrm{rec}}^{-}
  :=
  \inf\{t\ge t_0:\mathbb E[Z_t]\ge Z_\rho^{-}\},
\]
where \(Z_\rho^{-}\) is the lower edge of the inner recovery band. For an
insertion shock that leaves \(Z_{t_0}=z_+>Z_{\mathrm{high}}\), we define
\[
  \tau_{\mathrm{rec}}^{+}
  :=
  \inf\{t\ge t_0:\mathbb E[Z_t]\le Z_\rho^{+}\},
\]
where \(Z_\rho^{+}\) is the upper edge of the inner band. Thus we compare
threshold crossings of the mean population, not sample-path hitting times.

Theorem~\ref{thm:overshoot} states that the mean response cannot be faster than
the maximal graph-limited corrective rate. For deletion shocks, the fastest
possible mean growth is bounded by the envelope-certified rate
\(\Gamma_{+}^{\max}\), which gives
\[
  \tau_{\mathrm{rec}}^{-}-t_0
  \ge
  \frac{
  \log\!\left(Z_\rho^{-}/z_-\right)
  }{
  \log(1+\Gamma_{+}^{\max})
  }.
\]
For insertion shocks, the fastest possible mean contraction is bounded by
\(\Gamma_{-}^{\max}\), which gives
\[
  \tau_{\mathrm{rec}}^{+}-t_0
  \ge
  \frac{
  \log\!\left(z_+/Z_\rho^{+}\right)
  }{
  -\log(1-\Gamma_{-}^{\max})
  }.
\]
The numerical experiment evaluates these lower bounds for several deletion
targets \(z_-\in\{1,2,3,4\}\) and insertion targets
\(z_+\in\{11,13,15,17,20\}\).

\begin{figure}[t]
\centering
\includegraphics[width=\linewidth]{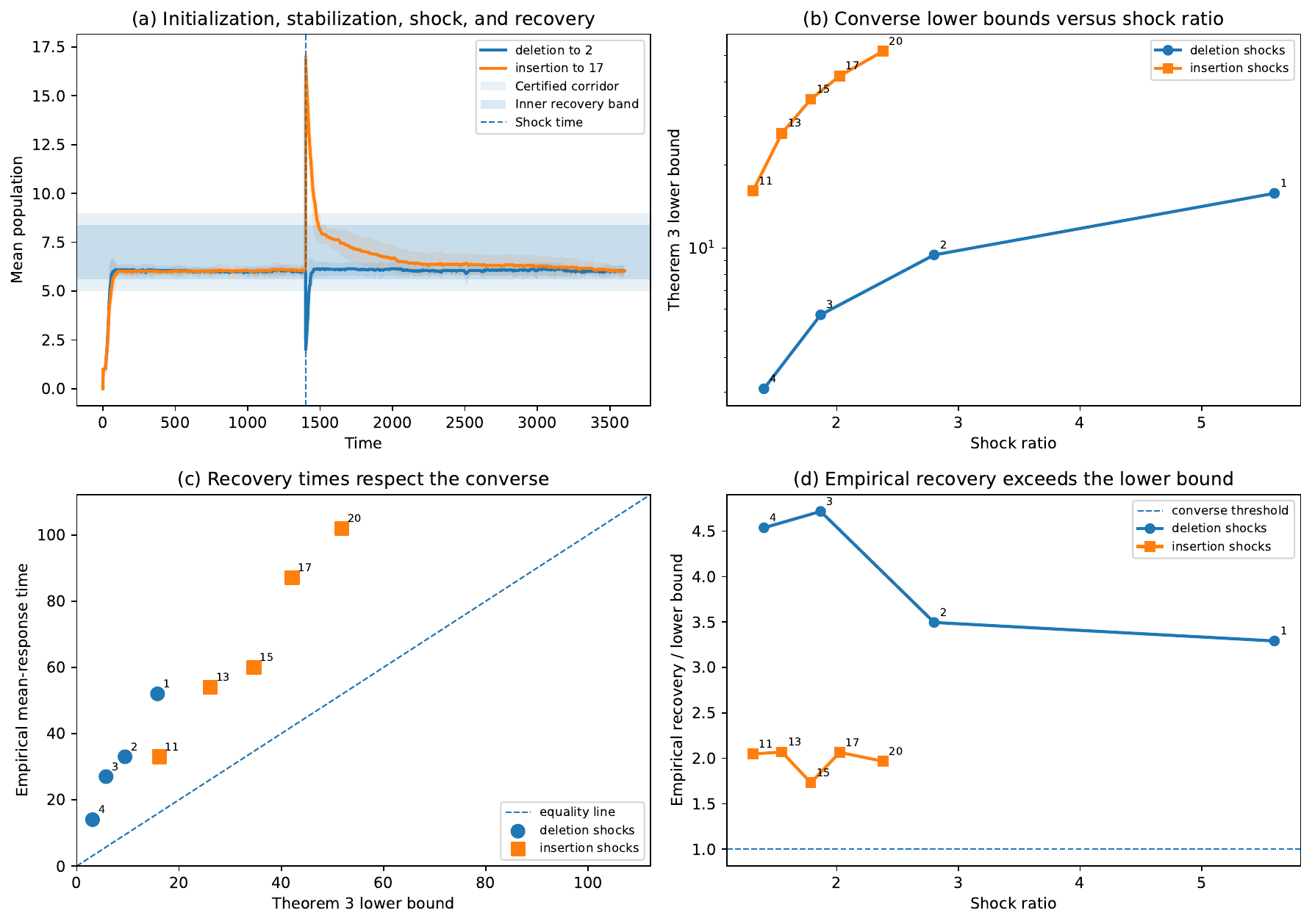}
\caption{Numerical illustration of the reaction-time limitation in
Theorem~\ref{thm:overshoot}. Panel (a) shows the full mean-population evolution:
initialization from one live token, stabilization inside the certified corridor, external
shock, and recovery. The shaded regions show the certified corridor and the
inner recovery band. Panel (b) plots the theorem lower bounds as functions of
the shock ratio for deletion and insertion shocks. Panel (c) compares the
empirical mean-response recovery times with the corresponding lower bounds; all
points lie above the equality line, showing that recovery is not faster than the converse bound. Panel (d) plots the empirical recovery time divided by the
theorem lower bound; all ratios are above one.}
\label{fig:numerical-theorem3-panel}
\end{figure}

The results are shown in Fig.~\ref{fig:numerical-theorem3-panel}. Panel (a)
shows representative mean trajectories. Starting from one live token, the process first
builds up and stabilizes inside the certified corridor. At the shock time, a
deletion shock pushes the mean population below the corridor, while an insertion
shock pushes it above the corridor. In both cases, the controller drives the
mean population back toward the inner recovery band. The insertion recovery is
slower because the shock ratio is larger and because high-population correction
is limited by the available termination and trap-deletion rates.

Panel (b) plots the lower bounds from Theorem~\ref{thm:overshoot} as functions
of the shock ratio. The monotone increase is consistent with the logarithmic
dependence predicted by the theorem: larger disturbances require longer recovery
times, even under the most optimistic corrective rate allowed by the graph.

Panel (c) compares the empirical mean-response recovery times with the
corresponding theorem lower bounds. The dashed line is the equality line. All
points lie above this line, for both deletion and insertion shocks. Hence the
observed recovery is never faster than the lower bound predicted by
Theorem~\ref{thm:overshoot}. This gives a direct numerical depiction of the reaction-time converse.

Panel (d) presents the same comparison as a ratio,
\[
  \frac{
  \text{empirical mean-response recovery time}
  }{
  \text{Theorem~\ref{thm:overshoot} lower bound}
  }.
\]
All ratios remain above one. This shows that the empirical recovery time
exceeds the graph-induced speed limit for every examined shock magnitude. The
deletion-shock ratios are larger because deletion recovery is comparatively
fast in absolute time, while the conservative lower bound remains small. The
insertion-shock ratios are closer to one, indicating that the insertion
experiments more tightly approach the reaction-time limitation.

Overall, Fig.~\ref{fig:numerical-theorem3-panel} illustrates the reaction-time behavior described by Theorem~\ref{thm:overshoot}: after stabilization,
both deletion and insertion shocks recover only after a delay that exceeds the
corresponding theorem lower bound.

\subsection{Illustration of Corollary~\ref{cor:overshoot-strategy}: Hysteresis Recovery}
\label{subsec:numerical-corollary1}

We finally illustrate the constructive hysteresis guarantee in
Corollary~\ref{cor:overshoot-strategy}. The previous subsection illustrated the
reaction-time converse of Theorem~\ref{thm:overshoot}, showing that recovery
cannot be faster than the graph-limited availability of eligible corrective
visits. Corollary~\ref{cor:overshoot-strategy} gives the complementary
achievability statement: a two-threshold hysteresis controller can drive the
population back toward the corridor, with recovery time scaling logarithmically
in the shock magnitude and inversely with the available drift margin.

The controller used in this experiment is the hysteresis controller described
after Theorem~\ref{thm:overshoot}. When the population is below the recovery
band, non-trap nodes apply the replenishment rule: a visiting token may fork if
the local age exceeds the threshold \(A_{\downarrow}\). When the population is
above the recovery band, non-trap nodes apply the suppression rule: a visiting
token may be deliberately terminated if the local age exceeds the threshold
\(A_{\uparrow}\). Inside the band, the controller remains idle. Thus corrective
actions are applied only when the population is outside the desired operating
range.

The numerical experiment uses the same heterogeneous graph and trap profile as
in the previous illustrations. {The process is initialized with one live token (\(Z_0=1\)) and} uses the non-extinction/regeneration convention of Section~\ref{sec:model}.
After this initialization phase, the mean population stabilizes inside the
certified corridor
\[
  [Z_{\mathrm{low}},Z_{\mathrm{high}}]=[5,9].
\]
We then apply repeated external shocks: deletion shocks push the population
below the corridor, and an insertion shock pushes it above the corridor. The
goal is to check whether the hysteresis rule repeatedly returns the population
toward the inner recovery band
\[
  \mathcal C_\rho=[Z_\rho^-,Z_\rho^+].
\]

The corollary predicts two recovery laws. After a deletion shock
\(Z_{t_0}=z_-<Z_{\mathrm{low}}\), the recovery time to the lower edge of the
inner band should scale as
\[
  \tau_{\mathrm{rec}}^{-}-t_0
  =
  O\!\left(
  \frac{1}{r_{\downarrow}}
  \log\frac{Z_\rho^-}{z_-}
  \right),
\]
where \(r_{\downarrow}\) is the positive low-population drift margin. After an
insertion shock \(Z_{t_0}=z_+>Z_{\mathrm{high}}\), the recovery time to the upper
edge of the inner band should scale as
\[
  \tau_{\mathrm{rec}}^{+}-t_0
  =
  O\!\left(
  \frac{1}{r_{\uparrow}}
  \log\frac{z_+}{Z_\rho^+}
  \right),
\]
where \(r_{\uparrow}\) is the positive high-population suppression margin.
Accordingly, the numerical illustration compares empirical mean-response recovery
times with the logarithmic predictors
\[
  \frac{\log(Z_\rho^-/z_-)}{\log(1+r_{\downarrow})}
  \qquad\text{and}\qquad
  \frac{\log(z_+/Z_\rho^+)}{-\log(1-r_{\uparrow})}.
\]
These are finite-sample versions of the two recovery scales in
Corollary~\ref{cor:overshoot-strategy}.

\begin{figure}[t]
\centering
\includegraphics[width=\linewidth]{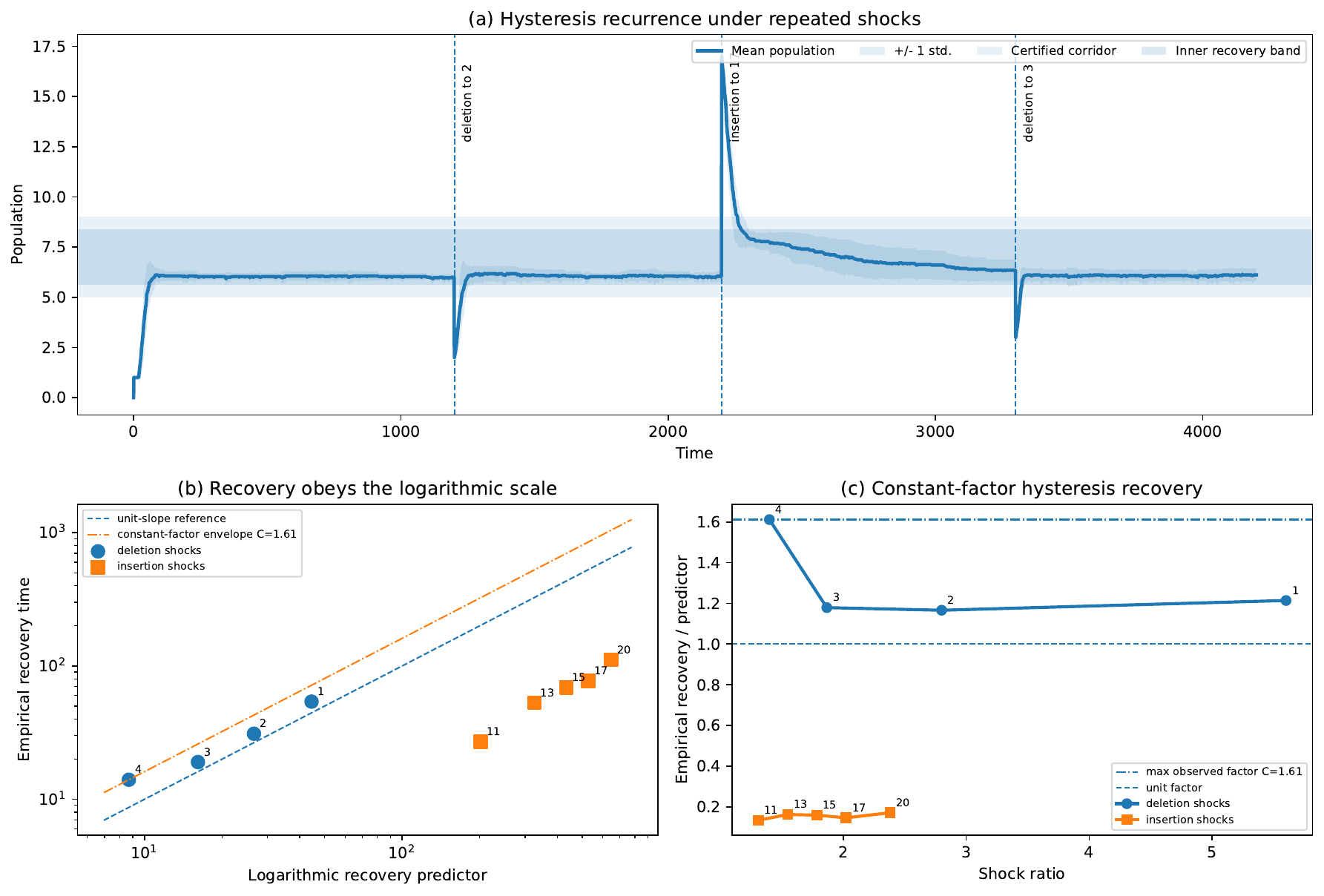}
\caption{Numerical illustration of the hysteresis recovery guarantee in
Corollary~\ref{cor:overshoot-strategy}. Panel (a) shows the full
mean-population evolution under the hysteresis controller: initialization from
one live token, stabilization inside the certified corridor, repeated deletion and
insertion shocks, and recovery toward the inner band. Panel (b) compares the
empirical mean-response recovery times with the logarithmic recovery predictors
appearing in Corollary~\ref{cor:overshoot-strategy}; all points are covered by
a single constant-factor envelope with maximum observed factor \(C=1.61\).
Panel (c) plots the empirical recovery time divided by the corresponding
predictor as a function of the shock ratio. The bounded ratios depict the constant-factor hysteresis recovery behavior predicted by the corollary.}
\label{fig:numerical-corollary1-panel}
\end{figure}

The results are shown in Fig.~\ref{fig:numerical-corollary1-panel}. Panel (a)
shows the full hysteresis behavior under repeated shocks. Starting from one live token, the population first builds up and stabilizes inside the certified corridor.
After the first deletion shock, the mean population drops below the corridor and
is pushed upward by the replenishment rule. After the insertion shock, the mean
population jumps above the corridor and is pushed downward by the suppression
rule. After the final deletion shock, the population again returns toward the
inner recovery band. This repeated return behavior is the recurrence mechanism
claimed in Corollary~\ref{cor:overshoot-strategy}.

Panel (b) depicts the logarithmic recovery scaling. The horizontal axis is the
corresponding logarithmic recovery predictor, while the vertical axis is the
empirical mean-response recovery time. The dashed line gives a unit-slope
reference, and the dash-dotted line gives the smallest constant-factor envelope
covering all observed points. The empirical recovery times remain below this
single envelope, with maximum observed constant factor approximately
\(C=1.61\). This is consistent with the \(O(\cdot)\) recovery statement of the
corollary: recovery time grows proportionally to the logarithmic shock predictor,
up to a bounded constant factor.

Panel (c) presents the same information through the ratio
\[
  \frac{
  \text{empirical mean-response recovery time}
  }{
  \text{logarithmic recovery predictor}
  }.
\]
The ratios remain uniformly bounded over all examined deletion and insertion
shocks. The deletion-shock ratios are larger but still bounded, while the
insertion-shock ratios are smaller because the empirical insertion recovery is
faster than the conservative predictor. As expected for an achievability bound, empirical recovery may be faster than the predictor; the plotted ratios remain within a bounded constant factor of the predicted logarithmic scale.

Overall, Fig.~\ref{fig:numerical-corollary1-panel} illustrates the constructive message of Corollary~\ref{cor:overshoot-strategy}. The same hysteresis mechanism
that keeps the process recurrent to the certified corridor also gives controlled
recovery after shocks. The empirical results show that the mean population
returns to the operating range after repeated deletion and insertion
disturbances, and that the measured recovery times follow the predicted
logarithmic scaling up to a modest constant factor.

\subsection{Summary of Numerical Illustrations}
\label{subsec:numerical-summary}

The numerical experiments illustrate the theoretical picture developed in
Section~\ref{sec:main-results}. On the heterogeneous modular graph, the
envelope-based viability and safety margins identify the certified operating
corridor \([Z_{\mathrm{low}},Z_{\mathrm{high}}]=[5,9]\). The empirical
conditional drift is positive below this corridor and negative above it, which
is exactly the qualitative drift structure required by
Theorem~\ref{Theorem~1}.

The finite-cost experiment depicts the complementary role of
Theorem~\ref{thm:comm-cost}. When the exposure-normalized graph-limited fork
capacity is below the trap exposure, the population cannot maintain a positive
finite operating level. When the balance condition is satisfied, the observed
trajectory is consistent with finite-cost stable operation. This illustrates the converse mechanism: the theorem gives a necessary balance condition rather than a standalone sufficiency result.

The shock experiments illustrate the reaction-time interpretation of
Theorem~\ref{thm:overshoot}. After the process has stabilized near the certified
corridor, deletion and insertion shocks are applied. In all examined cases, the
empirical mean-response recovery times remain above the corresponding theorem
lower bounds. Thus recovery is not faster than the graph-induced availability of
eligible corrective visits permits.

Finally, the hysteresis experiments illustrate the constructive recovery guarantee
of Corollary~\ref{cor:overshoot-strategy}. Under repeated deletion and insertion
shocks, the two-threshold controller repeatedly drives the mean population back
toward the inner recovery band. Moreover, the empirical recovery times remain
within a bounded constant factor of the logarithmic predictors appearing in the
corollary. This depicts the expected
\(O(r_{\downarrow}^{-1}\log(\cdot))\) and
\(O(r_{\uparrow}^{-1}\log(\cdot))\) recovery behavior.

Taken together, the experiments show that the return-time envelopes, drift
margins, finite-cost balance, reaction-time limitation, and hysteresis recovery
law form a coherent design principle for decentralized SRRW population
regulation on heterogeneous graphs.

\section{Conclusion}
\label{sec:conclusion}

This work developed a graph-aware framework for regulating populations of
self-regulating random walks on finite graphs. The main idea is that the
stationary distribution and return-time behavior of the underlying walk
determine both the exposure to trap nodes and the frequency of age-eligible
local control actions. Using population-dependent return-time envelopes, we
derived viability and safety conditions that identify operating corridors with
positive drift below the corridor and negative drift above it.

We also established two fundamental limitations. First, any viable and safe
finite-cost policy must have sufficient graph-limited fork intensity to balance
trap-induced absorption. Second, recovery after burst deletion or insertion
shocks cannot be faster than the graph permits eligible corrective visits. A
two-threshold hysteresis controller was then shown to be positive recurrent to
the prescribed corridor and to achieve recovery on the same logarithmic scale,
up to constant factors.

The worked examples and numerical illustrations show how the stationary
distribution, trap pressure, return-time quantities, drift margins, operating
corridor, and recovery scales can be evaluated on concrete graphs. Possible
extensions include sharper finite-population envelopes and the study of
nonreversible, adaptive, or time-varying mobility graphs.

\subsection*{Acknowledgment}
The authors used OpenAI's ChatGPT to assist with language editing,
presentation, and the organization of portions of the manuscript, as well as
with the implementation and refinement of numerical simulation code. All
theoretical developments, mathematical arguments, numerical results, and
conclusions were independently verified by the authors, who take full
responsibility for the content of the manuscript.
\subsection*{Organization of the paper}
\bibliographystyle{IEEEtran}
\bibliography{ref}

\appendices\label{sec:proofs}
\begin{figure}[t]
\centering
\resizebox{0.96\textwidth}{!}{%
\begin{tikzpicture}[
    font=\scriptsize,
    >=Latex,
    node distance=4.5mm and 12mm,
    box/.style={
        draw,
        rounded corners=3pt,
        thick,
        align=center,
        inner sep=3.2pt,
        minimum height=7mm
    },
    modelbox/.style={
        box, draw=black!70, fill=gray!12,
        text width=4.9cm
    },
    resultbox/.style={
        box, draw=blue!65!black, fill=blue!9,
        text width=4.9cm
    },
    remarkbox/.style={
        box, draw=teal!70!black, fill=teal!9,
        text width=4.9cm
    },
    examplebox/.style={
        box, draw=violet!70!black, fill=violet!10,
        text width=4.9cm
    },
    validationbox/.style={
        box, draw=red!70!black, fill=red!8,
        text width=4.9cm
    },
    proofbox/.style={
        box, draw=orange!80!black, fill=orange!12,
        text width=5.05cm
    },
    lemmabox/.style={
        box, draw=green!45!black, fill=green!9,
        text width=5.05cm
    },
    arrow/.style={->, thick},
    group/.style={
        draw=black!45,
        dashed,
        rounded corners=6pt,
        inner sep=5pt
    }
]

\node[modelbox] (m0) at (0,0)
{System model\\
lazy walk, $\boldsymbol{\pi}$, traps, $\Ld$,\\
local ages, W--AC policies};

\node[resultbox, below=of m0] (m1)
{Theorem~\ref{Theorem~1}\\
viability and safety envelopes};

\node[remarkbox, below=of m1] (m2)
{Remark~\ref{rem:corridor-characterization}\\
margin-certified operating corridor};

\node[resultbox, below=of m2] (m3)
{Theorem~\ref{thm:comm-cost}\\
finite-cost feasibility converse};

\node[resultbox, below=of m3] (m4)
{Theorem~\ref{thm:overshoot}\\
reaction-time limitation};

\node[resultbox, below=of m4] (m5)
{Corollary~\ref{cor:overshoot-strategy}\\
hysteresis recovery};

\node[examplebox, below=of m5] (m6)
{Worked examples\\
explicit corridors, drift checks,\\
and hand-checkable calculations};

\node[validationbox, below=of m6] (m7)
{Numerical validation\\
heterogeneous graph simulations,\\
finite-cost balance, shock recovery};

\draw[arrow] (m0) -- (m1);
\draw[arrow] (m1) -- (m2);
\draw[arrow] (m2) -- (m3);
\draw[arrow] (m3) -- (m4);
\draw[arrow] (m4) -- (m5);
\draw[arrow] (m5) -- (m6);
\draw[arrow] (m6) -- (m7);

\node[group,
    fit=(m0)(m1)(m2)(m3)(m4)(m5)(m6)(m7),
    label={[font=\small]above:Main paper}
] {};

\node[modelbox, right=18mm of m0] (a0)
{Appendix structure\\
proofs and supporting estimates};

\node[proofbox, below=of a0] (a1)
{Appendix A\\
Proof of Lemma~\ref{lem:stationary}};

\node[proofbox, below=of a1] (a2)
{Appendix B\\
Proof of Theorem~\ref{Theorem~1}};

\node[lemmabox, below=of a2] (a3)
{Key ingredients in Appendix B\\
return-time tails, population-age envelopes,\\
blockwise averaging};

\node[proofbox, below=of a3] (a4)
{Appendix C\\
Proof of Theorem~\ref{thm:comm-cost}};

\node[proofbox, below=of a4] (a5)
{Appendix D\\
Proof of Theorem~\ref{thm:overshoot}};

\node[proofbox, below=of a5] (a6)
{Appendix E\\
Proof of Corollary~\ref{cor:overshoot-strategy}};

\draw[arrow] (a0) -- (a1);
\draw[arrow] (a1) -- (a2);
\draw[arrow] (a2) -- (a3);
\draw[arrow] (a3) -- (a4);
\draw[arrow] (a4) -- (a5);
\draw[arrow] (a5) -- (a6);

\node[group,
    fit=(a0)(a1)(a2)(a3)(a4)(a5)(a6),
    label={[font=\small]above:Appendices}
] {};

\draw[arrow, dashed] (m1.east) -- (a2.west);
\draw[arrow, dashed] (m3.east) -- (a4.west);
\draw[arrow, dashed] (m4.east) -- (a5.west);
\draw[arrow, dashed] (m5.east) -- (a6.west);

\end{tikzpicture}%
}
\caption{Roadmap of the main argument, numerical validation, and appendix
proofs. The main text starts from the SRRW model, builds graph-dependent
viability and safety envelopes, uses the induced margins to characterize the
operating corridor, derives finite-cost and reaction-time limits, and gives a
hysteresis-based recovery guarantee. The worked examples make the envelope
calculations explicit on small graphs, while the numerical validation tests the
same drift-corridor, finite-cost, reaction-time, and hysteresis-recovery
mechanisms on a heterogeneous modular graph. The appendices provide the
corresponding proofs and the return-time, population-age, and block-averaging
estimates used in the recurrence arguments.}
\label{fig:appendix-roadmap-clean}
\end{figure}
\clearpage

\section{Proof of Lemma~1}\label{Proof~of~Lemma~1}

\begin{proof}
We verify reversibility with respect to $\boldsymbol{\pi}$, i.e.,
\begin{equation}\label{eq:detailed-balance}
  \boldsymbol{\pi}(u)\mathbf{P}_{uv}
  =
  \boldsymbol{\pi}(v)\mathbf{P}_{vu},
  \qquad\forall\,u,v\in\mathcal{V}.
\end{equation}
For the simple random walk, if $(u,v)\in\mathcal{E}$ then
\[
  \boldsymbol{\pi}(u)\mathbf{P}_{uv}
  =
  \frac{\deg(u)}{2|\mathcal{E}|}\cdot\frac{1}{\deg(u)}
  =
  \frac{1}{2|\mathcal{E}|}
  =
  \frac{\deg(v)}{2|\mathcal{E}|}\cdot\frac{1}{\deg(v)}
  =
  \boldsymbol{\pi}(v)\mathbf{P}_{vu},
\]
and if $(u,v)\notin\mathcal{E}$ both sides are zero. For the weighted
reversible walk, symmetry $w_{uv}=w_{vu}$ gives
\[
  \boldsymbol{\pi}(u)\mathbf{P}_{uv}
  =
  \frac{w_u}{\sum_x w_x}\cdot\frac{w_{uv}}{w_u}
  =
  \frac{w_{uv}}{\sum_x w_x}
  =
  \frac{w_v}{\sum_x w_x}\cdot\frac{w_{vu}}{w_v}
  =
  \boldsymbol{\pi}(v)\mathbf{P}_{vu}.
\]
Hence \eqref{eq:detailed-balance} holds in both cases, establishing
reversibility. Summing \eqref{eq:detailed-balance} over $u$ for a fixed $v$
yields
\[
  \sum_u \boldsymbol{\pi}(u)\mathbf{P}_{uv}
  =
  \boldsymbol{\pi}(v)\sum_u\mathbf{P}_{vu}
  =
  \boldsymbol{\pi}(v),
\]
because $\sum_u\mathbf{P}_{vu}=1$. Therefore
$\boldsymbol{\pi}^{\top}\mathbf{P}=\boldsymbol{\pi}^{\top}$. Since
$\mathcal{G}$ is connected, the chain is irreducible.

The non-lazy walk may be periodic, for example on a bipartite graph. The lazy
modification
\[
  \mathbf{P}'=\epsilon\mathbf{I}+(1-\epsilon)\mathbf{P},
  \qquad 0<\epsilon<1,
\]
is aperiodic and preserves the same stationary distribution because
\[
  \boldsymbol{\pi}^{\top}\mathbf{P}'
  =
  \epsilon\boldsymbol{\pi}^{\top}
  +(1-\epsilon)\boldsymbol{\pi}^{\top}\mathbf{P}
  =
  \boldsymbol{\pi}^{\top}.
\]
Thus $\mathbf{P}'$ is irreducible and aperiodic, and it admits the unique
stationary distribution $\boldsymbol{\pi}$. Standard finite-state Markov-chain
theory then gives
\[
  \alpha^\top(\mathbf{P}')^t\to\boldsymbol{\pi}^{\top}
\]
for every initial distribution $\alpha$.

Finally, conditioned on a fixed population size $Z_t$, if the $Z_t$ token
locations are independent and each is distributed according to
$\boldsymbol{\pi}$, then each token location is an independent categorical draw
with cell probabilities $\boldsymbol{\pi}$. Hence the occupancy vector
$(N_u(t))_{u\in\mathcal{V}}$ conditioned on $Z_t$ follows
\[
  (N_u(t))_{u\in\mathcal{V}}
  \sim
  \mathrm{Multinomial}(Z_t,\boldsymbol{\pi}).
\]
\end{proof}
 
\section{Proof of Theorem~1}
\label{Proof-of-Theorem-1}

We first record the two return-time estimates and the block-averaging estimate
used in the proof. The theorem itself is then proved item by item, following
the order of the statement.

Let $\mathbf{P}'$ be the lazy reversible random-walk kernel on
$\mathcal{G}=(\mathcal{V},\mathcal{E})$, with stationary distribution
$\boldsymbol{\pi}$. The system starts from one live token,
\[
  Z_0=1,
\]
placed at an arbitrary node, or drawn from an arbitrary initial distribution.
The initial condition affects only the transient phase. All drift and envelope
quantities below are evaluated after a burn-in period exceeding the mixing time.
Moreover, as stated in Theorem~\ref{Theorem~1}, the recurrence statement is
understood on the non-extinct communicating class, or under the regeneration
convention that prevents $Z_t=0$ from being absorbing.

Fix $\varepsilon\in(0,1/8]$, and let
$T_{\mathrm{mix}}=T_{\mathrm{mix}}(\varepsilon)$ denote the total-variation
mixing time of the lazy walk:
\begin{align}
\label{mixingtime}
  T_{\mathrm{mix}}(\varepsilon)
  :=
  \min\Big\{
  t\ge0:
  \max_{\alpha}
  \big\|
  \alpha^\top(\mathbf{P}')^t-\boldsymbol{\pi}^{\top}
  \big\|_{\mathrm{TV}}
  \le \varepsilon
  \Big\}.
\end{align}
Thus, after $T_{\mathrm{mix}}$ steps, every token location is within
$\varepsilon$ total-variation distance of $\boldsymbol{\pi}$, uniformly over its
starting node. The resulting stationary-approximation errors are absorbed into
the block-level constants below.

\begin{lemma}[Exponential envelopes for return-time tails]
\label{lem:ret-tail}
On a finite, lazy, reversible Markov chain, for every node
$u\in\mathcal{V}$ there exist constants $c_-(u),c_+(u)>0$ such that, for all
$A\ge1$,
\[
   e^{-c_+(u)A\boldsymbol{\pi}(u)}
   \le
   \Pr_u\{T_u^+\ge A\}
   \le
   e^{-c_-(u)A\boldsymbol{\pi}(u)}.
\]
\end{lemma}

\begin{proof}
Since the chain is finite, lazy, and irreducible, it satisfies a uniform
minorization: there exist $t_0\ge1$ and $\varepsilon_0\in(0,1)$ such that
\begin{equation}
\label{eq:minorization}
  (\mathbf{P}')^{t_0}(x,y)
  \ge
  \varepsilon_0\boldsymbol{\pi}(y),
  \qquad x,y\in\mathcal{V}.
\end{equation}
Fix $u\in\mathcal{V}$. Starting from any state, the probability of hitting $u$
within the next $t_0$ steps is at least
$\varepsilon_0\boldsymbol{\pi}(u)$. By the strong Markov property, partitioning
time into blocks of length $t_0$ gives
\[
  \Pr_u\{T_u^+\ge (m+1)t_0\}
  \le
  \big(1-\varepsilon_0\boldsymbol{\pi}(u)\big)^m
  \le
  e^{-m\varepsilon_0\boldsymbol{\pi}(u)}.
\]
Hence, after adjusting constants to absorb the finite prefix
$1\le A<2t_0$, there exists $c_-(u)>0$ such that
\[
  \Pr_u\{T_u^+\ge A\}
  \le
  e^{-c_-(u)A\boldsymbol{\pi}(u)},
  \qquad A\ge1.
\]

For the lower bound, we only need an exponential lower envelope with a
node-dependent constant. If the chain is nontrivial, choose a state
$v\neq u$ with $\mathbf{P}'(u,v)>0$. Laziness gives
$\mathbf{P}'(v,v)\ge \epsilon>0$. Starting from $u$, one admissible way to avoid
returning to $u$ for $A$ steps is to move from $u$ to $v$ at the first step and
then stay at $v$ for the remaining steps. Thus, for all $A\ge2$,
\[
  \Pr_u\{T_u^+\ge A\}
  \ge
  \mathbf{P}'(u,v)\,\epsilon^{A-2}.
\]
Since $\mathbf{P}'(u,v)>0$ is fixed for the node $u$, there exists a finite
constant $C_u<\infty$ such that, after absorbing the finite prefix,
\[
  \Pr_u\{T_u^+\ge A\}
  \ge
  e^{-C_u A},
  \qquad A\ge1.
\]
Because $\boldsymbol{\pi}(u)>0$, we may write
$C_u=c_+(u)\boldsymbol{\pi}(u)$ for some finite $c_+(u)>0$. Hence
\[
  \Pr_u\{T_u^+\ge A\}
  \ge
  e^{-c_+(u)A\boldsymbol{\pi}(u)},
  \qquad A\ge1.
\]
The one-state chain is degenerate for the present non-trap age-control problem
and can be excluded without loss. Combining the two estimates proves the lemma.
\end{proof}

\begin{lemma}[Population-level aggregate-age envelope]
\label{lem:population-age-envelope}
Consider \(z\ge1\) independent copies of the lazy reversible random walk on
\(\mathcal G\), evolving in stationarity with common stationary distribution
\(\boldsymbol{\pi}\). Let \(A_u^{(z)}(t)\) denote the aggregate local age of
node \(u\), namely the time since node \(u\) was last visited by any of the
\(z\) walks strictly before time \(t\). The age is read immediately before the
arrivals at time \(t\) reset it. Then, for every node \(u\in\mathcal V\), there
exist constants \(0<c_-(u)\le c_+(u)<\infty\), depending only on the lazy walk
and on \(u\), such that, for every integer \(A\ge1\) and every integer
\(z\ge1\),
\[
  \exp\{-c_+(u)zA\boldsymbol{\pi}(u)\}
  \le
  \Pr\{A_u^{(z)}(t)\ge A\mid \text{a tagged walk visits }u\text{ at }t\}
  \le
  \exp\{-c_-(u)zA\boldsymbol{\pi}(u)\}.
\]
Consequently,
\[
  \mathcal{L}_{\pi,z,\mathrm{nt}}^{+}(A)
  \le
  \sum_{u\in\mathcal V\setminus\mathcal P_{\mathrm{trap}}}
  \boldsymbol{\pi}(u)
  \Pr\{A_u^{(z)}(t)\ge A\mid \text{a tagged walk visits }u\text{ at }t\}
  \le
  \mathcal{L}_{\pi,z,\mathrm{nt}}^{-}(A),
\]
where
\[
  \mathcal{L}_{\pi,z,\mathrm{nt}}^{\pm}(A)
  :=
  \sum_{u\in\mathcal V\setminus\mathcal P_{\mathrm{trap}}}
  \boldsymbol{\pi}(u)
  \exp\{-c_{\pm}(u)zA\boldsymbol{\pi}(u)\}.
\]
\end{lemma}

\begin{proof}
Fix a node \(u\in\mathcal V\). Work in stationarity and condition on the event
that a tagged walk, say walk \(1\), visits \(u\) at time \(t=0\). The aggregate
age is read immediately before the arrivals at time \(0\) reset the age. Hence
the event \(A_u^{(z)}(0)\ge A\) means that none of the \(z\) walks visited
\(u\) during the preceding age window, strictly before time \(0\). Simultaneous
visits at time \(0\) do not affect the age that is read at that time.

By reversibility, the past trajectory of the tagged walk before time \(0\),
conditioned on \(X_1(0)=u\), has the same distribution as a future trajectory
started from \(u\). Therefore the probability that the tagged walk did not visit
\(u\) during the preceding age window is, up to the harmless endpoint convention
for the window,
\[
  r_u(A):=\Pr_u\{T_u^+\ge A\}.
\]
Changing between \(>\!A\), \(\ge A\), or \(A+1\) only changes constants and is
absorbed into the exponential envelopes.

For each other walk \(i=2,\ldots,z\), stationarity and independence imply that
the probability that walk \(i\) did not visit \(u\) during the same preceding
window is
\[
  h_u(A):=\Pr_{\boldsymbol{\pi}}\{\tau_u>A\},
\]
where
\[
  \tau_u:=\inf\{s\ge1:X_s=u\}
\]
for a walk initialized from \(\boldsymbol{\pi}\). The \(z\) walks are
independent, and only the tagged walk is conditioned to be at \(u\) at time
\(0\). Hence
\[
  \Pr\{A_u^{(z)}(0)\ge A\mid X_1(0)=u\}
  =
  r_u(A)\,h_u(A)^{z-1}.
\]

By Lemma~\ref{lem:ret-tail}, there exist constants
\(a_-(u),a_+(u)>0\) such that, for all \(A\ge1\),
\[
  e^{-a_+(u)A\boldsymbol{\pi}(u)}
  \le
  r_u(A)
  \le
  e^{-a_-(u)A\boldsymbol{\pi}(u)} .
\]
The stationary hitting tail \(h_u(A)\) satisfies bounds of the same exponential
form. The upper bound follows from the same minorization argument as in
Lemma~\ref{lem:ret-tail}: since the chain is finite, lazy, and irreducible,
there exists \(b_-(u)>0\) such that, for all \(A\ge1\),
\[
  h_u(A)
  =
  \Pr_{\boldsymbol{\pi}}\{\tau_u>A\}
  \le
  e^{-b_-(u)A\boldsymbol{\pi}(u)}.
\]
For the lower bound, choose any state \(v\neq u\) with
\(\boldsymbol{\pi}(v)>0\). Such a state exists unless the chain is the trivial
one-state chain, in which case the statement is vacuous for non-trap
age-based control. By laziness, starting from \(v\), the walk can avoid \(u\)
for \(A\) steps by staying at \(v\) throughout the interval, an event with
probability at least \(\epsilon^A\). Therefore
\[
  h_u(A)
  \ge
  \boldsymbol{\pi}(v)\epsilon^A.
\]
Since \(\boldsymbol{\pi}(v)>0\) and \(\boldsymbol{\pi}(u)>0\), there exists
\(b_+(u)<\infty\) such that, after absorbing the finite prefix,
\[
  h_u(A)\ge e^{-b_+(u)A\boldsymbol{\pi}(u)},
  \qquad A\ge1.
\]

Combining these bounds gives
\[
  e^{-a_+(u)A\boldsymbol{\pi}(u)}
  e^{-(z-1)b_+(u)A\boldsymbol{\pi}(u)}
  \le
  r_u(A)h_u(A)^{z-1}
\]
and
\[
  r_u(A)h_u(A)^{z-1}
  \le
  e^{-a_-(u)A\boldsymbol{\pi}(u)}
  e^{-(z-1)b_-(u)A\boldsymbol{\pi}(u)}.
\]
Choose
\[
  c_+(u):=a_+(u)+b_+(u),
  \qquad
  c_-(u):=\min\{a_-(u),b_-(u)\}.
\]
Then, for every \(z\ge1\),
\[
  a_+(u)+(z-1)b_+(u)\le zc_+(u),
\]
and
\[
  a_-(u)+(z-1)b_-(u)\ge zc_-(u).
\]
Therefore, for all \(A\ge1\) and \(z\ge1\),
\[
  e^{-c_+(u)zA\boldsymbol{\pi}(u)}
  \le
  r_u(A)h_u(A)^{z-1}
  \le
  e^{-c_-(u)zA\boldsymbol{\pi}(u)}.
\]
If necessary, replace \(c_+(u)\) by \(\max\{c_+(u),c_-(u)\}\) to ensure
\(c_-(u)\le c_+(u)\); this only weakens the lower bound and preserves the
estimate.

Thus
\[
  e^{-c_+(u)zA\boldsymbol{\pi}(u)}
  \le
  \Pr\{A_u^{(z)}(0)\ge A\mid X_1(0)=u\}
  \le
  e^{-c_-(u)zA\boldsymbol{\pi}(u)}.
\]
By stationarity, the same estimate holds at any time \(t\). Multiplying by
\(\boldsymbol{\pi}(u)\) and summing over
\(u\in\mathcal V\setminus\mathcal P_{\mathrm{trap}}\) gives the claimed
non-trap envelope bounds.
\end{proof}

Lemma~\ref{lem:population-age-envelope} is applied to the frozen-population
block-stationary comparison process; the discrepancy between this comparison
process and the controlled process within a block is absorbed into the
block-averaging error in Lemma~\ref{lem:block-Oz}. Without loss of generality,
we take \(c_+(u)\ge c_-(u)\) for every \(u\in\mathcal V\), by replacing
\(c_+(u)\) with \(\max\{c_+(u),c_-(u)\}\) if necessary. This only decreases the
lower exponential envelope and therefore preserves the lower bound.

For non-trap nodes, define
\[
  \mathcal{L}_{\pi,z,\mathrm{nt}}^{\pm}(A)
  :=
  \sum_{u\in\mathcal{V}\setminus\mathcal{P}_{\mathrm{trap}}}
  \boldsymbol{\pi}(u)
  e^{-c_{\pm}(u)zA\boldsymbol{\pi}(u)}.
\]
Then
\[
  \mathcal{L}_{\pi,z,\mathrm{nt}}^{+}(A)
  \le
  \mathcal{L}_{\pi,z,\mathrm{nt}}^{-}(A),
  \qquad A\ge0,\quad z\ge1.
\]

Let \(R_{\pi,z,\mathrm{nt}}(A)\) denote the true non-trap eligibility
probability in the frozen-population block-stationary comparison process:
\[
  R_{\pi,z,\mathrm{nt}}(A)
  :=
  \sum_{u\in\mathcal V\setminus\mathcal P_{\mathrm{trap}}}
  \boldsymbol{\pi}(u)
  \Pr\{A_u^{(z)}(t)\ge A
  \mid \text{a tagged walk visits }u\text{ at }t\}.
\]
Lemma~\ref{lem:population-age-envelope} gives
\begin{equation}
\label{eq:eligibility-envelope}
  \mathcal{L}_{\pi,z,\mathrm{nt}}^{+}(A)
  \le
  R_{\pi,z,\mathrm{nt}}(A)
  \le
  \mathcal{L}_{\pi,z,\mathrm{nt}}^{-}(A).
\end{equation}

For a general W--AC policy, let \(A_{\mathrm{eff}}^{\mathrm{ub}}\) be the upper
effective age in the regime under consideration. Since the per-visit fork
probability is capped by \(q\), the achieved fork intensity always satisfies
the universal upper bound
\begin{equation}
\label{eq:pfork-upper-envelope}
  p_{\mathrm{fork}}(z)
  \le
  q\,R_{\pi,z,\mathrm{nt}}(A_{\mathrm{eff}}^{\mathrm{ub}})
  \le
  q\,\mathcal{L}_{\pi,z,\mathrm{nt}}^{-}
  (A_{\mathrm{eff}}^{\mathrm{ub}}).
\end{equation}
A matching lower bound is not universal for all W--AC policies, because a policy
may choose a fork probability below the cap on eligible visits. When the policy
realizes the corresponding lower-envelope fork intensity, let
\(A_{\mathrm{eff}}^{\mathrm{lb}}\) be the lower-bound effective age. Then
\begin{equation}
\label{eq:pfork-lower-envelope}
  p_{\mathrm{fork}}(z)
  \ge
  q\,R_{\pi,z,\mathrm{nt}}(A_{\mathrm{eff}}^{\mathrm{lb}})
  \ge
  q\,\mathcal{L}_{\pi,z,\mathrm{nt}}^{+}
  (A_{\mathrm{eff}}^{\mathrm{lb}}).
\end{equation}

For the block construction, define
\[
  A_{\mathrm H}
  :=
  \max\{1,A_{\mathrm{eff}}^{\mathrm{ub}},A_{\mathrm{eff}}^{\mathrm{lb}}\},
\]
with the convention that, if a lower-bound effective age is not invoked, one may set
\(A_{\mathrm{eff}}^{\mathrm{lb}}=A_{\mathrm{eff}}^{\mathrm{ub}}\). Choose a
block length
\[
  B:=T_{\mathrm{mix}}+\kappa A_{\mathrm H},
  \qquad \kappa\ge4,
\]
and define block times \(t_k:=kB\) and block populations \(Z_k:=Z_{t_k}\). Let
\[
  \Ld
  :=
  \sum_{u\in\mathcal{P}_{\mathrm{trap}}}
  \zeta(u)\boldsymbol{\pi}(u)
\]
be the absorption pressure.

We use the one-action-per-visit discipline: at a non-trap visit, at most one of
\textsf{fork}, \textsf{terminate}, or \textsf{pass} occurs. Trap-induced
deletions are accounted for separately.

Strictly speaking, the population process \(Z_k\) alone need not be Markov,
because future dynamics also depend on token locations, saturated node ages,
and local controller states. Let \(Y_{t_k}\) denote the full system state at the
beginning of block \(k\), and let \(Z(Y_{t_k})=z\) be its population component.
In the estimates below, conditioning on \(Z_k=z\) is shorthand for conditioning
on a full block state with population component \(z\); all constants are
uniform over such full states.

\begin{lemma}[Block averaging]
\label{lem:block-Oz}
Fix a block \(k\) and condition on the full block state \(Y_{t_k}=y\), whose
population component is \(Z(y)=z\). Let \(S_{\mathrm{death}}\),
\(S_{\mathrm{fork}}\), and \(S_{\mathrm{terminate}}\) denote respectively the
total numbers of trap-induced deaths, intentional forks, and intentional
terminations during the block. Then there exists a finite constant
\(c_0(B)<\infty\), independent of \(z\) and uniform over all full states \(y\)
with population component \(z\), such that
\[
  \Big|
  \E[S_{\mathrm{death}}\mid Y_{t_k}=y]-zB\Ld
  \Big|
  \le c_0(B)z,
\]
\[
  \Big|
  \E[S_{\mathrm{fork}}\mid Y_{t_k}=y]-zB\,p_{\mathrm{fork}}(z)
  \Big|
  \le c_0(B)z,
\]
and
\[
  \Big|
  \E[S_{\mathrm{terminate}}\mid Y_{t_k}=y]-zB\,K_{\mathrm{terminate}}(z)
  \Big|
  \le c_0(B)z.
\]
Here \(p_{\mathrm{fork}}(z)\) and \(K_{\mathrm{terminate}}(z)\) denote the
corresponding per-visit means in the frozen-population block-stationary
comparison system at population level \(z\).
\end{lemma}

\begin{proof}
During the block, compare the controlled process with a frozen-population
stationary comparison process consisting of \(z\) independent walks initialized
from \(\boldsymbol{\pi}\) and using the same local action rules. In the
comparison process, there are exactly \(zB\) nominal token-visits during the
block, and the per-visit means are, by definition, \(\Ld\),
\(p_{\mathrm{fork}}(z)\), and \(K_{\mathrm{terminate}}(z)\).

The actual process differs from the comparison process because the locations at
the beginning of the block need not be stationary and because population updates
may occur within the block. The first discrepancy is controlled by the
total-variation mixing bound after the initial mixing portion of the block.
The second discrepancy is controlled by the one-action-per-visit rule: over a
fixed block length, each visit can create at most one token or remove at most
one token, so the difference between the actual number of visits and the
frozen-population nominal count is bounded in expectation by a constant times
\(z\), where the constant may depend on the fixed block length \(B\) but not on
\(z\).

For trap deaths, the stationary per-visit function is
\[
  g_{\mathrm{death}}(u)
  :=
  \mathbf{1}\{u\in\mathcal{P}_{\mathrm{trap}}\}\zeta(u),
\]
whose \(\boldsymbol{\pi}\)-mean is
\[
  \E_{\boldsymbol{\pi}}[g_{\mathrm{death}}(U)]
  =
  \sum_{u\in\mathcal{P}_{\mathrm{trap}}}
  \zeta(u)\boldsymbol{\pi}(u)
  =
  \Ld.
\]
The same argument applies to the bounded fork and deliberate-termination
indicators, which are functions of the visited node, the saturated local age,
and the controller's local randomness. Since the ages are saturated for the
recurrence argument, the bounded-population part of the full state space is
finite. Hence the same block-stationary coupling gives the stated \(O(z)\)
errors for all three event counts.
\end{proof}

By the population update,
\[
  Z_{k+1}
  =
  Z_k
  +
  S_{\mathrm{fork}}
  -
  S_{\mathrm{death}}
  -
  S_{\mathrm{terminate}}.
\]
Using Lemma~\ref{lem:block-Oz}, for every full block state \(Y_{t_k}=y\) with
population component \(Z(y)=z\), we obtain
\begin{equation}
\label{eq:block-drift}
  \E[Z_{k+1}-Z_k\mid Y_{t_k}=y]
  =
  zB
  \Big(
  p_{\mathrm{fork}}(z)
  -
  \Ld
  -
  K_{\mathrm{terminate}}(z)
  \Big)
  \pm c_1(B)z,
\end{equation}
for some finite constant \(c_1(B)\) independent of \(z\). Define the normalized
block error
\[
  \varepsilon_B:=\frac{c_1(B)}{B}.
\]
For readability, we keep writing the drift as a function of \(z\), but the
conditioning is always on the full block state.

\begin{proof}[Proof of Theorem~\ref{Theorem~1}]
We prove the six statements in the theorem in order.

\paragraph*{Proof of item \textup{(i)}}
In the low-population regime, deliberate terminations are inactive, so
\(K_{\mathrm{terminate}}(z)=0\). The achieved per-visit drift is therefore
\[
  p_{\mathrm{fork}}(z)-\Ld .
\]
For the population to recover from a low-population state, this drift cannot be
negative. Hence a necessary achieved-rate condition for low-population recovery is
\[
  p_{\mathrm{fork}}(z)-\Ld\ge0,
\]
or equivalently
\[
  p_{\mathrm{fork}}(z)\ge \Ld .
\]
This proves \eqref{eq:thm1-V0}.

\paragraph*{Proof of item \textup{(ii)}}
In the high-population regime, the net per-visit drift is
\[
  d(z)
  =
  p_{\mathrm{fork}}(z)-\Ld-K_{\mathrm{terminate}}(z).
\]
If this quantity is positive, then, at population level \(z\), births dominate
trap losses and deliberate terminations. Thus the population has positive
achieved drift and cannot be certified as safe at that level. Therefore a
necessary achieved-rate condition for high-population safety is
\[
  d(z)\le0,
\]
that is,
\[
  p_{\mathrm{fork}}(z)-\Ld-K_{\mathrm{terminate}}(z)\le0.
\]
This proves \eqref{eq:thm1-S0}.

\paragraph*{Proof of item \textup{(iii)}}
For any W--AC policy, the probability of forking at a non-trap visit is at most
the global cap \(q\). Moreover, by definition of
\(A_{\mathrm{eff}}^{\mathrm{ub}}\), every fork-eligible visit is included in the
event that the local age is at least \(A_{\mathrm{eff}}^{\mathrm{ub}}\).
Therefore
\[
  p_{\mathrm{fork}}(z)
  \le
  q\,R_{\pi,z,\mathrm{nt}}(A_{\mathrm{eff}}^{\mathrm{ub}}).
\]
Using the upper envelope in \eqref{eq:eligibility-envelope}, we get
\[
  p_{\mathrm{fork}}(z)
  \le
  q\,\mathcal L_{\pi,z,\mathrm{nt}}^{-}
  (A_{\mathrm{eff}}^{\mathrm{ub}}).
\]
Consequently, if
\[
  q\,\mathcal L_{\pi,z,\mathrm{nt}}^{-}
  (A_{\mathrm{eff}}^{\mathrm{ub}})
  -
  \Ld
  -
  K_{\mathrm{terminate}}(z)
  \le 0,
\]
then
\[
  p_{\mathrm{fork}}(z)-\Ld-K_{\mathrm{terminate}}(z)\le0.
\]
Thus the achieved-rate safety condition in item \textup{(ii)} is satisfied.
This proves the universal upper-envelope bound and the safety certificate
\eqref{eq:thm1-S-envelope}.

\paragraph*{Proof of item \textup{(iv)}}
A lower-envelope viability certificate is not automatic for every W--AC policy,
because a policy may choose not to fork at the cap on eligible visits. However,
if the policy realizes the lower-envelope fork intensity with lower-bound
effective age \(A_{\mathrm{eff}}^{\mathrm{lb}}\), then
\eqref{eq:pfork-lower-envelope} gives
\[
  p_{\mathrm{fork}}(z)
  \ge
  q\,\mathcal L_{\pi,z,\mathrm{nt}}^{+}
  (A_{\mathrm{eff}}^{\mathrm{lb}}).
\]
Therefore, if
\[
  q\,\mathcal L_{\pi,z,\mathrm{nt}}^{+}
  (A_{\mathrm{eff}}^{\mathrm{lb}})
  \ge
  \Ld,
\]
then
\[
  p_{\mathrm{fork}}(z)\ge\Ld.
\]
Thus the achieved-rate viability condition in item \textup{(i)} is satisfied.
For policies that do not realize the lower envelope, the condition must instead
be checked directly through \(p_{\mathrm{fork}}(z)\ge\Ld\). This proves
\eqref{eq:thm1-V-envelope} and the accompanying qualification.

\paragraph*{Proof of item \textup{(v)}}
First suppose that, in a low-population regime,
\[
  q\,\mathcal L_{\pi,z,\mathrm{nt}}^{-}
  (A_{\mathrm{eff}}^{\mathrm{ub}})
  <
  \Ld.
\]
By the universal upper-envelope bound proved in item \textup{(iii)},
\[
  p_{\mathrm{fork}}(z)
  \le
  q\,\mathcal L_{\pi,z,\mathrm{nt}}^{-}
  (A_{\mathrm{eff}}^{\mathrm{ub}})
  <
  \Ld.
\]
Hence \(p_{\mathrm{fork}}(z)\ge\Ld\) cannot hold. Therefore no W--AC policy with
fork cap \(q\) and that upper-bound effective age can satisfy the low-population
achieved-rate viability condition at that population level.

Second, if \(d(z)>0\) in a high-population regime, then
\[
  p_{\mathrm{fork}}(z)-\Ld-K_{\mathrm{terminate}}(z)>0,
\]
which is the negation of the achieved-rate safety condition in item
\textup{(ii)}. Thus safety fails at that population level.

Finally, suppose the policy is known to realize a lower-envelope fork intensity
and
\[
  q\,\mathcal L_{\pi,z,\mathrm{nt}}^{+}
  (A_{\mathrm{eff}}^{\mathrm{lb}})
  -
  \Ld
  -
  K_{\mathrm{terminate}}(z)
  >0.
\]
Using the lower-envelope bound \eqref{eq:pfork-lower-envelope}, we get
\[
  p_{\mathrm{fork}}(z)-\Ld-K_{\mathrm{terminate}}(z)>0.
\]
Therefore \(d(z)>0\), so the achieved-rate safety condition fails. This proves
all infeasibility statements in item \textup{(v)}.

\paragraph*{Proof of item \textup{(vi)}}
Assume that there exist
\(0<Z_{\mathrm{low}}<Z_{\mathrm{high}}<\infty\) and
\(\eta_-,\eta_+>0\) such that
\[
  p_{\mathrm{fork}}(z)-\Ld
  \ge
  \varepsilon_B+\eta_-,
  \qquad 0<z<Z_{\mathrm{low}},
\]
and
\[
  p_{\mathrm{fork}}(z)-\Ld-K_{\mathrm{terminate}}(z)
  \le
  -\varepsilon_B-\eta_+,
  \qquad z>Z_{\mathrm{high}}.
\]

Below the corridor, deliberate terminations are inactive. Hence
\eqref{eq:block-drift} gives, for \(0<z<Z_{\mathrm{low}}\),
\[
  \E[Z_{k+1}-Z_k\mid Y_{t_k}=y]
  \ge
  zB(\varepsilon_B+\eta_-)-c_1(B)z
  =
  B\eta_- z.
\]
Thus the expected block drift is strictly positive below the corridor:
\begin{equation}
\label{eq:drift-in}
  \E[Z_{k+1}-Z_k\mid Y_{t_k}=y]
  \ge
  B\eta_- z,
  \qquad 0<z<Z_{\mathrm{low}}.
\end{equation}

Above the corridor, \eqref{eq:block-drift} gives, for
\(z>Z_{\mathrm{high}}\),
\[
  \E[Z_{k+1}-Z_k\mid Y_{t_k}=y]
  \le
  -zB(\varepsilon_B+\eta_+)+c_1(B)z
  =
  -B\eta_+z.
\]
Thus the expected block drift is strictly negative above the corridor:
\begin{equation}
\label{eq:drift-out}
  \E[Z_{k+1}-Z_k\mid Y_{t_k}=y]
  \le
  -B\eta_+ z,
  \qquad z>Z_{\mathrm{high}}.
\end{equation}

The envelope implications stated in the theorem follow directly from items
\textup{(iii)} and \textup{(iv)}. Indeed, the lower-envelope inequality
\[
  q\,\mathcal L_{\pi,z,\mathrm{nt}}^{+}
  (A_{\mathrm{eff}}^{\mathrm{lb}})
  -
  \Ld
  \ge
  \varepsilon_B+\eta_-,
  \qquad 0<z<Z_{\mathrm{low}},
\]
implies the achieved low-population margin, whenever the policy realizes the
lower-envelope fork intensity. Similarly, the upper-envelope inequality
\[
  q\,\mathcal L_{\pi,z,\mathrm{nt}}^{-}
  (A_{\mathrm{eff}}^{\mathrm{ub}})
  -
  \Ld
  -
  K_{\mathrm{terminate}}(z)
  \le
  -\varepsilon_B-\eta_+,
  \qquad z>Z_{\mathrm{high}},
\]
implies the achieved high-population margin for every W--AC policy.

It remains to prove recurrence. Define the population corridor
\[
  \mathcal{C}
  :=
  \{z:\,Z_{\mathrm{low}}\le z\le Z_{\mathrm{high}}\}.
\]
For the Foster--Lyapunov argument, use the population Lyapunov function
\[
  V(y):=Z(y),
\]
where \(Z(y)\) is the population component of the full state \(y\). If
\(Z(y)>Z_{\mathrm{high}}\), then \eqref{eq:drift-out} gives
\[
  \E[V(Y_{t_{k+1}})-V(Y_{t_k})\mid Y_{t_k}=y]
  =
  \E[Z_{k+1}-Z_k\mid Y_{t_k}=y]
  \le
  -B\eta_+ Z(y).
\]
In particular, there exists \(\delta>0\) such that
\begin{equation}
\label{eq:FL-drift}
  \E[V(Y_{t_{k+1}})-V(Y_{t_k})\mid Y_{t_k}=y]
  \le
  -\delta,
  \qquad Z(y)>Z_{\mathrm{high}}.
\end{equation}

The low-population region \(0<Z(y)<Z_{\mathrm{low}}\) is bounded. Therefore it
does not need to be controlled by the unbounded part of the Lyapunov drift.
Instead, it is included in the bounded petite set used in the recurrence
argument. The positive drift in \eqref{eq:drift-in} ensures that this bounded
low-population region cannot form a closed recurrent class disjoint from the
corridor.

It remains to justify the petite-set condition. The theorem assumes finite age
saturation,
\[
  \bar A_u(t)=\min\{A_u(t),A_{\max}\},
\]
with \(A_{\max}<\infty\), and the process is considered on the non-extinct
communicating class or with regeneration at zero. Therefore, on every bounded
population set, the full system state has only finitely many token locations,
finitely many saturated age configurations, and finitely many local controller
states. Hence bounded-population sets are petite for the block-sampled
full-state chain.

In particular, the set of full states whose population component satisfies
\[
  1\le Z(y)\le Z_{\mathrm{high}}
\]
is petite. The drift condition \eqref{eq:FL-drift}, together with this
petite-set property, implies positive recurrence of the block-sampled full
system-state chain to this bounded set; see, e.g.,
\cite[Thm.~13.0.1]{MeynTweedie2009}.

It remains only to pass from this bounded recurrent set to the prescribed
corridor. Let
\[
  \mathcal C_{\mathrm{low}}
  :=
  \{y:\,1\le Z(y)<Z_{\mathrm{low}}\}
\]
be the bounded low-population region, and let
\[
  \mathcal C_{\mathrm{cor}}
  :=
  \{y:\,Z_{\mathrm{low}}\le Z(y)\le Z_{\mathrm{high}}\}
\]
be the corridor in the full state space. By finite age saturation,
\(\mathcal C_{\mathrm{low}}\) contains only finitely many full states. Moreover,
for every full state \(y\in\mathcal C_{\mathrm{low}}\), the low-population
drift estimate \eqref{eq:drift-in} gives
\[
  \E[Z_{k+1}-Z_k\mid Y_{t_k}=y]
  \ge
  B\eta_- Z(y)
  >
  0.
\]
Therefore no closed communicating class of the block-sampled chain can be
contained entirely in \(\mathcal C_{\mathrm{low}}\). Indeed, if such a finite
closed class existed, averaging the last strict positive drift inequality with
respect to its invariant distribution would give a strictly positive stationary
drift of the bounded function \(Z\), which is impossible.

Consequently, every recurrent class inside the bounded petite set
\[
  \{y:\,1\le Z(y)\le Z_{\mathrm{high}}\}
\]
must intersect the corridor \(\mathcal C_{\mathrm{cor}}\). Since the
low-population region is finite, the mean entrance time from
\(\mathcal C_{\mathrm{low}}\) into a recurrent class intersecting the corridor
is finite, and the mean hitting time of the corridor inside that finite
recurrent class is also finite. Hence the mean return time to the corridor is
finite.

Consequently, the population component returns to
\([Z_{\mathrm{low}},Z_{\mathrm{high}}]\) infinitely often with finite mean return
time and spends a nonzero fraction of time inside the corridor in steady state.

Since the original process and the block-sampled process differ only by the
finite block length \(B\), positive recurrence of the block chain transfers to
the original-time process up to the finite multiplicative time factor \(B\).
Thus the full system-state process is positive recurrent to the prescribed
corridor. This proves item \textup{(vi)} and completes the proof of
Theorem~\ref{Theorem~1}.
\end{proof}

\section{Proof of Theorem~\ref{thm:comm-cost}}
\label{proof-theorem2}

\begin{proof}
Fix a finite connected graph \(\mathcal{G}=(\mathcal{V},\mathcal{E})\), a
lazy reversible random walk \(\mathbf{P}'\) with stationary distribution
\(\boldsymbol{\pi}\), a trap set \(\mathcal{P}_{\mathrm{trap}}\) with
termination intensities \(\zeta(\cdot)\), and an arbitrary decentralized
age-based policy \(\{\mu_u\}_{u\in\mathcal{V}}\) satisfying the hypotheses of
Theorem~\ref{thm:comm-cost}. Let \(Z_t\) denote the number of live tokens at
time \(t\), and recall the stationary absorption pressure
\[
  \Ld
  =
  \sum_{u\in \mathcal{P}_{\mathrm{trap}}}
  \zeta(u)\boldsymbol{\pi}(u).
\]
Let \(A_{\mathrm{eff}}^{\mathrm{ub}}\) be the upper-bound effective triggering
age associated with the policy in the operating regime under consideration, and
let
\[
  q
  =
  \sup_{u\in\mathcal{V}\setminus\mathcal{P}_{\mathrm{trap}},\,a\ge0}
  q_u(a)
  \in(0,1]
\]
be the global per-visit fork cap.

Under the one-action-per-visit rule, each live token performs one random-walk
transition per time step. Therefore the communication load at time \(t\) is
exactly
\[
  C_{\mathrm{comm}}(t)=Z_t,
\]
and the long-run average communication cost is the long-run mean population:
\[
  \overline C_{\mathrm{comm}}
  =
  \limsup_{T\to\infty}
  \frac1T\sum_{t=0}^{T-1}\E[Z_t]
  =
  \overline Z .
\]
By hypothesis, the policy has finite nonzero mean population in steady
operation, so
\[
  0<\overline Z<\infty .
\]

We now prove the five statements of the theorem in order.

\paragraph*{Proof of item \textup{(i)}}
The quantity
\[
  \overline{\mathcal L}_{\mathrm{fork}}^{-}
  :=
  \limsup_{T\to\infty}
  \frac{
  \sum_{t=0}^{T-1}
  \E\!\left[
  Z_t\,
  \mathcal L_{\pi,Z_t,\mathrm{nt}}^{-}
  (A_{\mathrm{eff}}^{\mathrm{ub}})
  \right]
  }{
  \sum_{t=0}^{T-1}\E[Z_t]
  }
\]
is the exposure-weighted upper envelope of the graph-limited fork opportunities.
The weighting by \(Z_t\) is necessary because, at time \(t\), there are \(Z_t\)
live tokens, and hence \(Z_t\) possible token visits. Dividing by
\(\sum_{t=0}^{T-1}\E[Z_t]\) normalizes by the cumulative population exposure.

By the universal upper-envelope bound from Theorem~\ref{Theorem~1}, conditional
on the population level \(Z_t\), the achieved per-visit fork probability is
bounded by
\[
  p_{\mathrm{fork}}(Z_t)
  \le
  q\,\mathcal{L}_{\pi,Z_t,\mathrm{nt}}^{-}
  (A_{\mathrm{eff}}^{\mathrm{ub}}).
\]
Consequently, in the same block-stationary envelope sense,
\[
  \E[\Phi_{\mathrm{fork}}(t)]
  \le
  q\,
  \E\!\left[
    Z_t\,
    \mathcal{L}_{\pi,Z_t,\mathrm{nt}}^{-}
    (A_{\mathrm{eff}}^{\mathrm{ub}})
  \right].
\]
Summing over \(t=0,\ldots,T-1\) and dividing by
\[
  H_T:=\sum_{t=0}^{T-1}\E[Z_t]
\]
gives
\begin{equation}
  \limsup_{T\to\infty}
  \frac{
  \sum_{t=0}^{T-1}\E[\Phi_{\mathrm{fork}}(t)]
  }{
  \sum_{t=0}^{T-1}\E[Z_t]
  }
  \le
  q\,
  \overline{\mathcal{L}}_{\mathrm{fork}}^{-}.
  \label{eq:T2-fork-upper-ratio}
\end{equation}
This proves that the long-run fork rate per live-token exposure is bounded by
the exposure-weighted upper fork envelope.

\paragraph*{Proof of item \textup{(ii)}}
The actual long-run trap exposure per live token is defined by
\[
  \widehat{\Lambda}_{\mathrm{death}}
  :=
  \liminf_{T\to\infty}
  \frac{
  \sum_{t=0}^{T-1}\E[\Phi_{\mathrm{death}}(t)]
  }{
  \sum_{t=0}^{T-1}\E[Z_t]
  } .
\]
This is the cumulative expected number of trap-induced deaths divided by the
cumulative population exposure. Since \(0<\overline Z<\infty\) and \(Z_t\ge0\),
we have \(H_T\to\infty\), so this normalization is well defined in the long run.
This proves item \textup{(ii)}.

\paragraph*{Proof of item \textup{(iii)}}
The population evolves according to
\[
  Z_{t+1}-Z_t
  =
  \Phi_{\mathrm{fork}}(t)
  -
  \Phi_{\mathrm{terminate}}(t)
  -
  \Phi_{\mathrm{death}}(t),
\]
where \(\Phi_{\mathrm{fork}}(t)\), \(\Phi_{\mathrm{terminate}}(t)\), and
\(\Phi_{\mathrm{death}}(t)\) denote, respectively, the total numbers of fork
events, deliberate terminations, and trap-induced deaths at time \(t\). Summing
from \(t=0\) to \(T-1\) and taking expectations gives
\[
  \sum_{t=0}^{T-1}\E[\Phi_{\mathrm{fork}}(t)]
  =
  \sum_{t=0}^{T-1}\E[\Phi_{\mathrm{terminate}}(t)]
  +
  \sum_{t=0}^{T-1}\E[\Phi_{\mathrm{death}}(t)]
  +
  \E[Z_T-Z_0].
\]
Since deliberate terminations are nonnegative,
\begin{equation}
  \sum_{t=0}^{T-1}\E[\Phi_{\mathrm{fork}}(t)]
  \ge
  \sum_{t=0}^{T-1}\E[\Phi_{\mathrm{death}}(t)]
  +
  \E[Z_T-Z_0].
  \label{eq:T2-balance-lower}
\end{equation}

Divide \eqref{eq:T2-balance-lower} by
\[
  H_T=\sum_{t=0}^{T-1}\E[Z_t].
\]
Because \(H_T\to\infty\) and \(Z_T\ge0\),
\[
  \frac{\E[Z_T-Z_0]}{H_T}
  \ge
  -\frac{Z_0}{H_T}
  \longrightarrow 0 .
\]
Taking the lower limit therefore yields
\begin{equation}
  \liminf_{T\to\infty}
  \frac{
  \sum_{t=0}^{T-1}\E[\Phi_{\mathrm{fork}}(t)]
  }{
  H_T
  }
  \ge
  \widehat{\Lambda}_{\mathrm{death}} .
  \label{eq:T2-fork-lower-ratio-general}
\end{equation}
Thus every viable finite-cost policy must generate enough forks, per live token
on average, to compensate the actual long-run trap exposure.

Combining \eqref{eq:T2-fork-lower-ratio-general} with the upper fork-rate bound
\eqref{eq:T2-fork-upper-ratio}, we obtain
\[
  \widehat{\Lambda}_{\mathrm{death}}
  \le
  q\,\overline{\mathcal{L}}_{\mathrm{fork}}^{-}.
\]
Equivalently,
\[
  q\,\overline{\mathcal{L}}_{\mathrm{fork}}^{-}
  \ge
  \widehat{\Lambda}_{\mathrm{death}} .
\]
This proves the necessary finite-cost balance in item \textup{(iii)}.

\paragraph*{Proof of item \textup{(iv)}}
It remains to relate the actual exposure
\(\widehat{\Lambda}_{\mathrm{death}}\) to the stationary absorption pressure
\(\Ld\). In the stationary-envelope formulation used throughout
Theorem~\ref{Theorem~1}, the per-token trap-death probability is evaluated
after the mixing burn-in. Equivalently, using the same block-stationary
comparison as in Lemma~\ref{lem:block-Oz}, the long-run trap-death exposure
satisfies
\[
  \widehat{\Lambda}_{\mathrm{death}}=\Ld,
\]
up to the stationary-averaging error already absorbed into the block-level
constants. Substituting this identity into the actual-exposure converse from
item \textup{(iii)} gives
\[
  q\,\overline{\mathcal{L}}_{\mathrm{fork}}^{-}
  \ge
  \Ld .
\]
This is exactly \eqref{eq:commcost-lower}. Hence item \textup{(iv)} follows.

Moreover, if
\[
  q\,\overline{\mathcal{L}}_{\mathrm{fork}}^{-}<\Ld,
\]
then the maximal long-run fork intensity allowed by the age envelope is
strictly smaller than the stationary trap-induced absorption pressure. In that
case, no decentralized age-based policy in this stationary-envelope class can
simultaneously maintain viability, safety, and finite steady-state
communication cost.

\paragraph*{Proof of item \textup{(v)}}
Finally, suppose that the steady-state population is concentrated around a
nominal operating level \(z_\star\), so that
\[
  \overline{\mathcal{L}}_{\mathrm{fork}}^{-}
  =
  \mathcal{L}_{\pi,z_\star,\mathrm{nt}}^{-}
  (A_{\mathrm{eff}}^{\mathrm{ub}})
  +o(1).
\]
Substituting this concentration relation into the stationary-envelope condition
from item \textup{(iv)} gives
\[
  q\,
  \Big(
  \mathcal{L}_{\pi,z_\star,\mathrm{nt}}^{-}
  (A_{\mathrm{eff}}^{\mathrm{ub}})
  +o(1)
  \Big)
  \ge
  \Ld .
\]
Equivalently,
\[
  q\,\mathcal{L}_{\pi,z_\star,\mathrm{nt}}^{-}
  (A_{\mathrm{eff}}^{\mathrm{ub}})
  \ge
  \Ld
  \quad
  \text{up to the corresponding }o(1)\text{ error}.
\]
Thus, at the nominal operating population level, the graph must provide enough
fork-eligible visits to balance the absorption pressure generated by the trap
set.

Since both fork opportunities and trap deaths scale linearly with the number of
live tokens, the population factor cancels in the exposure-normalized balance.
Thus the result is a finite-cost feasibility converse, not a nontrivial lower
bound on
\(\overline C_{\mathrm{comm}}=\overline Z\). This proves item \textup{(v)} and
completes the proof of Theorem~\ref{thm:comm-cost}.
\end{proof}
\section{Proof of Theorem~3}
\label{proof-of-theorem-3}

\begin{proof}
Fix a finite, connected graph $\mathcal{G}=(\mathcal{V},\mathcal{E})$ and
consider the lazy reversible random walk with stationary distribution
$\boldsymbol{\pi}$. Let $\{Z_t\}_{t\ge0}$ be the total population under an
arbitrary local age-based policy satisfying the viability and safety conditions
\textbf{(V)}--\textbf{(S)} of Theorem~\ref{Theorem~1}. Let
$[Z_{\mathrm{low}},Z_{\mathrm{high}}]$ be the corresponding operating corridor,
and let $Z_0\in[Z_{\mathrm{low}},Z_{\mathrm{high}}]$ denote the nominal
steady-state population level.

Recall from Lemma~\ref{lem:ret-tail} that, for every node $u$, the return time
$T_u^+$ satisfies
\[
  \exp\!\big(-c_+(u)A\boldsymbol{\pi}(u)\big)
  \le
  \Pr\{T_u^+\ge A\}
  \le
  \exp\!\big(-c_-(u)A\boldsymbol{\pi}(u)\big),
  \qquad A\ge1.
\]
Define the global return-time exponent
\[
  \lambda_r
  :=
  \min_{u\in\mathcal{V}} c_-(u),
\]
so that
\[
  \Pr\{T_u^+\ge A\}
  \le
  e^{-\lambda_r A\boldsymbol{\pi}(u)},
  \qquad
  \forall u\in\mathcal{V},\ A\ge1.
\]
The quantity $\lambda_r$ gives a graph-dependent scale for the availability of
large local ages.

Let
\[
  \mathcal{C}_\rho=[Z_\rho^{-},Z_\rho^{+}]
\]
be the inner band, where
\[
  Z_\rho^{-}
  =
  (1-\rho)Z_{\mathrm{low}}+\rho Z_{\mathrm{high}},
  \qquad
  Z_\rho^{+}
  =
  \rho Z_{\mathrm{low}}+(1-\rho)Z_{\mathrm{high}}.
\]
Since $0<\rho<1/2$, this interval lies strictly inside
$[Z_{\mathrm{low}},Z_{\mathrm{high}}]$.

We work at the same envelope-limited mean-response level as in the statement
of the theorem. Thus, before the mean trajectory re-enters the inner band, we
compare the population with the fastest response allowed by the visit-triggered,
age-gated, one-action-per-visit dynamics.

On the low-population side, each live token can generate at most one visit per
step, and each visit can create at most one additional token. Hence the mean
multiplicative growth is upper bounded by the maximal envelope-certified fork
intensity
\[
  \Gamma_{+}^{\max}
  :=
  \sup_{1\le z<Z_{\mathrm{low}}}
  q\,\mathcal{L}^{-}_{\pi,z,\mathrm{nt}}(A_\downarrow).
\]
Consequently, during the low-population recovery phase,
\begin{equation}
\label{eq:T3-low-growth-bound}
  \E[Z_{t+1}]
  \le
  (1+\Gamma_{+}^{\max})\E[Z_t].
\end{equation}
This bound is deliberately optimistic: it ignores trap deaths and deliberate
terminations and therefore describes the fastest possible mean replenishment
compatible with the local age-gated fork envelope.

Similarly, on the high-population side, a live token can be removed only through
trap deletion or deliberate termination, and the one-action-per-visit rule
allows at most one removal per visit. By definition of
$\Gamma_{-}^{\max}$,
\[
  \Gamma_{-}^{\max}
  :=
  \sup_{z>Z_{\mathrm{high}}}
  \bigl(
  \Ld+\overline K_{\mathrm{terminate}}(z)
  \bigr),
  \qquad 0\le \Gamma_{-}^{\max}<1,
\]
the mean multiplicative contraction is upper bounded by
$\Gamma_{-}^{\max}$. Thus, during the high-population recovery phase,
\begin{equation}
\label{eq:T3-high-contraction-bound}
  \E[Z_{t+1}]
  \ge
  (1-\Gamma_{-}^{\max})\E[Z_t].
\end{equation}
Again, this is the most favorable contraction comparison: it ignores any forks,
which can only increase the population and therefore slow down suppression.

We now prove the three statements of the theorem in order.

\paragraph*{Proof of item \textup{(i)}}
Suppose that at time \(t=t_0\) an exogenous deletion shock leaves
\[
  Z_{t_0}=z_-<Z_{\mathrm{low}}.
\]
Let
\[
  \mu_t:=\E[Z_t].
\]
Before the mean population reaches the lower edge \(Z_\rho^{-}\) of the inner
band, \eqref{eq:T3-low-growth-bound} gives
\[
  \mu_{t+1}
  \le
  (1+\Gamma_{+}^{\max})\mu_t.
\]
Iterating from time \(t_0\) gives, for every integer \(s\ge0\) before recovery,
\[
  \mu_{t_0+s}
  \le
  z_-(1+\Gamma_{+}^{\max})^s.
\]
Therefore, in order to have
\[
  \mu_{t_0+s}\ge Z_\rho^{-},
\]
it is necessary that
\[
  z_-(1+\Gamma_{+}^{\max})^s
  \ge
  Z_\rho^{-}.
\]
Equivalently,
\[
  s
  \ge
  \frac{
  \log\!\left(Z_\rho^{-}/z_-\right)
  }{
  \log(1+\Gamma_{+}^{\max})
  }.
\]
Thus, with
\[
  \tau_{\mathrm{rec}}^{-}
  :=
  \inf\{t\ge t_0:\E[Z_t]\ge Z_\rho^{-}\},
\]
we obtain
\begin{equation}
\label{eq:T3-deletion-lower}
  \tau_{\mathrm{rec}}^{-}-t_0
  \ge
  \frac{1}{\log(1+\Gamma_{+}^{\max})}
  \log\!\left(\frac{Z_\rho^{-}}{z_-}\right).
\end{equation}
If \(\Gamma_{+}^{\max}=0\), the denominator is interpreted in the usual limiting
sense and the lower bound is infinite unless \(z_-\ge Z_\rho^{-}\), which is
consistent with the fact that no envelope-certified growth is available. This
proves the deletion-shock lower bound.

\paragraph*{Proof of item \textup{(ii)}}
Suppose that at time \(t=t_0\) an exogenous insertion shock leaves
\[
  Z_{t_0}=z_+>Z_{\mathrm{high}}.
\]
Before the mean population reaches the upper edge \(Z_\rho^{+}\) of the inner
band, \eqref{eq:T3-high-contraction-bound} gives
\[
  \mu_{t+1}
  \ge
  (1-\Gamma_{-}^{\max})\mu_t.
\]
Iterating from time \(t_0\) gives, for every integer \(s\ge0\) before recovery,
\[
  \mu_{t_0+s}
  \ge
  z_+(1-\Gamma_{-}^{\max})^s.
\]
Therefore, in order to have
\[
  \mu_{t_0+s}\le Z_\rho^{+},
\]
it is necessary that
\[
  z_+(1-\Gamma_{-}^{\max})^s
  \le
  Z_\rho^{+}.
\]
Since \(0\le\Gamma_{-}^{\max}<1\), this is equivalent to
\[
  s
  \ge
  \frac{
  \log\!\left(z_+/Z_\rho^{+}\right)
  }{
  -\log(1-\Gamma_{-}^{\max})
  }.
\]
Thus, with
\[
  \tau_{\mathrm{rec}}^{+}
  :=
  \inf\{t\ge t_0:\E[Z_t]\le Z_\rho^{+}\},
\]
we obtain
\begin{equation}
\label{eq:T3-insertion-lower}
  \tau_{\mathrm{rec}}^{+}-t_0
  \ge
  \frac{1}{-\log(1-\Gamma_{-}^{\max})}
  \log\!\left(\frac{z_+}{Z_\rho^{+}}\right).
\end{equation}
If \(\Gamma_{-}^{\max}=0\), the denominator is again interpreted in the limiting
sense and the lower bound is infinite unless \(z_+\le Z_\rho^{+}\). This proves
the insertion-shock lower bound.

\paragraph*{Proof of item \textup{(iii)}}
Suppose that
\[
  \Gamma_{+}^{\max}\le C_+\lambda_r,
  \qquad
  \Gamma_{-}^{\max}\le C_-\lambda_r,
  \qquad
  C_-\lambda_r<1,
\]
for constants \(C_+,C_->0\).

For the deletion case, use
\[
  \log(1+x)\le x,
  \qquad x\ge0.
\]
Thus, whenever \(\Gamma_{+}^{\max}>0\),
\[
  \frac1{\log(1+\Gamma_{+}^{\max})}
  \ge
  \frac1{\Gamma_{+}^{\max}}
  \ge
  \frac1{C_+\lambda_r}.
\]
Substituting this into \eqref{eq:T3-deletion-lower} gives
\[
  \tau_{\mathrm{rec}}^{-}-t_0
  =
  \Omega\!\left(
    \frac{1}{\lambda_r}
    \log\frac{Z_\rho^{-}}{z_-}
  \right).
\]

For the insertion case, use
\[
  -\log(1-x)
  \le
  \frac{x}{1-x},
  \qquad 0\le x<1.
\]
With \(x=\Gamma_{-}^{\max}\le C_-\lambda_r<1\), this gives
\[
  \frac1{-\log(1-\Gamma_{-}^{\max})}
  \ge
  \frac{1-\Gamma_{-}^{\max}}{\Gamma_{-}^{\max}}
  \ge
  \frac{1-C_-\lambda_r}{C_-\lambda_r}.
\]
Hence, in the regime where \(C_-\lambda_r\) is bounded away from one,
\[
  \tau_{\mathrm{rec}}^{+}-t_0
  =
  \Omega\!\left(
    \frac{1}{\lambda_r}
    \log\frac{z_+}{Z_\rho^{+}}
  \right).
\]

Combining the deletion and insertion cases gives the claimed mean
reaction-time limitation. The result is a graph-induced speed limit on recovery:
after an exogenous burst creates a population deficit or surplus, the return of
the mean population to the inner corridor is limited by the envelope-certified
availability of local corrective visits. When the response envelopes are of
order \(\lambda_r\), this becomes the natural graph-dependent scale
\(\lambda_r^{-1}\log(\cdot)\). The theorem does not claim that the instantaneous
undershoot or overshoot is bounded independently of the burst size, since the
burst itself may create an arbitrarily large excursion.
\end{proof}

\section{Proof of Corollary~\ref{cor:overshoot-strategy}}
\label{Proof-of-corollary-1}
\begin{proof}
The hysteresis controller satisfies the viability and safety inequalities
\textbf{(V)}--\textbf{(S)} of Theorem~\ref{Theorem~1} with strictly positive
margins. Equivalently, in the terminology of
Remark~\ref{rem:corridor-characterization}, the low-population viability margin
is positive below the corridor, while the high-population safety margin is
positive above the corridor. Thus, below the corridor, the mean fork rate
exceeds the trap-deletion rate, while above the corridor the combined
trap-deletion and deliberate-termination rate exceeds the fork rate.

Because the hysteresis controller uses two possible activation ages, define
\[
  A_{\mathrm{H}}
  :=
  \max\{1,A_{\downarrow},A_{\uparrow}\}.
\]
Let
\[
  B:=T_{\mathrm{mix}}+\kappa A_{\mathrm{H}},
  \qquad \kappa\ge4,
\]
and consider the block times \(t_k=kB\) and the block population
\[
  Z_k:=Z_{t_k}.
\]
As in the proof of Theorem~\ref{Theorem~1}, the population process alone need
not be Markov; conditioning on \(Z_k=z\) below should be understood as
conditioning on the full system state \(Y_{t_k}=y\) at time \(t_k\), whose
population component is \(z\). We suppress this notation when no confusion can
arise.

Lemma~\ref{lem:block-Oz} gives the block drift approximation
\[
  \E[Z_{k+1}-Z_k\mid Y_{t_k}=y]
  =
  zB\big(
  p_{\mathrm{fork}}(z)
  -
  \Ld
  -
  K_{\mathrm{terminate}}(z)
  \big)
  \pm c_1(B) z,
\]
where \(c_1(B)<\infty\) is independent of \(z\). Define the normalized block
error
\[
  \varepsilon_B:=\frac{c_1(B)}{B}.
\]

\paragraph*{Proof of item \textup{(i)}}
By assumption, there exist \(r_{\downarrow},r_{\uparrow}>0\) such that
\[
  p_{\mathrm{fork}}(z)-\Ld
  \ge
  \varepsilon_B+r_{\downarrow},
  \qquad
  z<Z_{\mathrm{low}},
\]
and
\[
  p_{\mathrm{fork}}(z)-\Ld-K_{\mathrm{terminate}}(z)
  \le
  -\varepsilon_B-r_{\uparrow},
  \qquad
  z>Z_{\mathrm{high}}.
\]
These are exactly the block-corrected low-population viability and
high-population safety margins required by Theorem~\ref{Theorem~1}. Therefore,
by Theorem~\ref{Theorem~1}, the full system-state process under the hysteresis
controller is positive recurrent to the corridor
\([Z_{\mathrm{low}},Z_{\mathrm{high}}]\). This proves item~\textup{(i)}.

\paragraph*{Proof of item \textup{(ii)}}
When \(Z_t<Z_{\mathrm{low}}\), the hysteresis controller activates its
low-population fork threshold \(A_{\downarrow}\). This is the aggressive
feasible triggering age used for replenishment below the corridor. By the
assumed positive low-population margin, interpreted at the same block scale as
Theorem~\ref{Theorem~1}, deliberate terminations are inactive and
\[
  p_{\mathrm{fork}}(z)-\Ld
  \ge
  \varepsilon_B+r_{\downarrow},
  \qquad z<Z_{\mathrm{low}}.
\]
Therefore, for \(z<Z_{\mathrm{low}}\),
\[
  \E[Z_{k+1}-Z_k\mid Y_{t_k}=y]
  \ge
  zB(\varepsilon_B+r_{\downarrow})-c_1(B)z
  =
  Br_{\downarrow}z .
\]
Equivalently,
\begin{equation}
  \E[Z_{k+1}\mid Y_{t_k}=y]
  \ge
  (1+Br_{\downarrow})z,
  \qquad z<Z_{\mathrm{low}}.
  \label{eq:cor-low-mult-drift}
\end{equation}

Suppose a deletion at time \(t_0\) leaves the system with population
\[
  z_-:=Z_{t_0}<Z_{\mathrm{low}}.
\]
For simplicity, take \(t_0\) to be a block time; if the shock occurs inside a
block, the bounds below change only by an additive \(O(B)\) term. Let
\(Z_\rho^{-}\) be the lower endpoint of the inner band \(\mathcal{C}_\rho\).

The following estimate is a block-scale mean-response estimate, not a pathwise
hitting-time estimate. Define
\[
  \tau_{\mathrm{rec}}^{-}
  :=
  \inf\{t\ge t_0:\E[Z_t]\ge Z_\rho^{-}\}.
\]
While the block-scale mean-response comparison remains below \(Z_\rho^{-}\) and
uses the low-population response regime, \eqref{eq:cor-low-mult-drift} gives
\[
  \E[Z_k]
  \ge
  z_-(1+Br_{\downarrow})^k .
\]
Hence, if the mean population has not yet reached \(Z_\rho^{-}\) after \(k\)
blocks, it is necessary that
\[
  z_-(1+Br_{\downarrow})^k < Z_\rho^{-}.
\]
Therefore the number of blocks needed for the mean-response comparison to reach
\(Z_\rho^{-}\) is at most
\[
  k_{\downarrow}
  =
  O\!\left(
  \frac{1}{\log(1+Br_{\downarrow})}
  \log\frac{Z_\rho^{-}}{z_-}
  \right).
\]
Returning to the original time scale multiplies this by \(B\), and therefore
\[
  \tau_{\mathrm{rec}}^{-}-t_0
  =
  O\!\left(
  \frac{B}{\log(1+Br_{\downarrow})}
  \log\frac{Z_\rho^{-}}{z_-}
  \right)
  +O(B).
\]
Since \(B\) is fixed once the graph and controller are fixed, this may be
written as
\[
  \tau_{\mathrm{rec}}^{-}-t_0
  =
  O\!\left(
  \frac{1}{r_{\downarrow}}
  \log\frac{Z_\rho^{-}}{z_-}
  \right),
\]
with constants depending only on the fixed block construction. This proves
item~\textup{(ii)}.

\paragraph*{Proof of item \textup{(iii)}}
When \(Z_t>Z_{\mathrm{high}}\), the hysteresis controller activates its
high-population termination/suppression threshold \(A_{\uparrow}\). This is the
conservative feasible triggering age used for suppression above the corridor.
By the assumed positive high-population margin, again interpreted at the block
scale of Theorem~\ref{Theorem~1},
\[
  p_{\mathrm{fork}}(z)-\Ld-K_{\mathrm{terminate}}(z)
  \le
  -\varepsilon_B-r_{\uparrow},
  \qquad z>Z_{\mathrm{high}}.
\]
Thus, for \(z>Z_{\mathrm{high}}\),
\[
  \E[Z_{k+1}-Z_k\mid Y_{t_k}=y]
  \le
  -zB(\varepsilon_B+r_{\uparrow})+c_1(B)z
  =
  -Br_{\uparrow}z .
\]
Set
\[
  \eta_{\uparrow}
  :=
  \min\left\{\frac12,Br_{\uparrow}\right\}
  \in(0,1).
\]
Then, possibly weakening the previous inequality,
\begin{equation}
  \E[Z_{k+1}\mid Y_{t_k}=y]
  \le
  (1-\eta_{\uparrow})z,
  \qquad z>Z_{\mathrm{high}}.
  \label{eq:cor-high-mult-drift}
\end{equation}

Suppose an insertion at time \(t_0\) leaves the system with population
\[
  z_+:=Z_{t_0}>Z_{\mathrm{high}}.
\]
Again, if \(t_0\) is not exactly a block time, this only changes the final
bound by an additive \(O(B)\) term. Let \(Z_\rho^{+}\) be the upper endpoint of
the inner band \(\mathcal{C}_\rho\). Define the block-scale mean recovery time
\[
  \tau_{\mathrm{rec}}^{+}
  :=
  \inf\{t\ge t_0:\E[Z_t]\le Z_\rho^{+}\}.
\]
While the block-scale mean-response comparison remains above \(Z_\rho^{+}\) and
uses the high-population response regime, \eqref{eq:cor-high-mult-drift} gives
\[
  \E[Z_k]
  \le
  z_+(1-\eta_{\uparrow})^k .
\]
Hence, if the mean population has not yet reached \(Z_\rho^{+}\) after \(k\)
blocks, it is necessary that
\[
  z_+(1-\eta_{\uparrow})^k > Z_\rho^{+}.
\]
Therefore the number of blocks needed for the mean-response comparison to reach
\(Z_\rho^{+}\) is at most
\[
  k_{\uparrow}
  =
  O\!\left(
  \frac{1}{-\log(1-\eta_{\uparrow})}
  \log\frac{z_+}{Z_\rho^{+}}
  \right).
\]
Returning to the original time scale gives
\[
  \tau_{\mathrm{rec}}^{+}-t_0
  =
  O\!\left(
  \frac{B}{-\log(1-\eta_{\uparrow})}
  \log\frac{z_+}{Z_\rho^{+}}
  \right)
  +O(B).
\]
Since \(B\) is fixed once the graph and controller are fixed, this becomes
\[
  \tau_{\mathrm{rec}}^{+}-t_0
  =
  O\!\left(
  \frac{1}{r_{\uparrow}}
  \log\frac{z_+}{Z_\rho^{+}}
  \right),
\]
again with constants depending only on the fixed block construction. This
proves item~\textup{(iii)}.

\paragraph*{Proof of item \textup{(iv)}}
Combining items~\textup{(ii)} and~\textup{(iii)}, the hysteresis controller
returns the mean-response trajectory to the inner band in logarithmic time in
the corresponding shock ratio:
\[
  \tau_{\mathrm{rec}}^{-}-t_0
  =
  O\!\left(
  \frac{1}{r_{\downarrow}}
  \log\frac{Z_\rho^{-}}{z_-}
  \right),
  \qquad
  \tau_{\mathrm{rec}}^{+}-t_0
  =
  O\!\left(
  \frac{1}{r_{\uparrow}}
  \log\frac{z_+}{Z_\rho^{+}}
  \right).
\]
In regimes where the available drift margins are comparable to the
envelope-certified response rates in Theorem~\ref{thm:overshoot}, these upper
bounds match the reaction limitations of Theorem~\ref{thm:overshoot} up to
constant factors. In particular, when these quantities are of order
\(\lambda_r\), the recovery time has the natural graph-dependent scale
\(\lambda_r^{-1}\log(\cdot)\). This proves item~\textup{(iv)} and completes the
proof of Corollary~\ref{cor:overshoot-strategy}.
\end{proof}

\end{document}